\documentclass[11pt]{article}

\usepackage{varioref}   
\usepackage[all]{xy}    
\usepackage{graphics}
\usepackage{graphicx}
\usepackage[centertags]{amsmath}
\usepackage{amssymb}
\usepackage{amsthm}

 \newtheorem{definition}{Definition}[section]
 
 \newtheorem{example}[definition]{Example}

\theoremstyle{plain}      

 \newtheorem{proposition}[definition]{Proposition}
 \newtheorem{theorem}[definition]{Theorem}
 \newtheorem{corollary}[definition]{Corollary}
 \newtheorem{lemma}[definition]{Lemma}
 \newtheorem{remark}[definition]{Remark}

\newcommand{\CC}{\mathbb{C}}
\newcommand{\RR}{\mathbb{R}}
\newcommand{\ZZ}{\mathbb{Z}}
\newcommand{\HH}{\mathbb{H}} 

\newcommand{\PP}{\mathbb{P}}
\newcommand{\DD}{\mathbb D}

\newcommand{\Tgq}{\overline{T}_g}
\newcommand{\Tgn}{T_{g,n}}
\newcommand{\Tgnp}{T_{g,n+1}}
\newcommand{\Tgnq}{\overline{T}_{g,n}}
\newcommand{\Tgqa}{\overline{T}_g(\alpha)}
\newcommand{\Tgqas}{\overline{T}_g(\alpha')}
\newcommand{\Tgnpq}{\overline{T}_{g,n+1}}
\newcommand{\Tgd}{\hat T_{g}}
\newcommand{\Tgnd}{\hat T_{g,n}}
\newcommand{\Mgn}{M_{g,n}}
\newcommand{\Mgnq}{\overline{M}_{g,n}}
\newcommand{\Mgq}{\overline{M}_g}
\newcommand{\dMgn}{\partial\Mgn}
\newcommand{\dMg}{\partial M_g}
\newcommand{\Ggn}{\Gamma_{g,n}}
\newcommand{\Ggnp}{\Gamma_{g,n+1}}
\newcommand{\Gga}{\Gamma_g(\alpha)}

\newcommand{\Grq}{\bar\Gamma_\iota}

\newcommand{\Ggas}{\Gamma_g(\alpha')}
\newcommand{\Hga}{H_g(\alpha)}
\newcommand{\UKe}{U_{V,\varepsilon}}
\newcommand{\Cref}{X_{\mbox{\scriptsize ref}}}
\newcommand{\Crefo}{X^0_{\mbox{\scriptsize ref}}}
\newcommand{\Crefos}{X^{0,*}_{\mbox{\scriptsize ref}}}
\newcommand{\Uref}{U_{\mbox{\scriptsize ref}}}
\newcommand{\Oref}{\Omega_{\mbox{\scriptsize ref}}}
\newcommand{\pgn}{\pi_{g,n}}
\newcommand{\Cg}{{\cal C}_g}
\newcommand{\Cgn}{{\cal C}_{g,n}}
\newcommand{\Cgnull}{{\cal C}_{g,0}}
\newcommand{\Cgq}{\overline{\cal C}_g}
\newcommand{\Cgqa}{\overline{\cal C}_g(\alpha)}
\newcommand{\Cgnq}{\overline{\cal C}_{g,n}}
\newcommand{\OplusG}{\Omega^+(G_x)}
\newcommand{\Oplusgn}{\Omega^+_{g,n}}
\newcommand{\Oplusgnd}{\hat\Omega^+_{g,n}}
\newcommand{\Oplusgd}{\hat\Omega^+_g}
\newcommand{\Oplusgda}{\hat\Omega^+_g(\alpha)}

\newcommand{\Oplusx}{\Omega^+(x)}
\newcommand{\Oplusxd}{\hat\Omega^+(x)}
\newcommand{\OminG}{\Omega^-(G_x)}
\newcommand{\Sgq}{\overline{S}_g}
\newcommand{\Sgqa}{\overline{S}_g(\alpha)}
\newcommand{\Sgqas}{\overline{S}_g(\alpha')}
\newcommand{\Sgs}{\tilde S_g}
\newcommand{\Sgd}{\hat S_g}
\newcommand{\Sgqd}{\hat{\overline{ S}}_g}
\newcommand{\Stabi}{\mbox{Stab}(\Delta_\iota)}

\newcommand{\Stabio}{\mbox{Stab}^0(\Delta_\iota)}
\newcommand{\Stabid}{\mbox{Stab}(\Delta)}
\newcommand{\emStabid}{\mbox{\em Stab}(\Delta)}

\newcommand{\slzwei}{\mbox{SL}_2}
\newcommand{\pslzwei}{\mbox{PSL}_2}
\newcommand{\tgn}{T_{g,n}}

\newcommand{\tg}{T_g}

\newcommand{\id}{\mbox{id}}

\newcommand{\re}{\mbox{Re}}
\newcommand{\im}{\mbox{Im}}
\newcommand{\pt}{\mbox{\it pt}}

\newcommand{\qdiff}{Q_X}
\newcommand{\qdiffpunct}{Q_{X_0}}
\newcommand{\qeinh}{\Sigma_X}
\newcommand{\qeinhpunct}{\Sigma_{X_0}}
\newcommand{\Deltaquer}{\overline{\Delta}}
\newcommand{\tgquer}{\overline{T}_g}
\newcommand{\mg}{M_g}

\newcommand{\bpm}{\begin{pmatrix}}
\newcommand{\epm}{\end{pmatrix}}

\newcommand{\cquer}{\overline{C}}
\newcommand{\affplus}{\mbox{Aff}^+}

\newcommand{\sozwei}{\mbox{SO}_2}

\newcommand{\trafo}{f}
\newcommand{\fs}{\footnotesize}
\newcommand{\scs}{\scriptsize}
\newcommand{\sscs}{\scriptscriptstyle}
\newcommand{\qc}{\mbox{Diffeo}^+}

\newcommand{\Mod}{\Gamma_g}
\newcommand{\emaffplus}{\mbox{\em Aff}^+}
\newcommand{\proj}{\mbox{proj}}
\newcommand{\Gammaquer}{\bar{\Gamma}}
\newcommand{\Xstern}{X^*}
\newcommand{\iotaquer}{\bar{\iota}}
\newcommand{\iatilde}{i_{\tilde{A}}}
\newcommand{\patilde}{p_{\tilde{A}}}
\newcommand{\Xref}{X_{\mbox{\fs ref}}}
\newcommand{\Cquer}{\overline{C}}

\newcommand{\emid}{\mbox{\em id}}
\newcommand{\emsozwei}{\mbox{\em SO}_2}
\newcommand{\aut}{\mbox{Aut}}

\newcommand{\grad}{\ensuremath{^\circ}} 

\newcounter{diagramm}
\newcommand{\empslzwei}{\mbox{\em PSL}_2}
\newcommand{\emslzwei}{\mbox{\em SL}_2}
\newcommand{\into}{\hookrightarrow}

\newcommand{\Gammaquerm}{\Gammaquer^*}

\begin{document}

\title{On the boundary of Teichm\"uller disks in Teichm\"uller
  and in Schottky space}

\date{}
\author{{\it Frank Herrlich
 \footnote{
 e-mail: {\sf herrlich@math.uni-karlsruhe.de}, 
    },
    Gabriela Schmith\"usen}
  \footnote{
   e-mail: {\sf schmithuesen@math.uni-karlsruhe.de}}\\[3mm]
     \small
     Mathematisches Institut II, 
       Universit\"at Karlsruhe, 76128 Karlsruhe, Germany}

%

\maketitle

\begin{abstract} 

We study the boundary of  
Teichm\"uller disks in $\tgquer$, a partial compacti\-fication
of Teichm\"uller space, and their image in  
Schottky space.\\
We give a broad introduction to Teichm\"uller disks
and explain the relation between Teichm\"uller curves and 
Veech groups.\\ 
Furthermore, we describe Braungardt's construction
of $\tgquer$ and compare it with the Abikoff augmented Teichm\"uller space.
Following Masur, we give a  description of Strebel rays that makes it easy
to understand their end points on the boundary of $\tgquer$. 
This prepares the description of boundary points that 
a Teichm\"uller disk has, with a particular emphasis to the case that it leads 
to a 
Teichm\"uller curve.\\ 
Further on we turn to Schottky space and describe
two different approaches to obtain a partial compactification. 
We give an overview how 
the boundaries of Schottky space, Teichm\"uller space and moduli
space match together and how the actions of the diverse groups on them are 
linked. Finally we consider the image of Teichm\"uller disks in Schottky
space and show that one can choose the projection from Teichm\"uller space
to Schottky space in such a manner that the image of the Teichm\"uller disk
is a quotient by an infinite group.\\[5mm]
\end{abstract}

\newpage

\tableofcontents   


\section{Introduction} \label{intro} 
One of the original motivations that led to the discovery of Teichm\"uller space\index{Teichm\"uller space}
was to better understand the classification of Riemann surfaces. Riemann
himself already saw that the compact Riemann surfaces of genus $g$ with $n$
marked points on it depend on $3g-3+n$ complex parameters (if this number is
positive). More precisely, there is a complex analytic space $\Mgn$\index{moduli space} whose
points correspond in a natural way to the isomorphism classes of such Riemann
surfaces. $\Mgn$ is even an algebraic variety, but its geometry is not easy to
understand. Most of the basic properties are known today, but many
questions on the finer structure of $\Mgn$ are still open.\footnote{Although
  we consider this general setting in a large part of this paper, we shall
restrict ourselves in this introduction to the case $n=0$ and write, as usual,
$M_g$ instead of $M_{g,0}$ (and later $T_g$ instead of $T_{g,0}$).}\\

Many classification problems become more accessible if the objects are endowed
with an additional structure or {\it marking}\index{marked Riemann surface}. The general strategy is to
first classify the marked objects and then, in a second step, to try to
understand the equivalence relation that forgets the marking. The markings
that Teichm\"uller introduced for a compact Riemann surface $X$ consist of
orientation preserving diffeomorphisms $f:\Cref\to X$ from a reference Riemann
surface $\Cref$ to $X$. Markings $(X,f)$ and $(X',f')$ are considered the same
if $f'\circ f^{-1}$ is homotopic to a biholomorphic map. Thus different
markings of a fixed Riemann surface differ by a homotopy class of
diffeomorphisms of $\Cref$. In other words the {\it mapping class group}\index{mapping class group}\index{Teichm\"uller modular group} (or
{\it Teichm\"uller modular group})
\begin{equation}\label{mapgroup}
\Gamma_g=\mbox{Diffeo}^+(\Cref)/\mbox{Diffeo}^0(\Cref)
\end{equation}
acts on the set $T_g$ of all marked Riemann surfaces of genus $g$, and the
orbit space $T_g/\Gamma_g$ is equal to $M_g$ (here $\mbox{Diffeo}^+(\Cref)$
denotes the group of orientation preserving diffeomorphisms of $\Cref$ and
$\mbox{Diffeo}^0(\Cref)$ the subgroup of those that are homotopic to the
identity).\\

Teichm\"uller discovered that in each homotopy class of diffeomorphisms
between compact Riemann surfaces $X$ and $X'$ there is a unique ``extremal
mapping'', i.\,e.\ a quasiconformal map with minimal dilatation. The logarithm
of this dilatation puts a metric on $T_g$, the {\it Teichm\"uller
  metric}\index{Teichm\"uller metric}. With it $T_g$ is a complete metric space, diffeomorphic to
$\RR^{6g-6}$, and $\Gamma_g$ acts on $T_g$ by isometries. There is
also a structure as complex manifold on $T_g$, for which the elements of $\Gamma_g$ act holomorphically and thus
make the quotient map $T_g\to M_g$ into an analytic map between complex spaces.
\\

That the complex structure on $T_g$ is the ``right one'' for the
classification problem can be seen from the fact that there is a family
$\Cg$ of Riemann surfaces over $T_g$ which in a very precise sense is
universal\index{universal family}\index{Riemann surfaces!universal family of}. This family can be obtained as follows: By the uniformization
theorem, the universal covering of a compact Riemann surface $X$ of genus
$g\ge 2$ is (isomorphic to) the upper half plane $\HH$. Any marking
$f:\Cref\to X$ induces an isomorphism $f_*$ from $\pi_g = \pi_1(\Cref)$,
the fundamental group of the reference surface, to $\pi_1(X)$. We may
obtain a holomorphic action of $\pi_g$ on $T_g\times\HH$ as follows: for
$\gamma\in\Gamma_g$, $x=(X,f)\in T_g$ and $z\in\HH$ put
\[\gamma(x,z) = (x,f_*(\gamma)(z)),\]
where we identify $\pi_1(X)$ with the group of deck transformations of the
universal covering $\HH\to X$.
The quotient $\Cg = (T_g\times\HH)/\pi_g$ is a complex manifold with a natural
projection $p:\Cg\to T_g$; the fibre $p^{-1}(X,f)$ is isomorphic to
$X$. Moreover $p$ is proper and therefore $p:\Cg\to T_g$ is a {\it family of
  Riemann surfaces}. The representation of $\Cg$ as a quotient of a manifold
by an action of $\pi_g$ is called a {\it Teichm\"uller structure}\index{Teichm\"uller structure} on this
family. It follows from results of Bers on the uniformization of families (see e.\,g.\ \cite[Thm.\,XVII]{B}) that this family is {\it universal}, i.\,e.\ every other family of
Riemann surfaces of genus $g$ with a Teichm\"uller structure can be obtained
as a pullback from $p:\Cg\to T_g$. In a more fancy language: $T_g$ is a fine
moduli space\index{moduli space!fine} for Riemann surfaces of genus $g$ with Teichm\"uller
structure.\\

It follows by the same arguments that for any family $\pi:{\cal C}\to S$ of Riemann
surfaces (over some complex space $S$) there is an analytic map
$\mu=\mu_\pi:S\to M_g$, which maps $s\in S$ to the point in $M_g$ that
corresponds to the isomorphism class of the fibre
$\pi^{-1}(s)$. Unfortunately, $\Gamma_g$ does not act freely on $T_g$;
therefore the quotient of $\Cg$ by the action of $\Gamma_g$ does not give a
universal family over $M_g$: the fixed points of elements in $\Gamma_g$
correspond to automorphisms of the Riemann surface, and the fibre over $[X]\in
M_g$ in the family $\Cg/\Gamma_g\to M_g$ is the Riemann surface
$X/\mbox{Aut}(X)$ (whose genus is strictly less than $g$ if $\mbox{Aut}(X)$ is
nontrivial). As a consequence, $M_g$ is not a fine moduli space for Riemann
surfaces, but only a ``coarse''\index{moduli space!coarse} one (see e.\,g.\ \cite[1A]{HM} for a precise
definition of fine and coarse moduli spaces).\\

There are several equivalent ways to define markings of Riemann surfaces and to describe Teichm\"uller space. Instead of classes of diffeomorphisms $f:\Cref\to X$ often conjugacy classes of group isomorphisms $\pi_g\to\pi_1(X)$ are used as markings. For the purpose of this paper the approach to Teich\-m\"uller space via  Teichm\"uller deformations\index{Teichm\"uller deformation} is very well suited; it is developed in Section \ref{deform}. The starting point is the observation that a holomorphic quadratic differential $q$ on a Riemann surface $X$ defines a flat structure $\mu$ on $X^* = X-\{\mbox{zeroes of\ }q\}$. Composing the chart maps of $\mu$ with a certain (real) affine map yields a new point in $T_g$. Any point in $T_g$ is in a unique way such a Teichm\"uller deformation of a given base point $(\Cref,\id)$, cf.~Section \ref {sr}.\\

The main objects of interest in this article are  {\it Teichm\"uller embeddings}\index{Teichm\"uller embedding}, i.\,e.\ holomorphic isometric embeddings $\iota:\HH\to T_g$ (or $\iota:\DD\to T_g$), where $\HH$ (resp.\ $\DD$) is given
the hyperbolic metric and $T_g$ the Teichm\"uller metric, see Definition \ref{emb}. The restriction of
$\iota$ to a hyperbolic geodesic line in $\HH$ (or $\DD$) is then a (real) geodesic line
in $T_g$ in the usual sense. The image $\Delta_\iota$ of such an embedding $\iota$ is called a
{\it Teichm\"uller geodesic}\index{Teichm\"uller geodesic} or {\it Teichm\"uller disk}\index{Teichm\"uller disk} in $T_g$. There are
plenty of Teichm\"uller disks in $T_g$. To see this note first that the tangent space
 to $T_g$ at a point $x=(X,f)\in T_g$ is naturally isomorphic to the vector space
$Q_X=H^0(X,\Omega_X^{\otimes 2})$ of holomorphic quadratic differentials on X (this results from the Bers embedding of $T_g$ as a bounded open subdomain of $Q_X$). We shall explain in Section \ref{defdisks} in three different ways how one can, for a given holomorphic quadratic differential $q$ on a Riemann surface $X$, construct a Teichm\"uller embedding $\iota:\DD\to T_g$ with $\iota(0)=x$ and $\iota'(0)=q$. This shows that for any $x\in T_g$ and any
(complex) tangent vector at $x$ there is a Teichm\"uller disk passing
through $x$ in direction of the given tangent vector.\\

There are several natural and closely related objects attached to a Teich\-m\"uller disk $\Delta_\iota$ (or a Teichm\"uller embedding $\iota:\HH\to T_g$): The first is a discrete subgroup of $\pslzwei(\RR)$ called the {\it (projective) Veech group}\index{Veech group} $\Grq$, cf.~Section \ref{vg}. If $q$ is the quadratic differential on the Riemann surface $X$ by which $\iota$ is induced, $\Grq$ consists of the derivatives of those diffeomorphisms of $X$ that are affine with respect to the flat structure defined by $q$. Veech showed that this subgroup of $\pslzwei(\RR)$ is always discrete (\cite[Prop.~2.7]{V}).\\[1mm]
A second group naturally attached to $\iota$ is the {\it stabilizer}\index{Teichm\"uller disk!stabilizer}
\[\Stabi = \{\varphi\in\Gamma_g:\varphi(\Delta_\iota)=\Delta_\iota\}\]
of $\Delta_\iota$ in the Teichm\"uller modular group.
The pointwise stabilizer
\[\Stabio=\{\varphi\in\Gamma_g:\varphi|_{\Delta_\iota}=\id_{\Delta_\iota}\}\]
is a finite subgroup of $\Stabi$, and $\Stabi/\Stabio$ then is (via $\iota$) a
group of isometries of $\HH$ and thus a subgroup of $\pslzwei(\RR)$. This subgroup coincides with the projective Veech group $\Grq$, see Section \ref{lattice}.\\[1mm]
Given a Teichm\"uller embedding $\iota$ we are also interested in the
image $C_\iota$ of $\Delta_\iota$ in the moduli space $M_g$. 
The map $\Delta_\iota\to C_\iota$ obviously factors through the Riemann surface
$\HH/\Grq$ or rather through its mirror image $V_\iota$, see Section \ref{tc} and in particular \ref{lattice}. The typical case seems to be $\Stabi=\{\id\}$ (although it
is not trivial to give explicit examples). Much attention has been given in
recent years to the other extreme case that $\Grq$ is a lattice\index{lattice in $\pslzwei(\RR)$} in
$\pslzwei(\RR)$. Then $V_\iota$ is of finite hyperbolic volume and hence a
Riemann surface of finite type, or equivalently an algebraic curve. In this
case the induced map $V_\iota\to C_\iota$ is birational (see \cite{EG}), i.\,e.\
$V_\iota$ is the desingularization (or normalization) of $C_\iota$. It follows
from a result of Veech (\cite{V}) that $V_\iota$ (and hence also $C_\iota$)
cannot be projective. If $\Grq$ is a lattice, the affine curve $C_\iota$ is
called a {\it Teichm\"uller curve}\index{Teichm\"uller curve}, cf.~Sect.~\ref{tc}. 
First examples were given by Veech \cite{V}; 
in them, $\Grq$ is a hyperbolic triangle group. Later more examples with 
triangle groups as Veech groups were found, see \cite{HS} for a comprehensive 
overview and \cite{BM} for recent results. Explicit examples for Teichm\"uller curves 
also with non triangle groups as Veech groups can be found 
e.~g. in \cite{mcm}, \cite{C} and \cite{L}. M\"oller has shown
(\cite{Martin}) that every Teichm\"uller curve is, as a subvariety of $M_g$, defined
over a number field. This implies that there are at most countably many
Teichm\"uller curves.\\

A special class of Teichm\"uller curves is obtained by {\it origamis}\index{origami} (or {\it square-tiled surfaces}\index{square-tiled surface}). They arise from finite coverings of an elliptic curve that are ramified over only one point. Given such a covering $p:X\to E$, the quadratic differential $q=(p^*\omega_E)^2$ (where $\omega_E$ is the invariant holomorphic differential on $E$) induces a Teichm\"uller embedding whose Veech group is commensurable to $\slzwei(\ZZ)$, see \cite{GJ}. Lochak proposed in \cite{L} a combinatorial construction for such coverings (which led to the name ``origami''), and Schmith\"usen \cite{Algo} gave a group theoretic characterization of the Veech group. In \cite{Diss}, origamis and their Veech groups are systematically studied and numerous examples are presented. Origamis in genus 2 where $q$ has one zero are classified in \cite{HL}. Using the description of origamis by gluing squares it is not difficult to see that there are, for any $g\ge2$, infinitely many Teichm\"uller curves in $M_g$ that come from origamis. In genus 3 there
is even an explicit example of an origami curve that is intersected by infinitely many others, 
see \cite{WMS}. \\

We want to study boundary points of Teichm\"uller disks and Teichm\"uller
curves; by this we mean, for a Teichm\"uller embedding $\iota:\HH\to T_g$, the
closures of $\Delta_\iota$ and $C_\iota$ in suitable (partial) compactifications
of $T_g$ and $M_g$, respectively. For the moduli space we shall use the
compactification\index{moduli space!compactification of} $\Mgq$ by stable Riemann surfaces. Here a one-dimensional connected
compact complex space $X$ is called a {\it stable Riemann surface}\index{stable Riemann surface} if all
singular points of $X$ are ordinary double points, i.\,e.\ have a neighbourhood
isomorphic to $\{(z,w)\in\CC^2:z\cdot w=0,|z|<1,|w|<1\}$; moreover we require
that every irreducible component $L$ of $X$ that is isomorphic to the
projective line $\hat\CC = \PP^1(\CC)$ intersects $\overline{X-L}$ in at least
three points. It was shown by Deligne and Mumford\index{Deligne-Mumford compactification} (\cite{DM}) that stable
Riemann surfaces are classified by an irreducible compact variety $\Mgq$ that,
like $M_g$, has the quality of a coarse moduli space. In fact, with the
approach of Deligne-Mumford it is possible to classify stable algebraic curves
over an arbitrary ground field: they construct a proper scheme over $\ZZ$ of
which $\Mgq$ is the set of complex-valued points. Some years later, Knudsen
\cite{K} showed that $\Mgq$ is a projective variety.\\

If $\iota:\HH\to T_g$ is a Teichm\"uller embedding such that $C_\iota$ is a
Teichm\"uller curve, the closure $\bar C_\iota$ of $C_\iota$ in $\Mgq$ is
Zariski closed and therefore a projective curve. In particular, $\bar C_\iota -
C_\iota$ consists of finitely many points, called the {\it cusps}\index{Teichm\"uller curve!cusps of} of $C_\iota$. It is very
interesting to know, for a given Teichm\"uller curve $C_\iota$, the number of
cusps and the stable Riemann surfaces that correspond to the cusps. In the
case that $\iota$ is induced by an origami
there is an algorithm that determines (among other information) the precise
number of cusps of $C_\iota$, see \cite{Algo}.\\

The boundary\index{moduli space!boundary of} $\dMg = \Mgq-M_g$ is a divisor, i.\,e.\ a projective subvariety
of (complex) codimension 1. It has irreducible components
$D_0,D_1,\dots,D_{[\frac{g}{2}]}$; the points in $D_0$ correspond to
irreducible stable Riemann surfaces with a double point, while for
$i=1,\dots,[\frac{g}{2}]$, $D_i$ classifies stable Riemann surfaces consisting
of two nonsingular irreducible components that intersect transversally, one of
genus $i$ and the other of genus $g-i$. The combinatorial structure of the
intersections of the $D_i$ is best described in terms of the {\it intersection
  graph}\index{intersection graph}: For a stable Riemann surface $X$, we define a graph $\Gamma(X)$ as
follows: the vertices of $\Gamma(X)$ are the irreducible components of $X$,
the edges are the double points (connecting two irreducible components of $X$
which need not be distinct). For every graph $\Gamma$ let $\Mgq(\Gamma)$ be the
set of points in $\Mgq$ corresponding to stable Riemann surfaces with
intersection graph isomorphic to $\Gamma$. It is not hard to see that for
a given genus $g$, there are only finitely many graphs $\Gamma$ with nonempty
$\Mgq(\Gamma)$, and that the $\Mgq(\Gamma)$ are the strata of a stratification
of $\Mgq$. This means that each $\Mgq(\Gamma)$ is a locally closed subset of
$\Mgq$ (for the Zariski topology), that $\Mgq$ is the disjoint union of the
$\Mgq(\Gamma)$, and that the closure of each $\Mgq(\Gamma)$ is a finite union
of other $\Mgq(\Gamma')$. A natural question in our context is: which
$\Mgq(\Gamma)$ contain cusps of Teichm\"uller curves? In \cite{michi} Maier showed that if $\Gamma$ has no ``bridge'', i.\,e.\ no edge $e$ such that $\Gamma-e$ is disconnected. the stratum $\Mgq(\Gamma)$ contains
points on a compactified Teich\-m\"uller curve $\bar C_\iota$ with a Teichm\"uller embedding $\iota$ that corresponds to an origami. M\"oller and
Schmith\"usen observed that this condition on the graph is necessary if the Teich\-m\"uller curve comes from a quadratic differential which is the square of a holomorphic 1-form (or equivalently from a translation structure on $X^*$).\\

Most of our knowledge about cusps of Teichm\"uller curves comes from studying boundary points of Teichm\"uller disks in a suitable extension of Teichm\"uller space. Several
different boundaries for Teichm\"uller space with very different properties
have been studied, like the Thurston boundary or the one coming from the Bers embedding. In the framework of this paper we look for a space $\Tgq$
in which $T_g$ is open and dense such that the action of the group $\Gamma_g$
extends to an action on $\Tgq$, and the quotient space $\Tgq/\Gamma_g$ is
equal to $\Mgq$. Such a space is the ``augmented'' Teichm\"uller space\index{augmented Teichm\"uller space} $\Tgd$
 introduced by Abikoff \cite{Abideg}. The points in $\Tgd$
are equivalence classes of pairs $(X,f)$, where $X$ is a stable Riemann
surface of genus $g$ and $f:\Cref\to X$ is a {\it deformation}\index{deformation map}. This is a
continuous surjective map such that there are finitely many loops
$c_1,\dots,c_k$ on $\Cref$ with the property that $f$ is a homeomorphism outside the $c_i$
and maps each $c_i$ to a single point $P_i$ on $X$. Abikoff defined a topology
on this space and showed that the quotient for the natural action of
$\Gamma_g$ on the pairs $(X,f)$ is the moduli space $\Mgq$ as a topological
space.\\

In his thesis \cite{VDiss}, Braungardt\index{Braungardt's construction} introduced the concept of a covering of
a complex manifold $S$ with cusps over a divisor $D$. He showed that under
mild assumptions on $S$ there exists a universal covering $\tilde X$ of this
type which extends the usual holomorphic universal covering of $S-D$ by
attaching ``cusps'' over $D$. $\tilde X$ is no longer a complex manifold or a
complex space, but Braungardt introduced a natural notion of holomorphic
functions in a neighbourhood of a cusp and thus defined a sheaf ${\cal
  O}_{\tilde X}$ of rings (of holomorphic functions) on $\tilde X$. In this
way $\tilde X$ is a locally complex ringed space, and the quotient map $\tilde
X\to\tilde X/\pi_1(S-D)=S$ is analytic for this structure. When applied to
$S=\Mgq$ and $D=\dMg$, Braungardt showed that the universal covering $\Tgq$ of
$\Mgq$ with cusps over $\dMg$ is, as a topological space with an action of
$\Gamma_g$, homeomorphic to Abikoff's augmented Teichm\"uller space. We shall
reserve the symbol $\Tgq$ in this article exclusively for this space
(considered as a locally ringed space). In Chapter \ref{volker}, we review
Braungardt's construction and results.\\

Our key technique to investigate boundary points of Teichm\"uller disks is the use of {\it Strebel rays}\index{Strebel ray}, see Definition \ref{sray}. By this we mean a geodesic ray in $T_g$ that corresponds by the construction in Section \ref{sr} to a Strebel quadratic differential on the Riemann surface $X$ at the starting point of the ray. A Strebel differential decomposes $X$ into cylinders swept out by horizontal trajectories. Mainly following \cite{M} we give in Section \ref{srendpoint} two explicit descriptions of the marked Riemann surfaces $(X_K,f_K)$ (for $K>1$) on a Strebel ray. This allows us to identify the boundary point $(X_\infty,f_\infty)$ at the ``end''\index{Strebel ray!end point of} of the ray as the stable Riemann surface that is obtained by contracting on $X$ the core lines of the cylinders in a prescribed way, see Sections \ref{endpoint} and \ref{conv}.\\

In the case that the Teichm\"uller embedding $\iota$ leads to a Teichm\"uller curve $C_\iota$ we show in Section \ref{bpofdisks} that all boundary points\index{Teichm\"uller disk!boundary points of} of $\Delta_\iota$ are obtained in this way. This shows in particular that all cusps of Teichm\"uller curves\index{Teichm\"uller curve!cusps of} are obtained by contracting, on a corresponding Riemann surface, the center lines of the cylinders of a Strebel differential. For the proof of this result we show that  the Teichm\"uller embedding $\iota$ can be extended\index{Teichm\"uller embedding!extension of} to a continuous embedding $\bar\iota:\HH\cup\{\mbox{cusps of}\ \Gammaquerm_{\iota}\} 
     \hookrightarrow \Tgq$, see Prop.~\ref{iquer}. Moreover, if the Veech group $\Grq$ is a lattice in $\pslzwei(\RR)$, the image of $\iotaquer$ is the closure of $\Delta_\iota$ in $\Tgq$, see Prop.~\ref{iquerfortc}.\\

Since $\Tgq$ has these cusp singularities at the boundary that prevent it from
being an ordinary complex space, whereas the boundary of $\Mgq$ is a nice
divisor in a projective variety, it is interesting to look at spaces that lie
properly between Teichm\"uller and moduli space and to ask for a boundary that
fits somehow in between $\Tgq$ and $\Mgq$. An example of such a space is
provided by the {\it Schottky space}\index{Schottky space} which goes back to the paper \cite{S} of F.~Schottky
from 1887. He studied discontinuous groups that are freely generated by
M\"obius transformations $\gamma_1,\dots,\gamma_g$ (for some $g\ge1$) chosen
in such a way that there are disjoint closed Jordan domains $D_1,D_1',\dots,D_g,D_g'$
such that $\gamma_i$ maps $D_i$ onto the complement of the interior of
$D_i'$. The Riemann surface of such a {\it Schottky group}\index{Schottky group} is compact of genus
$g$. It can be shown that every compact Riemann surface $X$ admits such a {\it
  Schottky uniformization}\index{Schottky uniformization} $X =\Omega/\Gamma$ (with $\Omega\subset\PP^1(\CC)$
open and $\Gamma$ a Schottky group), see Section \ref{s-struc}. The covering
$\Omega\to X$ is called a {\it Schottky covering}\index{Schottky covering}. 
It is minimal for the property that $\Omega$ is planar, i.\,e.\ 
biholomorphic to a subdomain of $\PP^1(\CC)$; here minimality means that
each unramified holomorphic covering $Y \to X$ with a planar manifold $Y$ 
factors through $\Omega$.\\

Schottky coverings are classified by a complex manifold $S_g$ of dimension
$3g-3$, called the Schottky space. The natural map from $T_g$ to $M_g$ factors
through $S_g$, therefore there is a subgroup $\Gga$ of the mapping class group
$\Gamma_g$ such that $T_g/\Gga=S_g$. Unfortunately the subgroup $\Gga$ is not
normal and depends on the choice of a certain group homomorphism $\alpha$. As
a consequence the induced map $S_g\to M_g$ is not the quotient for a group
action. We review this classical but not so widely known material in Sections
\ref{s-struc} and \ref{s-teich}.\\

The concept of Schottky coverings can be extended to stable Riemann
surfaces. If the analogous construction as for ordinary Riemann surfaces is
applied to a surface $X$ with nodes, we obtain a covering space $\Omega$ which
is not planar, but on which nevertheless a free group $\Gamma$ acts by
holomorphic automorphisms with quotient space $\Omega/\Gamma=X$. Although the
groups $\Gamma$ are no longer subgroups of $\pslzwei(\CC)$, it is possible to
find parameters for them in almost the same way as for Schottky groups, namely
by cross ratios of fixed points. It then turns out that these
generalized Schottky coverings\index{Schottky covering!generalized} are classified by a complex mani\-fold $\Sgq$
(which contains $S_g$ as an open dense subset), see Section
\ref{s-stab}. This result was originally proved in \cite{GH}; here we show
that it can easily be derived from Braungardt's characterization of $\Tgq$ as
the universal covering of $\Mgq$ with cusps over $\dMg$, see Section \ref{s-ext}.\\

Finally we wonder what the image in $S_g$ of a Teichm\"uller disk
$\Delta_\iota$ in $T_g$ might look like. In the general case we have no
idea. Of course, the image may depend on the choice of the subgroup $\Gga$
that gives the map $T_g\to S_g$. In the special situation that $C_\iota$ is a
Teichm\"uller curve we prove that for suitable choice of $\alpha$, the image
of $\Delta_\iota$ in $S_g$ is not a disk, see Prop.~\ref{durchschnitt}.\\

{\bf Acknowledgments:} We would like to thank Volker Braungardt for allowing us to include his results on $\Tgq$, and for his helpful comments on an earlier version of Chapter \ref{volker}. We are also grateful to Pierre Lochak and Martin M\"oller for many valuable conversations on Teichm\"uller disks, Teichm\"uller curves, and their boundaries. Furthermore, we would like to thank 
Bill Abikoff for his useful suggestions that helped to improve the 
exposition considerably.

\section{Geodesic rays, Teichm\"uller disks and Teichm\"uller curves}\label{bp}

The aim of this section is to  introduce Teichm\"uller 
disks and Teichm\"uller curves.
We start by recalling in \ref{deform} the concept  of
Teich\-m\"uller deformations and using them we give
a definition for the Teichm\"uller space $\tg$ alternative to the one we 
gave in the introduction. 
This will help us to define geodesic rays
in the Teich\-m\"uller space in \ref{sr}.
In \ref{defdisks} 
we introduce  Teichm\"uller disks as complex version of geodesic rays
giving different alternative definitions. Finally in \ref{tc}
we introduce the Veech group and Teichm\"uller 
curves and summarize some facts about the interrelation between these objects.

\subsection{Teichm\"uller deformations}\label{deform}

As one of numerous possibilities, one can define the Teichm\"uller space
as the space of Teichm\"uller deformations\index{Teichm\"uller deformations}. 
We briefly recall
this concept here. It is described e.g. in \cite[Ch I, \S 3]{A}.

At the end of this subsection we extend it to the corresponding concept
for punctured Riemann surfaces and their Teichm\"uller space $\tgn$, 
cf. \cite[Chapter II, \S 1]{A}.\\

Let $X = \Xref$ be a fixed Riemann surface of genus $g \geq 2$ and $q$ be a 
holomorphic quadratic differential on $X$.  We refer to the zeros of $q$ 
as {\em critical points}\index{holomorphic quadratic differential!critical point of}, all other points are {\em regular points}\index{holomorphic quadratic differential!regular point of}.
Then on the surface 
\[X^* = X - \{P \in X |\; P \mbox{ is a critical point of }q\}\]
the differential $q$  naturally defines a 
{\em flat structure}\index{flat structure} $\mu$, i.e.\
an atlas such that all transition maps are of the form 
$z \mapsto \pm z + c$, with 
some constant $ c \in \CC$. The charts of $\mu$
in regular points 
$z_0$ are given as
\begin{equation}\;\label{natchart}z \,\mapsto\, \int_{z_0}^z\sqrt{q(\xi)}d\xi.
\end{equation}

One may deform this flat structure 
 by composing each 
chart with the map 
\begin{equation} \label{Kdeform}
x + iy \;\; \mapsto \;\;  Kx + iy \quad 
   = \; \frac{1}{2}(K+1)z + \frac{1}{2}(K-1)\overline{z}, 
   \quad (x,y \in \RR)    
\end{equation}
with $K$ an arbitrary real number $>1$. This defines a new 
flat structure on $X^*$ which can uniquely be extended to a holomorphic 
structure on $X$.\\ 

We call $X_K$ the Riemann surface that we obtain this way,
$X_1 = X$ the surface with the original complex structure and
$f_K: X_1 \to X_K$ the map that is topologically the identity.
The map $f_K$ is a Teichm\"uller map and has constant complex dilatation \index{complex dilatation}
\[
k(z) = \frac{(f_K)_{\overline{z}}}{(f_K)_{z}} \;\; = \;\; \frac{K-1}{K+1}.\]
Its maximal dilatation\index{maximal dilatation} 
$\sup_{z \in X}\frac{1+|k(z)|}{1-|k(z)|}$ (as a quasiconformal map)  
is equal to $K$. 
\begin{definition}
Let $q$ be a holomorphic quadratic differential on $X$
and $K \in \RR_{>1}$.
The pair $(X_K,f_K)$ as defined above is called the 
{\em Teichm\"uller deformation}\index{Teichm\"uller deformation} of $X$ of 
constant dilatation $K$ with respect to $q$.
\end{definition}
The pair $(X_K, f_K)$ defines a point in the Teichm\"uller space 
$\tg$ which for simplicity we also denote as $(X_K,f_K)$.
Since the constant dilatation of $f_K$ is equal to $K$,  
the Teichm\"uller distance\index{Teichm\"uller distance} between the points $(X_1,\id)$ and $(X_K, f_K)$ 
of $\tg$ is $\log(K)$.\\ 

If two holomorphic quadratic differentials on $X$ 
are positive scalar multiples of each other, they define, for each $K$,
the same point in $\tg$. Thus one restricts to 
differentials with norm $1$. By Teichm\"uller's existence 
and uniqueness theorems, see e.g. \cite[Chapter I, (3.5), (4.2)]{A}, 
one can show that each point in $T_g$ is uniquely obtained as a
Teichm\"uller deformation.
If $\qdiff$ is the vector space of
all holomorphic quadratic differentials on $X$ and if $\qeinh$ 
is the unit sphere
in $\qdiff$, one may thus write
\begin{equation} \label{tgdeform}
\{(X,q,k)| \; q \in \qeinh, k \in (0,1)\} \cup  \{0\} \;= \; T_g.
\end{equation}
and the identification of the two sets is done by the map:
\begin{eqnarray} \label{tgtotdef}
   (X,q,k) & \mapsto & 
      (X_K,f_K) \quad \mbox{ with }
      K = \frac{1+k}{1-k} \Leftrightarrow k = \frac{K-1}{K+1}  \quad
                  \mbox{ and }\nonumber \\
      0   & \mapsto &  \mbox{ the base point } (X,\id)
\end{eqnarray}
$(X_K,f_K)$ depends of course by its definition on the differential $q$.
In the following we will denote the base point also by $(X,q,0)$.\\
The map (\ref{tgtotdef}) is a homeomorphism. 
Here on the 
left hand side of (\ref{tgdeform}) one takes the topology obtained
by identifying it with the open unit ball 
in $\qdiff$. It follows in particular, that $\tg$ is contractible.\\

Teichm\"uller deformations can be understood as 
{\em affine deformations}\index{affine deformation}
in the following sense:
Let us here and in the rest of the article identify $\CC$ with $\RR^2$
by the $\RR$-linear map sending $(1,i)$ to the standard basis of $\RR^2$.
Then the  map in (\ref{Kdeform}) is equal to the affine map
\[\bpm x\\y \epm \;\mapsto \; \bpm K & 0\\ 0&1\epm\cdot \bpm x\\y\epm \]
Since composing charts with a biholomorphic map does not change the point in 
Teich\-m\"uller
space, one obtains the same point $(X_K,f_K)$ in $\tg$ if 
one composes each chart of the flat structure $\mu$ on $X$ with 
the affine map
\begin{eqnarray}\label{XDK}
\bpm x\\y \epm \;\mapsto \; D_K\cdot\bpm x\\y\epm \;\; \mbox{ with } \,
   D_K =  \bpm \sqrt{K} & 0\\ 0&\frac{1}{\sqrt{K}}\epm \; \in \, \slzwei(\RR)
\end{eqnarray}
We will use the following notations which are compatible with those in 
Section~\ref{affdeforms} where we introduce the general concept 
of affine deformations.

\begin{definition}\label{defDK}
Let $X$ be a compact Riemann surface of genus $g$, $q$ a holomorphic quadratic
differential, $\mu$ the flat structure defined by $q$.
We call the flat structure defined by (\ref{XDK})\; $\mu_{D_K}$ 
and denote $(X,\mu)\circ D_K = (X,\mu_{D_K})$.
\end{definition}
Note that $(X,\mu_{D_K})$ is as Riemann surface isomorphic to $X_K$.
Thus the point $[(X,\mu_{D_K}),\id]$ in $\tg$ defined by the marking
$\id:X \to (X,\mu_{D_K})$ is equal to $(X_K,f_K)$.\\

Finally, let us turn to {\em Teichm\"uller deformations of punctured 
Riemann surfaces}\index{Teichm\"uller deformation!of punctured Riemann surface}:
The definition is done almost in the same way as in the case without 
punctures, see \cite[Chapter II, \S 1]{A}. 
Suppose that $g$ and $n$ are natural numbers with $3g-3+n>0$.
Let $X$ be a Riemann surface of genus $g$
with $n$ marked points $P_1$, \ldots, $P_n$,  and 
$\Xref = X_0 = X - \{P_1, \ldots, P_n\}$.\\ 
In this case, one uses {\em admissible holomorphic quadratic
differentials}\index{holomorphic quadratic differential!admissible} on $X_0$. They are by 
definition those meromorphic quadratic differentials 
on $X$ that restrict to a holomorphic quadratic differential on $X_0$
and have at each puncture either a simple pole or extend holomorphically across
the puncture, see \cite[Chapter II, (1.4)]{A}. The vector space of these
differentials is called $\qdiffpunct$. For $q \in \qdiffpunct$ we define the
{\em critical points}\index{holomorphic quadratic differential!critical point of} to be the marked points and all zeroes of $q$; the remaining
points are called {\em regular}\index{holomorphic quadratic differential!regular point of}. Now, the definition of Teichm\"uller deformation
is done exactly as before, just always 
replacing $\qdiff$ by $\qdiffpunct$. One obtains in the same way:
\begin{equation} \label{tgdeformpunct}
\{(X,q,k)| \; q \in \qeinhpunct, k \in (0,1)\} \cup  \{0\} \;= \; T_{g,n}.
\end{equation}
Here $\qeinhpunct$ is the unit ball in $\qdiffpunct$.

\subsection{Geodesic rays}\index{geodesic ray}
\label{sr}
Let $X = \Xref$ be a Riemann surface of genus g.
A holomorphic quadratic differential $q$ on $X$  
naturally defines a geodesic embedding of $\RR_{\geq 0}$ 
into $T_g$ with respect to the Teichm\"uller metric on $\tg$
as is described in the following.

\begin{definition}
Let $q$ be a holomorphic quadratic differential on $X$
and $\gamma$ the map:
\begin{equation}\label{geodesicray}
\gamma = \gamma_q: \left\{
  \begin{array}{lcl}
  [0,\infty) & \rightarrow & T_g\\
  t & \mapsto & (X_K,f_K) \; = \; (X,\mu_{D_K})
                  \; = \; (X,q,k) \\[1mm]
     &        &\phantom{(X_K,f_K) } 
                \mbox{with } \; K = e^t, \quad
                 k = \frac{K-1}{K+1}
  \end{array}\right.
\end{equation}
The image of $\gamma$ is called the {\em geodesic ray
in $\tg$ in direction of $q$ (or with respect to) $q$ starting at $(X,\emid)$}.
\end{definition}

Here we use the notation of the last section:
\[(X_K,f_K) \stackrel{\mbox{\fs Def. }\ref{defDK}}{=} (X,\mu_{D_K})  
           \stackrel{(\ref{tgtotdef})}{=} (X,q,k)\]
is the point in $\tg$
defined by the Teichm\"uller deformation of 
$X$ of dilatation $K$ with respect to $q$.
Recall from the last section that the distance between the 
two points $(X_K, f_K)$ 
and $(X,\id)$ in $\tg$
is $\log(K)$. Thus $\gamma$ is an isometric embedding.\\

In fact, from the description of $\tg$ given in (\ref{tgdeform})
one observes that all points in $\tg$ which have distance 
$\log(K)$ to the base point 
$(X,\id)$ are Teichm\"uller deformations of $X$ 
of constant dilatation $K$ with respect to a holomorphic quadratic 
differential. It follows that each isometric embedding of 
$[0,\infty)$ into $\tg$ is of the form (\ref{geodesicray}).

\subsection{Teichm\"uller disks}\label{defdisks}\index{Teichm\"uller disk}

In this section we define Teichm\"uller disks. They can be 
found defined under this name e.g. in \cite[p.~149/150]{n} and 
\cite[8.1-8.2]{GL}.  
One may find comprehensive overviews e.~g. in \cite{V} 
and \cite{EG}, or more recently \cite{mcm} and \cite{L},
to pick only a few of numerous references where they occur. 
We introduce them here
in detail comparing three different ways  
how to construct them. For completeness we have
included most of the proofs.

\begin{definition} \label{emb}
Let $3g-3+n > 0$.
A {\em Teichm\"uller disk $\Delta_{\iota}$} 
is the image of a holomorphic isometric
embedding
\[ \iota: \DD \hookrightarrow \tgn \]
of the complex unit disk  $\DD = \{z \in \CC| |z| < 1\}$ into the  
Teichm\"uller space. Here we take
the Poincar{\'e} metric of constant curvature $-1$ on $\DD$ and the 
Teichm\"uller metric on
$\tgn$. The embedding $\iota$ is also called {\em Teichm\"uller embedding}\index{Teichm\"uller embedding}. 
\end{definition}

Instead of the unit disk $\DD$ one may take as well the
upper half plane $\HH$ with the hyperbolic metric. We will switch between
these two models using the holomorphic isometry 
\begin{eqnarray} \label{uebergang} 
\trafo: \;\; \HH \rightarrow \DD,\quad t \mapsto \frac{i-t}{i+t}.
\end{eqnarray}
Thus Teichm\"uller disks are obtained equivalently
as images of holomorphic 
isometric embeddings $\HH \hookrightarrow \tgn$ of the upper half plane $\HH$
into the Teichm\"uller space $\tgn$.\\

How does one find such embeddings? 
Similarly as for geodesic rays, each holomorphic quadratic differential 
$q$ on a Riemann surface $X$ defines 
a Teich\-m\"uller disk. In the following we describe three alternative
constructions starting from such a differential
$q$ that all lead to the same Teichm\"uller disk $\Delta_q$. For simplicity
we only consider the case $n = 0$ and $g \geq 2$. However the same 
constructions
can be done in the general case of punctured surfaces.

\subsubsection[Teichm\"uller disks as a collection of geodesic rays]{
Teichm\"uller disks as a collection of geodesic rays}\label{coll}\index{Teichm\"uller disk}
\begin{definition}
Let $q$ be a holomorphic 
quadratic differential on a Riemann surface $X$ of genus $g$.
Let $\iota_1$ be the map
\begin{equation*}
\iota_1: \left\{
  \begin{array}{ccl}
   \DD  &\to& \tg\\
    z = r \cdot e^{i\varphi} &\mapsto& (X, e^{-i\varphi}\cdot q, r). 
  \end{array} \right.
\end{equation*}
\end{definition}

Here we use the definition of $\tg$ given by (\ref{tgdeform}).
Hence, $(X, e^{-i\varphi}\cdot q, r)$ is the point
defined by the Teichm\"uller deformation 
of $X$ of dilatation $K = \frac{1+r}{1-r}$ with respect to 
$q_{-\varphi} = e^{-i\varphi}\cdot q$. \\

We shall show in Proposition \ref{alleszs}  that 
$\iota_1$ is an isometric holomorphic
embedding, thus the image $\Delta_{\iota_1}$ of $\iota_1$ is 
a Teichm\"uller disk.\\

The map $\iota_1$ may be considered as a collection of geodesic rays 
in the following sense:
Let $\tau_{\varphi}$ be the geodesic ray in $\DD$ starting from $0$
in direction $\varphi$, i.e.:
\[ \tau_{\varphi}: 
\left \{\begin{array}{lcl}
 [0,\infty) &\to& \DD\\ 
    t         &\mapsto&  r(t)\cdot e^{i\varphi} 
   \quad \mbox{ with } r(t) = \frac{e^t - 1}{e^t + 1}
\end{array}
\right.
\]
Then $\iota_1 \circ \tau_{\varphi}: [0,\infty) \to \tg$ 
is equal to the map given in (\ref{geodesicray}) that defines the geodesic ray 
to the holomorphic quadratic differential $q_{-\varphi} = e^{-i\varphi}\cdot q$
on $X$.\\
Thus the Teichm\"uller disk $\Delta_{\iota_1}$ is the union of 
all geodesic rays defined by the differentials $e^{i\varphi}\cdot q$
with $\varphi \in [0,2\pi)$.
Furthermore,  $\iota_1 \circ \tau_{\varphi}$
is the parameterization by length of the restriction $\iota_1|_{R_{\varphi}}$
of $\iota_1$ to the ray $R_{\varphi} = \{r\cdot e^{i\varphi}|\,r \in [0,1)\}$.\\[3mm]

\subsubsection[Teichm\"uller disks by affine 
deformations]{Teichm\"uller disks by affine 
deformations\\[2mm]}\label{affdeforms}\index{Teichm\"uller disk}

We now describe a second approach that starting from
a holomorphic quadratic differential $q$
leads to the same Teichm\"uller disk
as in \ref{coll}.\\
Recall from Section \ref{deform} that a holomorphic
quadratic differential $q$  defines on 
$\Xstern = X - \{\mbox{zeroes of $q$}\}$ a flat structure 
$\mu$.
The group $\slzwei(\RR)$ acts on the flat structures 
of $\Xstern$ (as topological
surface) in the following way: 
Let $B \in \slzwei(\RR)$ and $\mu$ be a flat structure on 
$\Xstern$. Composing each chart of $\mu$ with the affine map
$z \mapsto B\cdot z$ gives a new flat structure on $\Xstern$ which we denote 
$B \circ (X,\mu)$
or $(X,\mu_B)$. In the special case $B = D_K$ 
we obtain the Teichm\"uller deformation of dilatation $K$, 
cf. Definition \ref{defDK}.

\begin{definition}\label{defdeform}
We call $(X,\mu_B) = B \circ (X,\mu)$ {\em affine deformation}\index{affine deformation}
of $(X,\mu)$ by the matrix $B$.
\end{definition}

\noindent
Note that for $B_1$, $B_2$ in $\slzwei(\RR)$ one may write
\[B_1\circ (B_2 \circ (X,\mu)) = B_1 \circ (X,\mu_{B_2}) = (X,\mu_{B_1B_2}) = 
B_1\cdot B_2 \circ(X,\mu).\]

\vspace*{2mm}

The flat structure $\mu_B$ defines in particular 
a complex structure on $X$. We identify here
the complex plane $\CC$ with $\RR^2$ as we already did in 
Section \ref{deform}.
In general the new complex structure will be 
different from the one defined by $\mu$. 
Taking the identity $\id:(X,\mu) \to (X,\mu_B)$ on 
$X$ as marking, we obtain a point $P_B = [(X,\mu_B),\id]$ 
in the Teichm\"uller space $\tg$. By abuse of notation we will
sometimes denote this point also just as $(X,\mu_B)$.\\
Thus one obtains the map
\[\hat{\iota}_2:\,\slzwei(\RR) \to \tg, \;\; B \; \mapsto \; 
   P_B = [(X,\mu_B),\id] = (X,\mu_B) \]
If however the matrix $A=U$ is in $\sozwei(\RR)$ the
map $\id: (X,\mu_B) \to (X,\mu_{U\cdot B})$ is holomorphic,
thus the point in Teichm\"uller
space is not changed, i.~e. 
\begin{equation}\label{PU}
U \in \sozwei(\RR) \;\; \Rightarrow \;\;
P_{UA} = P_A \mbox{ for all } A \in \slzwei(\RR)  
\end{equation}
Hence $\hat{\iota}_2$ 
induces a map 
\[\iota_2: \sozwei(\RR)\backslash\slzwei(\RR) \to \tg,\quad 
           \sozwei(\RR)\cdot B \; \mapsto \; 
         P_B = [(X,\mu_B),\id] = (X,\mu_B) .\] 

\vspace*{2mm}
Please note: The action of $\slzwei(\RR)$ on the flat structures
$\{(X,\mu_A)|\; A \in \slzwei(\RR)\}$ does  not  descend to the
image set $\{P_A|\; A \in \slzwei(\RR)\}$ in $\tg$;
in particular: 
$P_{U} = P_{I} \; \not\Rightarrow  \; P_{AU} = P_{A}$!\\

\noindent
{\bf The Teichm\"uller disk:}\\[1mm]
One may identify $\sozwei(\RR)\backslash\slzwei(\RR)$ 
with the upper half plane
$\HH$ in the following way:
Let $\slzwei(\RR)$ act by 
M\"obius transformations  
on the upper half plane $\HH$. This action is transitive and 
$\sozwei(\RR)$ is the stabilizing 
group of $i$. Thus the map 
\begin{equation}\label{pequ}
p:\;\slzwei(\RR) \to \HH, \;\; A \, \mapsto \, -\overline{A^{-1}(i)}
\end{equation}
induces a bijection $\sozwei(\RR)\backslash\slzwei(\RR) \to \HH$.
Its inverse map is induced by
\[\HH \to \slzwei(\RR),\quad t \, \mapsto \, 
  \frac{1}{\sqrt{\im(t)}}\bpm 1&\re(t)\\0&\im(t) \epm\] 
Composing $\iota_2$ from above with this bijection
one obtains a map from $\HH$ to $\tg$ which we also denote
by $\iota_2$.

\begin{definition}\label{defi2}
Let $q$ be a holomorphic quadratic differential on the
Riemann surface $X$ and $\mu$ the flat structure 
defined by $q$.
Let $\iota_2$ be the map
\[\iota_2:\, \HH \to \tg, \quad 
     t  \,\mapsto\,  P_{A_t} = [(X,\mu_{A_t}),\id]\] 
with $A_t$ chosen such that $-\overline{A_t^{-1}(i)} = t$.
\end{definition}

Note that the identification of $\sozwei(\RR)\backslash\slzwei(\RR)$ with 
$\HH$ 
given by $p$ may seem a bit ponderous, but
one has to compose $A \mapsto A^{-1}(i)$ with the reflection at
the imaginary axis in order that $\iota_2$ becomes holomorphic.
We will see this later in \ref{beltrams}.
In fact one has much more, as is stated in the
following proposition.

\begin{proposition} \label{alleszs}
The maps $\iota_1$ and $\iota_2$ are Teichm\"uller embeddings. 
They define the same Teichm\"uller disk 
\begin{eqnarray}\label{deltaq}
\Delta_q \; =  \; \Delta_{\iota_1} \;  = \;  \iota_1(\DD) 
    \; =  \; \Delta_{\iota_2}  \; = \;  
\iota_2(\HH).
\end{eqnarray}
\end{proposition}

\begin{proof}
The proof is given in the rest of Subsection \ref{affdeforms}
and in \ref{beltrams}:\\
In Proposition \ref{ioeinszwei} we show that 
$\iota_2 = \iota_1 \circ f$ with $f$ from (\ref{uebergang})
(see also Figure~\ref{gb}); thus it is sufficient
to show only for one of them that it is isometric,
and in the same manner for being holomorphic.\\
In Proposition \ref{remisom} it is  shown that $\iota_2$ 
is isometric. In Subsection \ref{beltrams}, it is shown 
that $\iota_1$ is holomorphic (see Corollary \ref{corhol}). 
For this purpose we introduce an
embedding $\iota_3:\DD \to \tg$, using Beltrami differentials,
for which it is not difficult to see that it is holomorphic, 
and show that it is equal 
to $\iota_1$.\\ 
That $\iota_1$ and $\iota_2$ define the same
Teichm\"uller disks then also follows also from Proposition 
\ref{ioeinszwei}.
\end{proof}

In fact the described constructions do not only give some special 
examples  but all Teichm\"uller disks 
are obtained as follows:
Each Teichm\"uller disk is equal to $\Delta_q$ as in
(\ref{deltaq}) for some  holomorphic quadratic differential $q$.
And all Teich\-m\"uller embeddings are of the form $\iota_1:\DD \into \tg$ 
or equivalently $\iota_2: \HH \into \tg$, see \cite[7.4]{GL}. \\

In order to see that $\iota_2$ from Definition \ref{defi2}
is isometric we first calculate the Teichm\"uller distance
between two affine deformations.\index{affine deformations!Teichm\"uller distance between two}\\[3mm]

\noindent
{\bf Teichm\"uller distance between two affine deformations:}\\[1mm]
In what follows we will constantly use the following fact
about matrices in $\slzwei(\RR)$:
\begin{remark}
Each matrix $A \in \emslzwei(\RR)$ with $A \not\in \emsozwei(\RR)$ 
can be decomposed uniquely
up to the minus signs as follows:
\begin{eqnarray}\label{decompose}
&A = U_1\cdot D_K \cdot U_2 \;\; \mbox{ with } U_1, U_2 \in \emsozwei(\RR),
\quad D_K = \bpm \sqrt{K} & 0 \\ 0 & \frac{1}{\sqrt{K}}\epm, \;\;K > 1.
\nonumber&\\
&\mbox{We may denote:} \quad 
  U_2 = U_{\theta} = 
    \bpm \cos(\theta)&-\sin(\theta)\\ \sin(\theta) & \cos(\theta) \epm.&
\end{eqnarray}
\end{remark}
This fact can e.g. be seen geometrically as follows: $\slzwei(\RR)$
acts transitively 
on the upper half plane $\HH$ by M\"obius transformations.
The point $i\in \HH$ can be mapped to $A(i)\neq i$ by first
doing a  stretching 
along the imaginary axis in direction $\infty$ and afterwards 
a rotation around $i$, i.~e. $A(i) = U_1(D_K(i))$
with suitably chosen $U_1 \in \sozwei(\RR)$ and $D_K$ with $K>1$ 
as in the remark. 
Since the stabilizer of $i$ in $\slzwei(\RR)$ is $\sozwei(\RR)$,
one has $A = U_1\cdot D_K \cdot U_2$ with $U_2$ also in $\sozwei(\RR)$.
A short calculation gives the uniqueness claim.\\

In the following proposition again let
$q$ be a holomorphic quadratic differential on 
$X = \Xref$ and $\mu$ the flat structure that $q$ defines.

\begin{proposition}\label{distance}
Let $A$ and $B$ be in $\emslzwei(\RR)$ with $A\cdot B^{-1} \not\in \emsozwei(\RR)$ and 
\[ A\cdot B^{-1} = U_1 \cdot D_K \cdot U_2 \]
with $U_1$, $U_2$ and $D_K$ as in (\ref{decompose}). Then 
the Teichm\"uller distance between the two points
$P_A = [(X,\mu_A),\emid]$ and $P_B = [(X,\mu_B),\emid]$ 
in $\tg$ is $\log(K)$.
\end{proposition}

\begin{proof}
We will proceed in three steps:\\

\noindent
{\bf a)} Suppose $B$ is the identity matrix $I$ and 
\[A \;=\; D_K \;=\; \bpm \sqrt{K} & 0\\ 0&\frac{1}{ \sqrt{K}} \epm 
  \;\;\mbox{ for some } 
K \in \RR_{>1}.\]
Thus we have in fact that $P_A = [(X,\mu_{D_K}),\id]$ is the point in $\tg$
defined by the Teichm\"uller deformation of dilatation $K$
with respect to $q$, see Definition \ref{defDK}. Hence the distance
between $P_A$ and the base point $(\Xref,\id) =  P_I$
is $\log(K)$.\\

\noindent
{\bf b)} Suppose again that $B = I$, but $A$ is an arbitrary
matrix in $\slzwei(\RR)$.\\ 
Thus \; $A = U_1 \cdot D_K \cdot U_2$ \; and
the map  $\id: (X,\mu) \to (X,\mu_A)$  is  the composition of three maps:
\[(X,\mu) \stackrel{\footnotesize\id}{\to} (X,\mu_{U_2}) 
          \stackrel{\footnotesize\id}{\to} (X,\mu_{D_KU_2})
          \stackrel{\footnotesize\id}{\to} (X,\mu_{U_1D_KU_2}) \]
Since the first and the third map are biholomorphic 
the Teichm\"uller distance is again $\log(K)$.\\
More precisely, write \; $U_2 = U_{\theta}$ \;
as in (\ref{decompose}).
Then $\mu_{U_2}$ is the flat structure obtained by
composing each chart with $z \mapsto e^{i\theta}\cdot z$. This
is equal to the flat structure defined by the quadratic differential  
$q_{2\theta} = (e^{i\theta})^2\cdot q$ which is holomorphic on 
the Riemann surface $X$.\\ 
Now,\; $\id:(X,\mu_{U_2}) \to (X,\mu_{D_KU_2})$\; is (up to the stretching
$z \mapsto \sqrt{K}\cdot z$)
the Teich\-m\"uller deformation of dilatation $K$ with respect to 
the holomorphic quadratic differential 
$q_{2\theta}$. Thus the distance between $P_A = P_{U_1D_KU_2} 
\stackrel{(\ref{PU})}{=} P_{D_KU_2}$ and
the base point $P_B = P_I$ is $\log(K)$.\\

\noindent
{\bf c)} Let now $A$, $B$ be arbitrary in $\slzwei(\RR)$. The Teichm\"uller
metric does not depend on the chosen base point. Thus we may consider
$P_B$ as base point and $P_A$ as coming from the affine deformation defined
by the matrix $A\cdot B^{-1}$. 
Then with the given decomposition $A\cdot B^{-1} = U_1 \cdot D_K \cdot U_2$ 
the distance is as in b) equal to $\log(K)$.
\end{proof}

\begin{proposition} \label{remisom}
$\iota_2$ is an isometric embedding
\end{proposition}

\begin{proof}
We denote by $\rho$ the Poincar{\'e} distance in $\HH$
and by $d_T$ the Teichm\"uller distance in $\tg$. Let
$t_1$ and $t_2$ be arbitrary distinct points in $\HH$. We 
may write $t_1 = p(A)$ and $t_2 = p(B)$ with  $A$, $B$ in
$\slzwei(\RR)$, $p$ as in (\ref{pequ}). Let 
$AB^{-1} = U_1D_KU_2$ the decomposition of $AB^{-1}$
from (\ref{decompose}). ($AB^{-1} \notin \sozwei(\RR)$ because 
$t_1 \neq t_2$)
\begin{eqnarray*} 
\rho(t_1,t_2) &=& \rho(-\overline{B^{-1}(i)}, -\overline{A^{-1}(i)}) 
  = \rho(B^{-1}(i), A^{-1}(i)) = \rho(AB^{-1}(i),i)\\
 &=& \rho(U_1D_KU_2(i),i) 
  =  \rho(U_1D_K(i),i) \stackrel{\star}{=} \rho(D_K(i),i)\\ 
 &=& \rho(Ki,i)\, =\, \log(K) \stackrel{\mbox{\fs Prop. \ref{distance}}}{=} 
d_T(P_B,P_A) = d_T(\iota_2(t_1),\iota_2(t_2)) 
\end{eqnarray*}
The equality $\star$ is given since $U_1$ is a hyperbolic
rotation with center $i$ and thus does not change the distance
to $i$. 
\end{proof}

Now we show that $\iota_1$ and $\iota_2$
are ``almost'' the same map.

\begin{proposition}\label{ioeinszwei}
$\iota_1$ and $\iota_2$ fit together. More precisely:
$\iota_1 \circ \trafo = \iota_2,$
with the isomorphism $\trafo:\HH \to \DD$  from (\ref{uebergang}).
\end{proposition}

The following diagram may be helpful while reading the
proof. Some parts will be explained only after the proof; in particular 
the space $B(X)$ of Beltrami differentials will be introduced in 
\ref{beltrams}.\\

\begin{minipage}{\linewidth}

\[\hspace*{-5mm} \xymatrix{
     &  \DD  \ar[rd]^{b} \ar@/^1.5cm/[rrd]^{\iota_1}
          \save[]+<-22mm,2mm> *\txt<8pc>{%
          ${\scriptstyle
                    z = r\cdot e^{\alpha i} = \frac{i-t}{i+t} 
                    = \frac{K-1}{K+1}\cdot e^{-2i\theta} \in }$}
         \restore
     & & \\
         \slzwei(\RR) 
           \ar[rd]^{p}
           \ar@/_2.7cm/[rrr]_{\hat{\iota}_2}
             \save[]+<9mm,4mm> *\txt<8pc>{%
             $ \scriptstyle \ni \;\; A \; =\; U_1D_KU_{\theta} $
             } 
             \restore
      & &
      B(X) \ar[r]^{\Phi} 
      \save[]+<4mm,7mm> *\txt<8pc> {
            $z\cdot\frac{\bar{q}}{|q|} = 
              \frac{i-t}{i+t}\frac{\bar{q}}{|q|}  
                  $}\restore            
      & 
        **[r]\hspace*{3mm}\tg \; \ni &{}
         \save[]+<15mm,0mm>*\txt<5cm>{%
          $\left\{
          \begin{minipage}{4cm} 
             $\iota_1(z)$ \newline
             $\scs =  (X,e^{-\alpha i}q,r)$\newline
             \quad \; $\scs = (X,e^{2\theta i}q,
                \frac{\sscs K-1}{\sscs K+1}) =$  
             \newline
             $\iota_2(t)$ \newline
             $\scs = P_A = [(X,\mu_A), \id]$  \newline
             $\scs = [D_K \circ U_2 \circ (X,\mu), \id]$
             \end{minipage}\right.$}
         \restore \\
     &  \hspace*{1cm}\HH  \hspace*{1cm}
          \ar[ru] \ar[uu]^{\trafo} 
        \save[]+<14mm,-2mm>*\txt<8pc>{%
        $\scs \ni t = -\overline{A^{-1}(i)}$}
        \restore     
     & & 
  }\]
\begin{center}
\refstepcounter{diagramm}{\it Figure \arabic{diagramm}}:
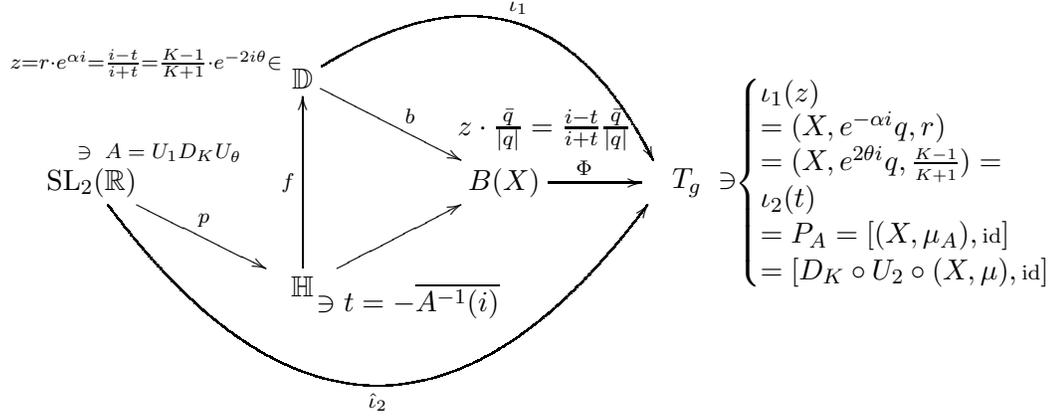
{\it Diagram for alternative definitions of the Teichm\"uller disk
$\Delta_q$}
\label{gb}
\end{center}
\end{minipage}
\vspace*{3mm}

\begin{proof}
We proceed in two steps:

\begin{enumerate}
\item 
Let $A \in \slzwei(\RR)$ be decomposed as in (\ref{decompose}):
$A = U_1\cdot D_K \cdot U_2$, \; $U_2 = U_{\theta}$.\\[1.5mm]
We show that $(f\circ p)(A) = r\cdot e^{-2i\theta}$ with $r = \frac{K-1}{K+1}$.
\item We show that $\iota_1(r\cdot e^{-2i\theta}) = \hat{\iota}_2(A)$.\\
\end{enumerate}

\noindent
{\bf Step 1:} 
One may express $t := p(A)$ in terms of $K$ and $\theta$ as follows:
\begin{eqnarray*}
t &=& -\overline{A^{-1}(i)} = -U_2^{-1}D_K^{-1}(\overline{i})
  = -U_2^{-1}(-\frac{i}{K}) 
  = -\frac{\cos(\theta)\cdot \frac{-i}{K} + \sin(\theta)}
          {-\sin(\theta)\cdot\frac{-i}{K} + \cos(\theta)}\\
  &=&  \frac{i\cos(\theta) - K\sin(\theta)}{i\sin(\theta)+K\cos(\theta)}
\end{eqnarray*}
Now one has:
\begin{eqnarray*}
f(p(A)) 
   &=& f(t) \,=\, \frac{-t+i}{t+i} 
   \,=\,  \frac{-i\cos(\theta) + K\sin(\theta) +i(i\sin(\theta)+K\cos(\theta))}{
          i\cos(\theta) - K\sin(\theta)  +i(i\sin(\theta)+K\cos(\theta))}\\
   && \hspace*{-20mm}
    \,=\, \;\,
      \frac{(K-1)[\sin(\theta)+i\cos(\theta)]}
           {(K+1)[-\sin(\theta)+i\cos(\theta)]}
    \,=\,   \frac{K-1}{K+1} \cdot
        \frac{-(\sin(\theta)+ i\cos(\theta))^2}
           {(\sin(\theta)-i\cos(\theta))(\sin(\theta)+i\cos(\theta))}  \\
   && \hspace*{-20mm}
    \,=\, \;\,
       \frac{K-1}{K+1}(\cos(\theta) - i\sin(\theta))^2
    \,=\,  \frac{K-1}{K+1}\cdot e^{-2i\theta}
\end{eqnarray*}
\noindent
{\bf Step 2:}
$\iota_1(r\cdot e^{-2i\theta}) = (X,e^{2i\theta}\cdot q, r) \in \tg$
is the point in the Teichm\"uller space that is obtained as 
Teichm\"uller deformation of dilatation $\frac{1+r}{1-r} = K$ with respect to
the quadratic differential $e^{2i\theta}\cdot q$. 
Recall from the proof of Proposition \ref{distance} that this is 
precisely the point in 
$\tg$ defined by the affine deformation 
$D_K \circ U_{\theta} \circ (X,\mu) = (X,\mu_{D_KU_{\theta}}) 
= (X, \mu_{D_KU_2})$. Thus
\begin{equation}\label{DUK}
(X,e^{2i\theta}\cdot q, r)  = P_{D_KU_{\theta}} =
P_{D_KU_2} \stackrel{(\ref{PU})}{=} 
P_{U_1D_KU_2} = P_A = 
   \hat{\iota}_2(A).
\end{equation}
\end{proof}

Using (\ref{DUK}) one may also describe the geodesic rays 
$\iota_1\circ \tau_{\varphi}$ from \ref{coll} in 
the Teichm\"uller disk 
$\Delta_q = \Delta_{\iota_1} = \Delta_{\iota_2}$ as follows.

\begin{corollary}\label{DKU}
Define $D_K, U_{\theta}$ as in (\ref{decompose}).
The map
\[
[0,\infty ) \;\to\; \tg, \quad t \;\mapsto  \;
   P_{D_KU_{\theta}} = [D_{K} \circ (X,\mu_{U_\theta}), \emid]
\;\mbox{ with } K = e^t
\]
is equal to $\iota_1\circ \tau_{-2\theta}$.\\ It is thus by \ref{coll}
the geodesic ray\index{geodesic ray} in direction of the quadratic differential 
$q_{2\theta} = e^{2\theta i}q$.
\end{corollary}

\begin{proof}
One has: \;\;
$ t \;\stackrel{\tau_{-2\theta}}{\mapsto}\; r(t)e^{-2\theta i} 
    \;\stackrel{\iota_1}{\mapsto}\;  (X,e^{2\theta i}\cdot q, r(t)) 
    \;\stackrel{(\ref{DUK})}{=}\; P_{D_KU_{\theta}}.
$
\end{proof}

Hence, geometrically one obtains the geodesic ray to $q_{\varphi}$
by rotating the flat structure by $U_{\frac{\varphi}{2}}$ and then stretching
in vertical direction with dilatation $K$.

\subsubsection[Beltrami differentials]{Beltrami differentials\\[2mm]}\label{beltrams}
In order to see that $\iota_1$ and $\iota_2$ are holomorphic
we introduce an alternative way to define $\iota_1$
using Beltrami differentials\index{Beltrami differential}. We keep this aspect
short and refer to e.g. \cite{n} for more details.\\
Let
\[M(X) = \{(X_1,f)|\; 
   \begin{array}[t]{l}
   X_1 \mbox{ Riemann surface },\\
   f:X\to X_1 \mbox{ is a quasiconformal homeomorphism}
   \}/\approx \end{array}\]
with $(X_1,f_1) \approx (X_2,f_2) \Leftrightarrow f_2\circ f_1^{-1}$ 
is biholomorphic.\\
One has a natural projection $M(X) \to \tg$. Furthermore $M(X)$ 
can be ca\-nonical\-ly identified with the open unit ball 
$B(X)$ in the Banach space 
$L^{\infty}_{(-1,1)}(X)$ of $(-1,1)$-forms by the bijection:
\[M(X) \to B(X),\;\; (X_1,f) \mapsto \mu_f, \]
where $\mu_f$ is the Beltrami differential (or complex dilatation)  
of $f$, cf. \cite[2.1.4]{n}.\\
Thus one obtains a projection 
$\Phi:  B(X) \to \tg$.
The map $\Phi$ is holomorphic (\cite[3.1]{n}). 
Furthermore, for each quadratic differential $q$ 
and for all $k \in (0,1)$ the form $k\frac{\bar{q}}{|q|}$ 
is in $B(X)$ 
(\cite[2.6.3]{n}) 
Thus one may define the map
\[\iota_3: \left\{
  \begin{array}{lclcl}
   \DD &\stackrel{b}{\to}& B(X) &\stackrel{\Phi}{\to}& \tg\\
    z  &\mapsto& z\cdot\frac{\bar{q}}{|q|} & \mapsto &
        \Phi(z\cdot\frac{\bar{q}}{|q|})
   \end{array} \right . 
\]
It is composition of two holomorphic maps and thus itself holomorphic.\\
We will show in the following remark that $\iota_3 = \Phi \circ b = \iota_1$, 
cf.\ Figure \ref{gb}.

\begin{remark}
For all $z_0 \in \DD: \iota_3(z_0) = \iota_1(z_0)$.
\end{remark}

\begin{proof}
Let $z_0 = r\cdot e^{i\alpha} \in \DD$ and $A \in \slzwei(\RR)$
with $f(p(A)) = z_0$.\\ 
Decompose $A = U_1D_KU_2$ as in (\ref{decompose}) with $U_2 = U_{\theta}$.
Then by Step 1 of the proof of  Proposition \ref{ioeinszwei}, 
$r = \frac{K-1}{K+1}$ \, and $\alpha = -2\theta$. Furthermore, 
by Proposition \ref{ioeinszwei} 
\[\iota_1(z) = \hat{\iota}_2(A) = [(X,\mu_A), \id] = [(X,\mu_{D_KU_2}),\id]\]
Let us calculate the Beltrami differential of the Teichm\"uller deformation 
$f = \id: X \to (X,\mu_{D_KU_2})$. We will see that it is equal to 
$z_0\cdot \frac{\bar{q}}{|q|}$. From this it follows that 
$\iota_1(z) = \iota_3(z)$.\\[2mm]
One has $f = g\circ h$ with $h = \id:X \to (X,\mu_{U_2})$ and 
$g = \id: (X,\mu_{U_2}) \to  (X,\mu_{D_KU_2})$. Locally in the  
charts of
the flat structure defined by $q$, the maps $g$ and $h$ are given by
\[g: z \mapsto K\cdot\re(z) + i\cdot\im(z) \;\quad\;\mbox{ and }\;\quad\;
  h: z \mapsto e^{i\theta}\cdot z.\] 
Thus in terms of these charts one has:
\begin{eqnarray*}
&&f_z \;=\; {g}_z\cdot h_z + g_{\bar{z}}\cdot \bar{h}_{z} 
      \;=\; e^{i\theta}\cdot g_z \;\;\qquad 
f_{\bar{z}} = g_z\cdot h_{\bar{z}} + g_{\bar{z}}\cdot \bar{h}_{\bar{z}}
    \;=\; e^{-i\theta}\cdot g_{\bar{z}} \hspace*{2cm}\\
&&\Rightarrow \quad \frac{f_{\bar{z}}}{f_z} 
    \;=\; e^{-2i\theta}\cdot\frac{g_{\bar{z}}}{g_z}
    \;=\; e^{-2i\theta}\cdot\frac{K-1}{K+1} 
    \;=\; e^{i\alpha}\cdot r \;=\; z_0
\end{eqnarray*}
Hence the Beltrami differential of $f$ is $z_0\cdot\frac{\bar{q}}{|q|}$.
\end{proof}

One obtains immediately the following conclusion.

\begin{corollary}\label{corhol}
$\iota_1 = \iota_3$ is holomorphic. By Proposition \ref{ioeinszwei} 
$\,\iota_2$ is also holomorphic.
\end{corollary}

\subsection[Teichm\"uller curves]{Teichm\"uller 
curves\\[2mm]}\label{tc}\index{Teichm\"uller curve}
In this section we introduce Teichm\"uller curves and recall some
properties of  them, in particular their relation to Veech groups. 
This was explored by Veech in his article \cite{V} and has been studied
by many authors since then. Overviews and further properties
can be found e.g. in \cite{mcm}, \cite{EG} or \cite{HS}. \\

Let $\iota: \DD \hookrightarrow \tg$ be a Teichm\"uller embedding and
$\Delta = \Delta_{\iota} = \iota(\DD)$ its image. We may consider
the image of $\Delta_{\iota}$ in the moduli space $\mg$ 
under the natural projection 
$\tg \to \mg$, cf. Chapter \ref{intro}. 
In general it will be something with a large closure. But 
occasionally it is an algebraic curve. Such a curve is
called Teichm\"uller curve.

\begin{definition} If the image of 
the Teichm\"uller disk $\Delta$ in 
the moduli space $\mg$ is an algebraic curve $C$, 
then $C$ is called {\em Teichm\"uller curve}.\\ 
A surface $(X,q)$, with a Riemann surface $X$ and
a holomorphic quadratic differential $q$
such that the Teichm\"uller
disk $\Delta = \Delta_q$ defined by $q$  
projects to a Teichm\"uller curve
is called {\em Veech surface}\index{Veech surface}.
\end{definition}

How can one decide whether a surface $(X,q)$ induces a Teichm\"uller 
curve or not? An answer to this question is given
by the Veech group, a subgroup of $\slzwei(\RR)$ associated to $(X,q)$.
This is explained in the following two subsections.

\subsubsection[Veech groups]{Veech groups\\[2mm]}\label{vg}\index{Veech group}
Let $X$ be a Riemann surface and $q$ a holomorphic quadratic 
differential on $X$. Let $\mu$ be the flat structure on $X$ defined by $q$.
One obtains a discrete subgroup of $\slzwei(\RR)$ 
as follows: Let $\affplus(X,\mu)$ be the group
of orientation preserving diffeomorphisms 
which are affine\index{affine diffeomorphism} with respect to the flat structure
$\mu$, i.e. diffeomorphisms which are in terms of a local chart $z$
of $\mu$  
given by
\[z  \mapsto A\cdot z  + t, \quad \mbox{ for some }
  A = \bpm a&b\\c&d\epm \in \slzwei(\RR), t \in \CC.\]
As above we identify the complex plane 
$\CC$ with $\RR^2$. Furthermore, we denote for $z=x+iy$: $A\cdot z = ax+by + i(cx + dy)$.\\
Since $\mu$ is a flat structure, up to change of sign the matrix $A$ does not 
depend on the charts. Thus one has a group homomorphism:
\[D: \;\affplus(X,\mu) \to \pslzwei(\RR), \;\; f \mapsto [A].\]
For simplicity we will denote the image $[A]$ of the matrix $A$
in $\pslzwei(\RR)$ often also just by $A$.
\begin{definition} \label{Veechgroup}
The image $\bar{\Gamma}(X,\mu) = D(\affplus(X,\mu))$ of $D$ 
is called the {\em projective Veech group} of $(X,\mu)$.
\end{definition}

We will denote the projective Veech group also by $\bar{\Gamma}(X,q)$ 
and $\bar{\Gamma}_{\iota}$, 
where $\iota:\DD \into \tg$ or $\iota:\HH \into \tg$ is the Teichm\"uller
embedding defined by $q$ as described in \ref{defdisks}.
$\bar{\Gamma}(X,\mu)$ is a discrete subgroup of $\pslzwei(\RR)$, 
see \cite[Prop. 2.7]{V}.

\subsubsection[The action of the Veech group on the Teichm\"uller disk]{
The action of the Veech group on the Teichm\"uller disk\\[2mm]}
Recall  that the projection 
$\tg \to \mg$ from the Teich\-m\"uller space to the moduli space
is given by the quotient
for the action of the mapping class group \index{mapping class group}
\[\Mod = \mbox{Diffeo}^+(X)/\mbox{Diffeo}_0(X),\]
cf. (\ref{mapgroup}) in the introduction.
The action of $\mbox{Diffeo}^+(X)$ on $\tg$ is given by
\begin{eqnarray*}
&\rho:&\qc(X) \;\to\;\aut(\tg) \cong \Gamma_q, \;\;
   \varphi  \mapsto  \rho_{\varphi}\\
&&\mbox{with } \; \rho_{\varphi}:\; \tg \to \tg, \;\; 
         (X_1,h)  \mapsto (X_1,h\circ\varphi^{-1}).
\end{eqnarray*}
The affine group\index{affine group} $\affplus(X,\mu)$ 
acts as subgroup of $\qc(X)$ on $\tg$.
The following remark (cf. \cite[Theorem 1]{EG}) determines
this action restricted to the Teichm\"uller disk 
\[\Delta = \Delta_q = 
\{P_B = [(X,\mu_B),\id] \in \tg| B \in \slzwei(\RR)\}.\]

\begin{remark} \label{actab}
$\emaffplus(X,\mu)$ stabilizes $\Delta$. Its action
on $\Delta$ is given by:
\begin{eqnarray*}
\varphi \in \emaffplus(X,\mu),\, B \in \emslzwei(\RR) &\Rightarrow&
\rho_{\varphi}(P_B)  = 
P_{BA^{-1}}\\ 
&&\mbox{ with } A \in \emslzwei(\RR)
    \mbox{ a preimage  of } 
D(\varphi) = [A].
\end{eqnarray*}
\end{remark}

\begin{proof}
Let $\varphi \in \affplus(X)$, $B \in \slzwei(\RR)$ and
$A \in \slzwei(\RR)$ be a preimage of $D(\varphi) = [A] \in \pslzwei(\RR)$.
In the following commutative diagram
\[
\xymatrix @-1pc {
(X,\mu)  \ar[rr]^{\varphi^{-1}} \ar[rrrrdd]_{\mbox{\fs id}}  && 
        (X,\mu)  \ar[rr]^{\mbox{\fs id}} && (X,\mu_B) \\&&&&\\
    &     &&              & (X,\mu_{BA^{-1}}) \ar[uu]
}
\]
the map $(X,\mu_{BA^{-1}}) \, \to \, (X,\mu_B)$ is, as a composition
of affine maps, itself affine. Its derivative is 
$D(\id \circ \varphi^{-1} \circ \id^{-1}) = 
BA^{-1}(B{A}^{-1})^{-1} = I$. 
Thus it
is biholomorphic and $\rho_{\varphi}([(X,\mu_B),\id])  = 
[(X,\mu_{BA^{-1}}),\id]$.
\end{proof}

It follows from Remark \ref{actab}
that $\affplus(X,\mu)$ is mapped by $\rho$
to Stab$(\Delta)$, the global stabilizer\index{Teichm\"uller disk!global stabilizer of} of $\Delta$
in $\Gamma_g$. Furthermore 
$\rho: \affplus(X,\mu) \to \mbox{Stab}(\Delta)\,\subseteq \Gamma_g$ is in fact
an isomorphism: It is injective, see \cite[Lemma~5.2]{EG}
and surjective, see \cite[Theorem~1]{EG}. Thus we 
have $\affplus(X,\mu) \cong \mbox{Stab}(\Delta)$.\\

From Remark \ref{actab} it also becomes clear that
the action of $\varphi \in \affplus(X,\mu)$ 
depends only on $D(\varphi)$. Thus one obtains in fact 
an action of the projective Veech group 
$\bar{\Gamma}(X,\mu)$ on $\Delta$.

\begin{corollary}
$\Gammaquer(X,\mu) \, \subseteq \, \empslzwei(\RR)$ acts on 
$\Delta = \{P_B  \in \tg|\, B \in \emslzwei(\RR)\}$
by:
\begin{equation}\label{abc}
\rho_{[A]}(P_B) = P_{BA^{-1}} \quad \mbox{ where 
$A$ is a preimage in $\emslzwei(\RR)$ of $[A]$}. 
\end{equation}
\end{corollary}

Finally one may use the Teichm\"uller embedding 
$\iota_2:\HH \to \tg$ defined by $q$  (cf. \ref{defi2})
in order to compare the action of $\bar{\Gamma}(X,\mu)$
on $\Delta = \Delta_{\iota} = \iota(\HH)$ with its action on 
$\HH$ via M\"obius 
transformations. One obtains the diagram in the following remark
(cf. \cite[Proposition 3.2.]{mcm}).

\begin{remark} \label{achtionstogether}
Let $A \in \empslzwei(\RR)$. Denote by $A:\,\HH \to \HH$ 
its action as M\"obius transformation on $\HH$. The following diagram
is commutative:
\[ \xymatrix @-1pc {
\HH \ar[rr]^{t \mapsto -\bar{t}} \ar[d]_{A}
  && \HH \ar[rr]^{\iota}\ar[d]^{RAR^{-1}} && \Delta \ar[d]^{\rho_A} \\
\HH \ar[rr]^{t \mapsto -\bar{t}} && \HH \ar[rr]^{\iota} && \Delta\\
}\]
\begin{center} 
\refstepcounter{diagramm}{\it Figure \arabic{diagramm}}
\label{diagramactions}
\end{center}
Here $R = \bpm -1 & 0 \\0&1\epm$, thus $R$ acts on $\PP^1(\CC)$
by $z \mapsto -z$.
\end{remark}
\begin{proof}
Let $t \in \HH$. Choose some $B \in \slzwei(\RR)$ with 
$-\overline{B^{-1}(i)} = -\bar{t}$,
thus $\iota(-\bar{t}) = P_B = [(X,\mu_B),\id]$ and using (\ref{abc})
we obtain the diagram:
\[ \xymatrix @-1pc {
t \ar@{|->}[rr]^{t \mapsto -\bar{t}} \ar@{|->}[d]_{A}
  && -\bar{t} \ar@{|->}[rr]^(0.3){\iota} 
  &&   P_B = [(X,\mu_B),\id] \ar@{|->}[d]^{\rho_A} \\
 A(t) \ar@{|->}[rr]^{t \mapsto -\bar{t}} && -\overline{A(t)}  && 
 P_{BA^{-1}}= [(X,\mu_{BA^{-1}}),\id]\\
}\]

The commutativity of the diagram in Figure \ref{diagramactions}
then follows from 
\begin{eqnarray*}
&&RAR^{-1}(-\bar{t}) = -A(\bar{t}) = -\overline{A(t)} \mbox{ and}\\
&&-\overline{(BA^{-1})^{-1}(i)} = 
  -A(\overline{B^{-1}(i)}) = -A(\bar{t}) = -\overline{A(t)} 
\mbox{, thus } \iota(-\overline{A(t)}) = P_{BA^{-1}}.
\end{eqnarray*}
\end{proof}

\subsubsection[Veech groups and Teichm\"uller 
  curves]{Veech groups and Teichm\"uller curves\\[2mm]}\label{lattice}\index{Veech group}\index{Teichm\"uller curve}
In Remark \ref{actab} we saw that the affine group $\affplus(X,\mu)$
maps isomorphically to the global stabilizer of the Teichm\"uller disk 
$\Delta$ in $\Mod$.
Denote by $\proj:\tg \to \mg$ the canonical projection. It then 
follows from Remark \ref{achtionstogether}
that the map 
\[ \proj\circ \iota:\;\;\HH \;\,\to\;\,  \proj(\Delta) \;\,\subseteq\; \mg\]
factors through $\HH/R\bar{\Gamma}(X,\mu)R^{-1}$.
We call 
\[\Gammaquerm(X,\mu) = R\bar{\Gamma}(X,\mu)R^{-1}\]
the {\em mirror projective Veech group}\index{mirror Veech group}, since $\HH/\Gammaquerm(X,\mu)$
is a mirror image of $\HH/\Gammaquer(X,\mu)$, and 
refer to it also as $\Gammaquerm(X,q)$ or $\Gammaquerm_{\iota}$.\\
$\HH/\Gammaquerm(X,\mu)$ is a surface of finite type and hence an algebraic
curve if and only if $\Gammaquerm(X,\mu)$ is a lattice in $\pslzwei(\RR)$.
Altogether one obtains the following situation 
(cf. \cite[Corollary 3.3]{mcm}).

\begin{corollary} \label{latticeproperty}
$(X,q)$ induces a Teichm\"uller curve $C$ if and only if $\Gammaquer(X,\mu)$
is a lattice in $\empslzwei(\RR)$. In this case the  
following diagram holds:
\[ \xymatrix @-1pc {
\HH \ar[rr]^{t \,\mapsto\, -\bar{t}} \ar[d]
  && \;\; \HH \;\; \ar[d] \ar[rr]^{\iota}
  && \; \Delta = \Delta_{\iota} \;\; \subseteq \;\; \tg 
     \ar@<-18pt>[d]^{\mbox{\fs \em proj}} \ar@<28pt>[d]^{\mbox{\fs \em proj}}\\
\HH/\Gammaquer(X,\mu) \ar[rr]^{\mbox{\fs antihol.}} 
  && \;\; \HH/\Gammaquerm(X,\mu) \; \ar[rr]^(.44){\mbox{\fs birat.}}
  && \quad \;\;  C 
  \quad \quad \subseteq \;\; \mg
}\]
In particular if $\Gammaquer(X,\mu)$
is a lattice, then 
\begin{itemize}
\item 
$\HH/\Gammaquerm(X,\mu)$ is the normalization
of the Teichm\"uller curve $C$,
\item  $\HH/\Gammaquer(X,\mu)$ is antiholomorphic
to $\HH/\Gammaquerm(X,\mu)$.
\end{itemize}
\end{corollary}

\section{Braungardt's construction of $\Tgnq$}
\label{volker}
Before we continue our study of Teichm\"uller disks and pass to the boundary, 
we want to explain the partial compactification $\Tgnq$\index{Teichm\"uller 
space!partial compactification $\Tgnq$} of the Teichm\"uller space $\Tgn$ that 
we shall use in the subsequent chapters. As mentioned in the introduction, 
$\Tgnq$ will be a locally ringed space which, as a topological space, 
coincides with Abikoff's augmented Teichm\"uller space\index{Teichm\"uller 
space!augmented}\index{augmented Teichm\"uller space} $\Tgnd$ (see the 
discussion following 
Proposition \ref{homeo}). The points of this space can be considered as marked 
stable Riemann surfaces $(X,f)$, where $f:\Cref\to X$ is a deformation map. 
The forgetful map $(X,f)\mapsto X$ defines a natural map from $\Tgnq$ to the 
moduli space $\Mgnq$ of stable $n$-pointed Riemann surfaces of genus $g$. 
This map extends the projection $\Tgn\to\Mgn$ and is in fact also the quotient 
map for the natural action of the mapping class group $\Ggn$. But the 
stabilizers of the boundary points are infinite, and at the boundary the 
topology of $\Tgnq$ is quite far from that of a manifold.\\

In his thesis \cite{VDiss}, V.~Braungardt gave a construction of $\Tgnq$ which
uses only the complex structure of $\Mgnq$ and the boundary divisor
$\partial\Mgn$. Moreover his construction endows $\Tgnq$ with the structure of a
locally ringed space and he shows that it is a fine moduli space for 
``marked'' stable Riemann surfaces.
In this chapter we give a brief account of his
approach.

\subsection{Coverings with cusps}
\label{v-cov}
The basic idea of Braungardt's construction\index{Braungardt's construction} 
is to study, for a complex manifold $S$,
quotient maps $W\to W/G = S$  that have ``cusps'' over a divisor $D$
 in $S$. This
concept, which will be explained in this section, generalizes the familiar 
ramified coverings. The key result is that, in
the appropriate category of such quotient maps, there exists a universal
 object $p:\tilde W\to S$ with cusps
over $D$.\\

In general $\tilde W$ cannot be a complex manifold or even a complex space. 
Therefore we have to work in the larger category of {\it locally complex 
ringed spaces}\index{locally ringed space}, i.\,e.\ topological spaces $W$ 
endowed with a sheaf ${\cal O}_W$ of $\CC$-algebras (called the {\it structure 
sheaf\index{structure sheaf}}) such that at each point $x\in W$ the stalk  
${\cal O}_{W,x}$ is a local $\CC$-algebra. The basic properties of such spaces 
can be found e.\,g.\ in \cite[Ch.\,1, \S\,1]{GR} (where they are called 
$\CC$-ringed spaces).\\

In our situation Braungardt constructs a normal locally
complex ringed space $\tilde W$ such that the subspace $\tilde W_0=\tilde W - p^{-1}(D)$ is a complex manifold
and the restriction $p|_{\tilde W_0}:\tilde W_0\to S_0=S - D$ is the usual 
universal
covering.
\begin{example}\label{horo} 
The simplest example is well known and quite typical: Take
$S$ to be the unit disk $\DD=\{z\in\CC:|z|<1\}$ and $D=\{0\}$. The universal
covering of $S - D$ is, of course, $\exp:\HH\to\DD-\{0\}$, $z\mapsto e^{2\pi
  iz}$. It turns out that the universal covering in Braungardt's sense is 
$\hat\HH=\HH\cup\{\mbox{i}\infty\}$ with the {\it horocycle topology}\index{horocycle topology}, 
i.\,e.\ the sets
$\HH_R=\{z\in\CC:\mbox{Im}\,z>R\}\cup\{\mbox{i}\infty\}$  
for $R>0$ form a basis of neighbourhoods of the point $\mbox{i}\infty$. Note
that this topology is not the one induced from the Euclidean topology if
$\HH\cup\{\infty\}$ is considered as a subset of the Riemann sphere
$\hat\CC$.\\[1mm]
$\hat\HH$ is given the structure of a normal complex ringed space by taking 
${\cal
  O}(U)$ to be the holomorphic functions on $U$ for open subsets $U$ of $\HH$,
and by defining
${\cal O}(\HH_R)$ to be the set of holomorphic functions on 
$\{z\in\CC:\mbox{Im}\,z>R\}$ that have a
continuous extension to $\mbox{i}\infty$. Clearly ${\cal O}(\HH_R)$ contains
all functions of the form $z\mapsto e^{2\pi iz/n}$ for all $n\ge1$.
\end{example}
We now give the precise definitions. We begin with the class of spaces that we 
need
(cf.~\cite[3.1.3]{VDiss}): 
\begin{definition}
\label{R-space}
Let $(W,{\cal O}_W)$ be a locally complex ringed space whose structure sheaf
 ${\cal O}_W$ is a subsheaf of the sheaf ${\cal C}^\infty(W,\CC)$ of continuous
 complex valued functions on $W$. \\[1mm]
{\bf a)} A subset $B\subset W$ is called {\it analytic} if there is an open
covering $(U_i)_{i\in I}$ of $W$ and for each $i\in I$ there are finitely many
elements $f_{i,1},\dots,f_{i,n_i}\in{\cal O}_W(U_i)$ such that $B\cap U_i$
is the zero set of $\{f_{i,1},\dots,f_{i,n_i}\}$.\\[1mm]
{\bf b)} We call $(W,{\cal O}_W)$ an {\it
R-space}\index{R-space} if, for every open $U\subseteq W$ and every proper 
closed analytic
subset $B\subset U$, a continuous function $f:U\to\CC$ is in ${\cal O}_W(U)$ if
and only if its restriction to $U-B$ is in ${\cal O}_W(U-B)$.
\end{definition}
Note that all complex spaces are R-spaces: The required property is just
Riemann's extension theorem\index{Riemann's extension theorem}, 
see \cite[Chapter 7]{GR}. 
\begin{definition}
\label{cusps}
Let $S$ be a complex manifold and $D\subset S$ a proper closed analytic
subset. Then a surjective morphism $p:W\to S$ from an R-space $(W,{\cal O}_W)$
to $S$ is called a 
{\it covering   with cusps over}\index{covering with cusps} $D$ if there is a 
group $G$ of 
automorphisms of $W$ (as locally complex ringed space) such that
\begin{itemize}
\item[(i)] $p$ is the quotient map $W\to W/G=S$,
\item[(ii)] $W_0=p^{-1}(S-D)$ is a complex manifold and $p|_{W_0}:W_0\to
  S_0=S-D$ is an unramified covering,  
\item[(iii)] for any $x\in W$ there is a basis of neighbourhoods $U_x$ that
  are {\it precisely invariant}\index{precisely invariant} under the stabilizer
 $G_x$ of $x$ in $G$ (i.\,e.\
  $G_x(U_x)=U_x$ and $g(U_x)\cap U_x=\emptyset$ for each $g\in G-G_x$).
\end{itemize}
\end{definition}
Note that, in particular, any ramified normal covering of complex manifolds is 
a covering in the sense of this definition (with cusps over the branch
locus). As mentioned before, the basic result is (see \cite[Satz 3.1.9]{VDiss})
\begin{theorem}
\label{univ}
{\em (i)} For any complex manifold $S$ and any proper closed analytic subset 
$D\subset
S$ there exists an initial object $p:(\tilde W,{\cal O}_{\tilde W})\to S$ in
the category of coverings of $S$ with cusps over $D$; it is
called the {\em universal covering  with cusps over}\index{covering with 
cusps!universal}
$D$. The restriction of $p$ to $\tilde W_0=p^{-1}(S_0)$ is the universal 
covering of $S_0$, and the group
$G=\mbox{\em Aut}(\tilde W/S)$ is the fundamental group $\pi_1(S_0)$.\\[1mm]
{\em (ii)} If $S'$ is an open submanifold of $S$ and $\tilde W'$ the universal
covering of $S'$ with cusps over $D'=D\cap S'$, then $\tilde W'/H'$ embeds as 
an open
subspace into $\tilde W$, where $H'$ is the kernel of the homomorphism 
$\pi_1(S'-D')\to\pi_1(S-D)=G$.
\end{theorem}
\begin{proof}
We only sketch the construction of the space $(\tilde W,{\cal O}_{\tilde W})$. 
The details that it 
satisfies all the required properties are worked out in \cite{VDipl}. For the
proof of (ii) we refer to \cite{VDiss}.\\
Let $S_0=S-D$, $G=\pi_1(S_0)$ and $p_0:W_0\to S_0$ the universal covering. 
$\tilde W$
is obtained from $W_0$ by ``filling in the holes above $D$'' in such a way
that the $G$-action extends from $W_0$ to $\tilde W$. More formally, the fibre 
$\tilde W_s$
of $\tilde W$ over any point $s\in S$ is constructed as follows: let $\mathfrak
U(s)$ be the set of open connected neighbourhoods of $s$ in $S$; for any $U\in
\mathfrak U(s)$ denote by $X(U)$ the set of connected components of
$p_0^{-1}(U)$. Then
\[ \tilde W_s = \{(x_U)_{U\in\mathfrak U(s)}:x_U\in X(U),\,x_U\cap x_{U'}\not=\emptyset\
\mbox{for all}\ U,U'\in\mathfrak U(s)\}.\]
Clearly $\tilde W_s=p_0^{-1}(s)$ for $s\in S_0$. Note that by definition, $G$ acts
transitively on each $\tilde W_s$, thus $\tilde W/G=S$. For any $x=(x_U)\in \tilde W$ define the
sets $x_U\cup\{x\}$, $U\in\mathfrak U(s)$, to be open neighbourhoods of
$x$. Finally define the structure sheaf by
\begin{align}
\label{calO}
{\cal O}_{\tilde W}(U) = \{f:U\to\CC\ \mbox{continuous}: f\ \mbox{holomorphic on}\
U\cap \tilde W_0\}
\end{align}
for any open subset $U$ of $\tilde W$.
\end{proof}
A key point in Braungardt's proof of Theorem \ref{univ} is the existence of
neighbourhoods $U$ for any point $a\in D$ such that the natural homomorphism
$\pi_1(U-D) \to \lim\limits_{\stackrel{\longrightarrow}{U'\in\mathfrak U(a)}}\,\pi_1(U'-D)$ is an isomorphism. He calls such
neighbourhoods {\it decent}\index{decent neighbourhood}.
The importance of this notion is that if $U$ is a decent neighbourhood of a
point $a\in D$ and $\bar x_U$ a connected component of $p^{-1}(U)$, then $\bar
x_U$ is 
precisely invariant under the stabilizer $G_x$ in $G$ of the unique point $x\in
\bar x_U\cap p^{-1}(a)$.\\

Decent neighbourhoods in the above sense do not exist in general for singular
complex spaces. For example, if $S$ is a stable Riemann surface and $s\in S$ a
node, $U-\{s\}$ is not even connected for small neighbourhoods $U$ of
$s$. Nevertheless the construction can be generalized to this case, and
the proof of the theorem carries over to this case as
Braungardt explains in \cite[Anm.\ 3.1.4]{VDiss}; we therefore have:
\begin{corollary}
\label{univstab}
Any stable Riemann surface\index{stable Riemann surface} has a universal covering with cusps over the nodes.
\end{corollary}
Near the inverse image of a node, the universal covering of a stable Riemann
surface looks like two copies of $\hat\HH$ glued together in the cusps. If such a
neighbourhood is embedded into the complex plane or $\PP^1(\CC)$ it is called a
{\it doubly cusped region}\index{doubly cusped region}, cf.\ \cite[VI.A.8]{Mask}.

\subsection{The cusped universal covering of $\Mgnq$}
\label{v-tgnq}
Let us now fix nonnegative integers $g$, $n$ such that $3g-3+n>0$. We want to 
construct the space $\Tgnq$ as the universal covering of $\Mgnq$ with cusps
over the compactification (or boundary) divisor $\dMgn$. But we cannot apply 
Theorem~\ref{univ} directly to $\Mgnq$ since it is not a
manifold, but only an orbifold (or smooth stack). Braungardt circumvents this
difficulty by
\begin{definition}
\label{covmgnq}
A morphism $p:Y\to\Mgnq$ of locally complex ringed spaces is called a {\it
  covering with cusps over} $D=\partial\Mgn$ if there is an open covering
  $(U_i)_{i\in I}$ of $\Mgnq$ and for each $i\in I$ a covering
  $q_i:U'_i\to U_i$ by a complex manifold $U'_i$ such that $p|_{p^{-1}(U_i)}$
  factors as
  $p^{-1}(U_i)\stackrel{\scriptsize{p'_i}}{\longrightarrow}U'_i\stackrel{\scriptsize{q_i}}{\longrightarrow}U_i$,
  where $p'_i$ is a covering with cusps over $q_i^{-1}(D)$ (in the sense of
  Definition \ref{cusps}).
\end{definition}
Then one can use Theorem \ref{univ} to prove
\begin{proposition}
\label{tgnq}
There is a universal covering $\Tgnq\to\Mgnq$\index{Teichm\"uller space!partial compactification $\Tgnq$}\index{moduli space!compactification $\Mgnq$} with cusps over
$\partial\Mgn$.
\end{proposition}
\begin{proof}
We first construct local universal coverings and then glue them together. 
For any $s\in\Mgnq$ choose an open neighbourhood $U$ and a covering
$q':U'\to U$ with a manifold $U'$. Let $\tilde W'$ be the universal covering
of $U'$ with cusps over $D'=q'^{-1}(D)$. Let $H'$ be the kernel of the
homomorphism from $\pi_1(U'-D')$ to $\Ggn$.
Theorem \ref{univ} (ii) suggests that the quotient $\tilde W'/H'$ 
should be an open part of the universal covering of $\Mgnq$. All that remains
to show is that the
$\tilde W'/H'$ glue together to a covering with cusps over $D$. This is done in \cite[3.2.1]{VDiss}
\end{proof}
Locally $\Tgnq$ looks like a product of a ball with some copies of the
universal covering $\hat\HH$ of $\DD$ with cusps over $\{0\}$ which 
was explained in Section \ref{v-cov}:
\begin{corollary}
\label{local}
Let $x\in\Tgnq$ correspond to a stable Riemann surface $X$ with $k$ nodes. Then $x$ has
a neighbourhood that is isomorphic to 
\[\hat\HH^k\times\DD^{3g-3+n-k}.\]
\end{corollary}
\begin{proof}
Let $s\in\Mgnq$ be the image point of $x$. The deformation theory of stable
Riemann surfaces gives us a map from $\DD^{3g-3+n}$ onto a neighbourhood of 
$s$ such that the inverse image of $D=\dMgn$ is the union of axes
$D'=\{(z_1,\dots,z_{3g-3+n}): z_1\cdot\ldots\cdot z_k = 0\}$, see 
\cite[Sect.~3B]{HM}.
The fundamental group of $\DD^{3g-3+n}-D'$ is a free abelian group on $k$
generators; they correspond to Dehn twists about the loops that are
contracted in $X$. Thus the homomorphism
$\pi_1(\DD^{3g-3+n}-D')\to\Ggn$ is injective. By Proposition~\ref{tgnq} and
its proof the universal covering $\tilde W$ of $\DD^{3g-3+n}$ with cusps over $D'$ is
therefore a neighbourhood of $x$. It is not hard to see that $\tilde W$ is of the
given form.
\end{proof}
Our next goal is to compare $\Tgnq$ to the {\it augmented} Teichm\"uller space
$\Tgnd$ introduced by Abikoff \cite{Abideg}. 
\begin{proposition}[cf.~\cite{VDiss}, Satz 3.4.2]
\label{homeo}
$\Tgnq$\index{Teichm\"uller space!partial compactification $\Tgnq$} is homeomorphic to the augmented Teichm\"uller space\index{augmented Teichm\"uller space} $\Tgnd$.
\end{proposition}
Before proving the proposition we summarize the definition and some properties
of $\Tgnd$: As a point set,
\begin{align}\label{tgndach}
 \Tgnd = \{(X,f): X\ \mbox{a stable Riemann surface of type}\ (g,n),\nonumber\\
f:\Cref\to X\ \mbox{a deformation}\}/\sim
\end{align}
As mentioned in the introduction, a deformation\index{deformation map} is a map that contracts some
disjoint loops on $\Cref$ to points (the nodes of $X$) and is a homeomorphism
otherwise. The equivalence relation is the same as for $\Tgn$:
$(X,f)\sim(X',f')$ if and only if there is a biholomorphic map $h:X\to X'$
such that $f'$ is homotopic to $h\circ f$.\\[1mm]
Abikoff puts a topology on $\Tgnd$ by
 defining neighbourhoods $\UKe$ of a point $(X,f)$
for a compact neighbourhood $V$ of the set of nodes in $X$ and
$\varepsilon>0$:
\begin{align}\label{uke}
\UKe = \{(X',f'):\exists\ \mbox{deformation}\ h:X'\to X,\ \ 
(1+\epsilon)\,\mbox{-quasiconformal}\nonumber\\ 
\mbox{on}\ h^{-1}(X-V),\ \mbox{such that}\ f\ \mbox{is homotopic to}\
 h\circ f'\}/\sim
\end{align}
 The action of the mapping class
group $\Ggn$ extends continuously to $\Tgnd$ (\cite[Thm.~4]{Abideg}), and the
orbit space $\Tgnd/\Ggn$ is $\Mgnq$ (as a topological space). 
\begin{proof}[Proof of Proposition \ref{homeo}]
Braungardt shows (see \cite[Hilfssatz 3.4.4]{VDiss}) that the stabilizer
of a point $(X,f)\in\Tgnd$ in $\Ggn$ is an extension of the free abelian group 
generated by the Dehn twists about the contracted loops by the holomorphic 
automorphism group Aut$(X)$ of $X$. For any $V$ and
$\varepsilon$, $\bigcap\limits_{\sigma\in\mbox{\scriptsize Aut}(X)}\sigma(\UKe)$ is
invariant under the stabilizer of $(X,f)$, and for sufficiently small $V$ and
$\varepsilon$, it is precisely invariant. Therefore the quotient map $\Tgnd\to
\Mgnq$ is a covering with cusps over $\dMgn$ in the sense of Definition 
\ref{covmgnq}, except that so far no structure
sheaf has been defined on $\Tgnd$. But this can be done in the same way as in
(\ref{calO}). 
The universal property of $\Tgnq$ then
yields a map $p:\Tgnq\to\Tgnd$ compatible with the action of $\Ggn$ on both
sides. To show that this map is an isomorphism we compare the stabilizers in
$\Ggn$ for the points in both spaces. For a point in $\Tgnd$ we just
described this stabilizer, and the proof of Corollary~\ref{local} shows that
for a corresponding point in $\Tgnq$ it is also an extension of $\ZZ^k$ by
Aut$(X)$. 
\end{proof}

\subsection{Teichm\"uller structures}
\label{v-structure}
In this section we explain how Braungardt extends the universal family\index{universal family}
of marked Riemann surfaces that is well known to exist over $\Tgn$ to a family
over $\Tgnq$ which still is universal for the appropriate notion of marking or
{\it Teich\-m\"uller structure}.\\

As above we fix a reference Riemann surface $\Cref$ of type $(g,n)$; let
$Q_1,\dots,Q_n$ be the marked points and $\Crefo = \Cref -
\{Q_1,\dots,Q_n\}$. Let us also fix a universal covering $\Uref\to\Crefo$ and
identify $\pgn=\pi_1(\Crefo)$ with the group Aut$(\Uref/\Crefo)$ of deck
transformations.\\

A classical construction of the family $\Cgn$ over $\Tgn$ goes as follows
(cf.~\cite{Bers}): For
every point $x=(X,P_1,\dots,P_n,f)\in\Tgn$ take a universal covering of $X^0 =
X-\{P_1,\dots,P_n\}$ and arrange them so that they form an $\HH$\,-bundle
$\Omega^+$ over $\Tgn$. Then $\Cgn$ is obtained as the quotient of 
$\Omega^+$ by the natural action of $\pgn$. More precisely, $\Omega^+$ is defined as follows: to $x\in\Tgn$ there corresponds the quasifuchsian group $G_x=w^\mu G(w^\mu)^{-1}$, where $G=\mbox{Aut}(\Uref/\Crefo)\cong\pgn$ and $w^\mu$ is the quasiconformal automorphism of $\PP^1(\CC)$ associated to $x$, see e.\,g.\ \cite[6.1.1]{IT}. The domain of discontinuity of $G_x$ consists of two connected components $\Omega^-(x)=w^\mu(L)$ (where $L$ is the lower half plane) and $\Omega^+(x)=w^\mu(\HH)$. Then $\Omega^+(x)/G_x=X^0$, whereas $\Omega^-(x)/G_x={\Crefos}$, the mirror image of $\Crefo$.\\

To extend this family we identify $\Tgnq$ with $\Tgnd$ by Corollary
\ref{homeo}. As explained in \cite{Abideg}, any point
$x=(X,P_1,\dots,P_n,f)\in\Tgnd-\Tgn$ corresponds to a {\it regular B-group}\index{regular B-group} 
$G_x$. This means that $G_x$ is a Kleinian group isomorphic to $\pgn$ whose domain of discontinuity $\Omega(G_x)$ has a
unique simply connected invariant component $\OminG$ such that $\OminG/G_x$ is
isomorphic to ${\Crefos}$. For the union $\OplusG=\Oplusx$ of the other 
components of $\Omega(G_x)$ it holds that $\OplusG/G_x\cong
X^0-\{\mbox{nodes}\}$. To every node in $X$ there corresponds a conjugacy
class of parabolic elements in $G_x$, each of which is
accidental\index{accidental parabolic element} (i.\,e.\ it becomes hyperbolic in the Fuchsian group $hG_xh^{-1}$,
where $h:\OminG\to\HH$ is a conformal map). Near a fixed point of such a
parabolic element, $\OplusG$ is a doubly cusped region\index{doubly cusped region}, cf.\ the remark at the end of
Section \ref{v-cov}. If we denote by $\Oplusxd$ the union of $\OplusG$ with
the fixed points of the parabolic elements in $G_x$ (accidental or not), then $\Oplusxd\to X$ is
the universal covering of $X$ with cusps over the nodes (cf.~Corollary
\ref{univstab}).  
\begin{definition}
\label{cgnq}
Let
\begin{align}
\Oplusgnd=\{(x,z)\in\Tgnq\times\PP^1(\CC):z\in\Oplusxd\}.\nonumber
\end{align}
On $\Oplusgnd$, $\pgn$ acts in such a way that for fixed $x\in\Tgnq$ the action on $\Oplusx$ is that of $G_x$. $\Cgnq=\Oplusgnd/\pgn$ is called the {\it universal family} over $\Tgnq$\index{universal family over $\Tgnq$}.
\end{definition}
Braungardt shows (\cite[Hilfssatz 4.2.1]{VDiss}) that
$\Oplusgn=\{(x,z)\in\Oplusgnd: x \in \tgn, \; z\in\Oplusx\}$ is an open subset of
$\Tgnq\times\PP^1(\CC)$ and hence has a well defined structure of a complex
ringed space. One can extend this structure sheaf to all of $\Oplusgnd$ in the
same way as in (\ref{calO}). Then clearly $\Cgnq$ is also a complex ringed
space, and the fibre over $x$ is isomorphic to the stable Riemann surface $X$
represented by $x$.\\

To justify the name ``universal'' family for $\Cgnq$ we introduces the notion
of a Teichm\"uller structure\index{Teichm\"uller structure}: For a single smooth Riemann surface $(X,P_1,\dots,P_n)$ of
type $(g,n)$, a Teichm\"uller structure is just a marking: so far in this article we used markings as classes of mappings $\Cref\to X$; equivalently a marking can be given as an isomorphism $\pgn\to\pi_1(X-\{P_1,\dots,P_n\})$ inducing an
isomorphism $\pi_g=\pi_1(\Cref)\to\pi_1(X)$ and respecting the orientation and the
conjugacy classes of the loops around the $Q_i$ resp.\ $P_i$. Yet another equivalent way to give a
marking is as a universal covering $U\to X^0$ together with an isomorphism
$\pgn\to \mbox{Aut}(U/X^0)$. This last characterization also works for a stable Riemann surface if
we take for $U$ a universal covering with cusps over the nodes. Before we
extend this definition to the relative situation we recall the notion of a
family of stable Riemann surfaces:
\begin{definition}
\label{famstab}
Let $S$ be a complex ringed space. A {\it family of stable Riemann surfaces}\index{stable Riemann surfaces!family of} of type $(g,n)$ over
$S$ is a complex ringed space ${\cal C}$ together with a proper flat map
$\pi:{\cal C}\to S$
such that the fibres $X_s=\pi^{-1}(s)$, $s\in S$, are stable Riemann surfaces
of genus $g$. In addition we are given $n$ disjoint sections $P_i:S\to {\cal
  C}$, $i=1,\dots,n$, of $\pi$ such that $P_i(s)$ is not a node on $X_s$. We denote by ${\cal
C}^0={\cal C}-\bigcup\limits_{i=1}^nP_i(S)$ the complement of the marked sections.
\end{definition}
\begin{definition}
\label{t-structure}
Let ${\cal C}/S$ be a family of stable Riemann surfaces of type $(g,n)$ over a complex ringed space $S$. A {\it
  Teichm\"uller structure} on ${\cal C}$ is a complex ringed space ${\cal U}$ together with a morphism ${\cal U}\to {\cal C}$ such
  that for every $s\in S$ the (restriction of the) fibre $U_s^0\to X_s^0$ is a universal covering
  with cusps over the nodes, together with an isomorphism
  $\pgn\to\mbox{Aut}({\cal U}/{\cal C}^0)$. 
\end{definition}
Putting everything together we obtain
\begin{theorem}
\label{fein}
$\Tgnq$ is a fine moduli space for stable Riemann surfaces with Teichm\"uller
structure. $\Cgnq\to\Tgnq$ is the universal family and $\Oplusgnd\to\Cgnq =
\Oplusgnd/\pgn$ is the universal Teichm\"uller structure.
\end{theorem}
Finally Braungardt gives a very elegant and conceptual description of $\Cgnq$
which extends a classical result of Bers (\cite[Thm.~9]{Bers}) to the boundary:
\begin{proposition}
\label{tgn+1}
$\Tgnpq/\pgn$ is in a natural way isomorphic to $\Cgnq$.
\end{proposition}
\begin{proof}
The kernel of the obvious homomorphism $\Ggnp\to\Ggn$ can be identified with
$\pgn$, which gives the action on $\Tgnpq$. The holomorphic map $\Tgnp\to\Tgn$
which forgets the last marked point extends to a map $\Tgnpq\to\Tgnq$ by a
general property of universal coverings with cusps. The difficult step in
Braungardt's proof is to show that the induced map $\Tgnpq/\pgn\to\Tgnq$ has the right fibres. For this
purpose he constructs a map $\Oplusgnd\to\Tgnpq$ and shows that it is
bijective and induces isomorphisms on the fibres over $\Tgnq$.
\end{proof}

\section{Boundary points of Teichm\"uller curves} 

The aim of this chapter is to study the boundary points
of the Teichm\"uller disks
and Teichm\"uller curves introduced in Chapter \ref{bp} in 
$\Tgq$ and $\Mgq$, respectively. 
Here and later, whenever we speak about $\tgquer$ and its boundary, we 
mean the bordification of the Teich\-m\"ul\-ler space
described in Chapter \ref{volker}.\\
In particular we will derive, in Section \ref{bpofdisks},
the following description   
of the boundary points of Teichm\"uller curves
(see Proposition~\ref{iquerfortc} and Corollary~\ref{finalcor}
for a more precise formulation): 
\begin{theorem}\label{thm}
One obtains the boundary points of a Teichm\"uller curve  
by contracting the centers of 
all cylinders in Strebel directions.
They are determined by the parabolic elements in 
the associated mirror Veech group.
\end{theorem}

This statement seems to be well known to the experts although
we are not aware of a published proof.\\
In Section \ref{srendpoint} we prepare for the proof of Theorem \ref{thm}
by introducing Strebel 
rays. They are special 
geodesic rays in Teichm\"uller space which always converge
to a point on the boundary. 
Following Masur \cite{M}, we describe 
this boundary point
quite explicitly using the affine structure of the quadratic differential $q$ that defines the Strebel ray.\\
In Section \ref{bpofdisks} we turn to the boundary points of 
Strebel rays that are contained in a Teichm\"uller disk.
In particular if the Teichm\"uller disk
descends to a Teichm\"uller curve in the moduli space, all its boundary points 
can be determined explicitly with the aid of the projective Veech group.
One obtains  Theorem \ref{thm} as a conclusion.

\subsection{Hitting the boundary via a Strebel ray}
\label{srendpoint}

In this section, we introduce Strebel rays and describe their end point
on the boundary of $\tgquer$. 
As before, everything might be done
as well for punctured surfaces and the moduli space $\tgn$ with $3g-3+n > 0$,
but for ease of notation, we restrict to the case $n=0$.\\

Let $X$ be a Riemann surface of genus $g \geq 2$, $q$ a
holomorphic quadratic differential on $X$.
Recall from Section \ref{deform} that 
with $q$ we have chosen a  
natural flat structure $\mu$ on the surface 
$X^* = X - \{\mbox{critical points  of $q$}\}$ whose charts were given
in (\ref{natchart}).
The maximal real curves in $\Xstern$
which are locally mapped by these charts to horizontal 
(resp. vertical) line segments 
are called {\em horizontal}\index{trajectory!horizontal}\index{trajectory} (resp. {\em vertical})\index{trajectory!vertical} {\em trajectories}. A 
trajectory is {\em critical}\index{trajectory!critical} if it ends in a critical point. Otherwise
it is {\em regular}\index{trajectory!regular}.

\begin{definition}
We say that a holomorphic quadratic differential $q$
is {\em Strebel}\index{Strebel differential},
if all regular horizontal trajectories are closed.
\end{definition}

Strebel differentials play an exceptional role 
in the following sense. Recall from Section \ref{sr} that 
each holomorphic quadratic differential defines a geodesic ray.
If $q$ is Strebel, then the geodesic ray defined by its negative $-q$ 
converges in $\tgquer$ to an end point on the boundary.
This is described more precisely in the following proposition which  
was proven by Masur in \cite{M}. We give a version of his proof 
with parts of the notation and arguments adapted to the context of 
our article.\\ 

Recall  also from \ref{sr} that we obtain the geodesic ray\index{geodesic ray} to $-q$ as the image 
of the isometric embedding
\begin{equation}\label{minusq}
\gamma = \gamma_{-q}: \left\{
  \begin{array}{lcl}
  [0,\infty) & \rightarrow & T_g\\
  t & \mapsto & (X_K,f_K) = [(X,\mu_{-q})\circ \bpm K&0\\0&1\epm, \id]
      \qquad \mbox{with } K = e^t .
   \end{array}\right.
\end{equation}
Note that here $(X_K,f_K)$ is the Teichm\"uller deformation of 
$X$ of dilatation $K$ with respect to $-q$. Furthermore, $\mu_{-q}$ is
the translation structure on $\Xstern$ defined by $-q$.

\begin{proposition} \label{rayminusq}
Suppose $q \neq 0$ is a Strebel differential.
For the geodesic ray defined by $\gamma_{-q}$ 
in $\tg$, one has:
\begin{itemize}
\item[a)] The ray converges towards a unique point $(X_{\infty},f_{\infty})$ 
on the boundary of 
  the Teich\-m\"uller space $\tg$.
\item[b)] One obtains this point by contracting the 
central lines of the horizontal cylinders defined by $q$ as is described in \ref{contract}.
\end{itemize}
\end{proposition}

\begin{definition}\label{sray}
In the previous Proposition, the geodesic ray 
defined by $-q$, i.~e. the image of $\gamma_{-q}$  in $\tg$,
is  called a {\it Strebel ray}\index{Strebel ray}.
\end{definition}

For the proof of Proposition~\ref{rayminusq} one may 
use two slightly different perspectives 
of the Strebel ray. 
They are described in \ref{patch}, \ref{stretch} and 
\ref{dan}, \ref{contract}.
In  \ref{endpoint} we describe the boundary point $(X_{\infty},f_{\infty})$.
In \ref{conv} we show that the Strebel ray in fact converges  
towards this point.\\

Throughout Section \ref{srendpoint}, we assume that the 
differential $q$  is Strebel.

\subsubsection[$X$ as patchwork of rectangles]{$X$ as patchwork of 
rectangles\\[2mm]}
\label{patch}

One may regard $X$ as a patchwork of rectangles\index{patchwork!of rectangles} in the complex plane, as is
described in the following.\\

Since $q$ is Strebel, the surface $X$, with the critical points and critical
horizontal trajectories removed, is swept out by closed horizontal trajectories.
More precisely, it follows from the work of Strebel (cf. \cite{st},
also see \cite[Theorem B]{M} which contains a list of the results
we use here) that
the surface $X$, except for the critical points and critical
horizontal trajectories, is covered by a finite number of  
maximal  {\em hori\-zontal cylinders}\index{cylinder!horizontal} $Z_1$, \ldots 
$Z_p$, i.~e. annuli that are swept out by closed horizontal trajectories. 
For each $Z_i$ one may choose a vertical trajectory $\beta_i$ 
joining opposite boundary components
of $Z_i$.
If we remove $\beta_i$ from $Z_i$, the remainder is mapped, by the natural 
chart $w_i$ defined by $\mu$ (see (\ref{natchart})), to an open rectangle $R_i$ 
in the complex plane.
The horizontal and vertical edges have lengths 
\[a_i = \int_{\alpha_i}|q(z)|^{\frac12}dz \;\; \mbox{ and } \;\;
  b_i = \int_{\beta_i}|q(z)|^{\frac12}dz,\]
where $\alpha_i$ is any closed horizontal trajectory in the cylinder $Z_i$.\\ 

One may extend $w_i^{-1}$ 
uniquely to a 
map from the closure $\bar{R}_i$ of  $R_i$
to the closure of the annulus $Z_i$. 
Then the two horizontal edges of $\bar{R}_i$ are
mapped to the two horizontal boundary components of $Z_i$ and
the two vertical edges are both mapped to $\beta_i$.
The critical points of $q$ 
that lie on the boundary of $Z_i$ define by their preimage
marked points on the horizontal edges  of $\bar{R}_i$
and decompose them into segments. \\

For each such segment $s$ on a horizontal edge 
of $\bar{R}_i$ its image on $X$ joins the annulus $Z_i$
to an annulus $Z_j$ possibly with $i=j$.\\
Thus the map $w_j\circ w_i^{-1}$ ($w_i^{-1}$ is the
extended map, $w_j$ is locally
the inverse map of the extended map $w_j^{-1}$)
is an {\em identification map}
between $s$ and a segment on a horizontal edge 
of $\bar{R}_j$. (Images of critical points have to be excluded.)\\
These identification maps are of the form $z \mapsto \pm z + c$
with a constant $c \in \CC$.\\

Conversely, given the closed rectangles
$\bar{R}_1$, \ldots, $\bar{R}_p$,
the marked points on their horizontal edges  
and these identification maps, we may recover the surface $X$ 
as follows: for each $i$ glue the two vertical edges of $\bar{R}_i$ by a 
translation and the horizontal edges (with the marked points removed) 
by the identification maps.
In this way, one obtains a surface $X^*$ with the flat structure
on it inherited from the euclidean plane $\CC$. 
By filling in the punctures at vertices, we obtain the original compact
Riemann surface $X$.\\

In this sense one may
consider $X$ as a patchwork 
of the rectangles
$\bar{R}_1$, \ldots, $\bar{R}_p$. This description
depends of course on the 
chosen holomorphic quadratic Strebel differential $q$.

\begin{example}{\em Two Riemann surfaces $X$ given as a patchwork of 
rectangles}:\label{examp}\\[2mm]
In the two examples in Figure \ref{l} and Figure
\ref{kk}, the two vertical edges of each rectangle are glued by
a translation, respectively.
Horizontal segments with the same name are glued. 
The direction of the arrow indicates whether the identification 
is a translation or a rotation by
$180\grad$. In the example in Figure \ref{l} 
one only has translations,
in the example in Figure \ref{kk} only rotations.\\[2mm]
In the first example the surface $X$ is of genus $2$ and all marked points
are identified and thus give only one point on $X$. In the second example
one obtains a surface of genus $0$ with four marked points
indicated by the four symbols $\bullet$, $\star$, $\circ$ and 
{\fs $\square$}.\\

\begin{minipage}[t]{6cm}
{\em Surface of genus 2 with 1 marked point:}\\
\setlength{\unitlength}{1cm}
\begin{picture}(6,4)
\put(0,0){\framebox(5,1.5){$R_1$}}
\put(0,2.5){\framebox(2,1){$R_2$}}
\put(-.1,-.1){\LARGE $\bullet$}
\put(1.9,-.1){\LARGE $\bullet$}
\put(.9,1.4){\LARGE $\bullet$}
\put(2.9,1.4){\LARGE $\bullet$}
\put(-.1,2.4){\LARGE $\bullet$}
\put(1.9,2.4){\LARGE $\bullet$}
\put(-.1,3.4){\LARGE $\bullet$}
\put(1.9,3.4){\LARGE $\bullet$}

\put(1,-.3){$s_1$}
\put(3.5,-.3){$s_2$}
\put(0.3,1.65){$s_2$}
\put(2,1.65){$s_3$}
\put(4,1.65){$s_2$}
\put(1,2.2){$s_3$}
\put(1,3.65){$s_1$}

\put(1.5,-.1){$>$}
\put(4,-.1){$>$}
\put(.6,1.42){$>$}
\put(2.6,1.42){$>$}
\put(4.4,1.42){$>$}
\put(1.5,2,4){$>$}
\put(1.5,3.42){$>$}
\put(-0.4, 1){$b_1$}
\put(4.3,-.5){$a_1$}
\put(-0.4, 3.1){$b_2$}
\put(1.7,2.2){$a_2$}
\end{picture}\\
\begin{center}
\refstepcounter{diagramm}{\it Figure \arabic{diagramm}}
\label{l}
\end{center}
\end{minipage}
\hspace*{5mm}
\begin{minipage}[t]{4cm}
{\em Surface of genus 0 with 4 marked points:}\\
\setlength{\unitlength}{1cm}
\begin{picture}(3,3.5)
\put(0,1.5){\framebox(4,1){$R_1$}}
\put(-.1,1.4){\LARGE $\bullet$}
\put(3.9,1.4){\LARGE $\bullet$}
\put(-.15,2.4){\LARGE $\star$}
\put(3.85,2.4){\LARGE $\star$}
\put(1.9,1.35){\LARGE $\circ$}
\put(1.75,2.37){ $\square$}
\put(1,1.2){$s_1$}
\put(3,1.2){$s_1$}
\put(1,2.65){$s_2$}
\put(3,2.65){$s_2$}
\put(1.4,1.4){$>$}
\put(2.55,1.4){$<$}
\put(1.4,2.45){$>$}
\put(2.55,2.45){$<$}
\put(-0.5, 2.2){$b_1$}
\put(3.55,1){$a_1$}
\end{picture}\\
\begin{center}
\refstepcounter{diagramm}{\it Figure \arabic{diagramm}}  
\label{kk}
\end{center}
\end{minipage}
\end{example}

\subsubsection[Stretching the cylinders]{Stretching 
the cylinders\\[2mm]} \label{stretch}

We will now redescribe the Strebel ray\index{Strebel ray} defined by  $-q$ 
by stretching the rectangles in the 'patchwork' from \ref{patch}
in the vertical direction.\\

The flat structure defined by  $-q = e^{\pi i}\cdot q$ is obtained from 
the flat structure $\mu$ defined by $q$
by composing each chart with a rotation by $\frac{\pi}{2}$.
Thus the deformation $(X_K,f_K)$ of dilatation $K$ with respect to $-q$
is equal to the affine deformation 
\[ \bpm K & 0\\ 0 & 1\epm \circ \bpm 0&-1\\1&0\epm \circ (X,\mu) 
  \;\; = \;\; \bpm 0 & -K\\ 1 & 0\epm \circ (X,\mu).
\]
This defines by (\ref{PU}) the same point in $\tg$ as the affine deformation
\[\bpm 0&1\\-1&0\epm \circ \bpm 0 & -K\\ 1 & 0\epm \circ (X,\mu)
  =   \bpm 1 & 0\\ 0 & K\epm \circ (X,\mu). \]
Thus the isometric embedding $\gamma = \gamma_{-q}$ in (\ref{minusq})
is equivalently given by
\begin{equation}\label{qminuszwei}
\gamma_{-q}:  \left\{
  \begin{array}{lcl}
  [0,\infty) & \rightarrow & T_g\\
  t & \mapsto & (X_K, f_K) \, = \, 
        [ \bpm 1 & 0\\ 0 & K\epm \circ (X,\mu),\;\; \id],\quad
   K = e^t 
  \end{array}\right.
\end{equation}

Recall again that here $(X_K,f_K) = (X_K^{-q},f_K^{-q})$ 
is the Teichm\"uller deformation with
respect to the differential $-q$.\\

Hence we obtain the point $\gamma_{-q}(t)$ as follows:
Each chart of $\mu$ is composed with the map 
$x + iy \mapsto x + iKy, \; (x,y \in \RR) $ with $K = e^t$,
and the marking is topologically the identity.
Now, let $X$ be given as a patchwork of the rectangles $\bar{R}_1$,
\ldots, $\bar{R}_p$ as in \ref{patch}. Then we obtain the surface
$X_K = X_K^{-q}$ in the following way:
We stretch each rectangle 
$\bar{R}_i$, which has horizontal and 
vertical edges of lengths $a_i$ and $b_i$, into a rectangle 
$\bar{R}_i(K)$ with horizontal 
and vertical edges of lengths $a_i$ and $K\cdot b_i$. The identification
maps of the horizontal segments are again translations or rotations 
identifying the same segments as before. The  surface
$X_K = X^{-q}_K$ then is the patchwork obtained from  $\bar{R}_1(K)$, \ldots,
$\bar{R}_p(K)$ as described in \ref{patch}.\\
On $\bar{R_i}$, the diffeomorphism $f_K = f^{-q}_K$ 
has image $\bar{R}_i(K)$ and is given by
\[x+iy \;\; \mapsto \;\; x + iKy.  \]
This glues to a well defined diffeomorphism on $X^*$,
which can be uniquely extended to $X$.

\begin{example} {\em $K$-stretched  surfaces:}\\
\hspace*{7mm}
\begin{minipage}[t]{6cm}
\setlength{\unitlength}{1cm}
\begin{picture}(6,6.7)
\put(0,0){\framebox(5,3){$R_1$}}
\put(0,4){\framebox(2,2){$R_2$}}
\put(-.1,-.1){\LARGE $\bullet$}
\put(1.9,-.1){\LARGE $\bullet$}
\put(.9,2.9){\LARGE $\bullet$}
\put(2.9,2.9){\LARGE $\bullet$}
\put(-.1,3.9){\LARGE $\bullet$}
\put(1.9,3.9){\LARGE $\bullet$}
\put(-.1,5.9){\LARGE $\bullet$}
\put(1.9,5.9){\LARGE $\bullet$}

\put(1,-.3){$s_1$}
\put(3.2,-.3){$s_2$}
\put(0.3,3.15){$s_2$}
\put(2,3.15){$s_3$}
\put(4,3.15){$s_2$}
\put(1,3.8){$s_3$}
\put(1,6.15){$s_1$}

\put(1.5,-.1){$>$}
\put(4,-.1){$>$}
\put(.6,2.92){$>$}
\put(2.6,2.92){$>$}
\put(4.4,2.92){$>$}
\put(1.5,3.9){$>$}
\put(1.5,5.92){$>$}
\put(-0.8, 2){$Kb_1$}
\put(4.3,-.55){$a_1$}
\put(-0.8, 5.5){$Kb_2$}
\put(1.7,3.6){$a_2$}
\end{picture}\\
\begin{center}
\refstepcounter{diagramm}{\it Figure \arabic{diagramm}}
\label{Kl}
\end{center}
\end{minipage}
\begin{minipage}[t]{4.5cm}
\setlength{\unitlength}{1cm}
\hspace*{5mm}
\begin{picture}(3,5.7)
\put(0,3.5){\framebox(4,2){$R_1$}}
\put(-.1,3.4){\LARGE $\bullet$}
\put(3.9,3.4){\LARGE $\bullet$}
\put(-.15,5.4){\LARGE $\star$}
\put(3.85,5.4){\LARGE $\star$}
\put(1.9,3.35){\LARGE $\circ$}
\put(1.75,5.37){ $\square$}
\put(1,3.2){$s_1$}
\put(3,3.2){$s_1$}
\put(1,5.65){$s_2$}
\put(3,5.65){$s_2$}
\put(1.4,3.4){$>$}
\put(2.55,3.4){$<$}
\put(1.4,5.45){$>$}
\put(2.55,5.45){$<$}
\put(-0.8, 5){$Kb_1$}
\put(3.5,3){$a_1$}
\end{picture}\\
\begin{center}
\refstepcounter{diagramm}{\it Figure \arabic{diagramm}}  
\label{Kkk}
\end{center}
\end{minipage}\\
One obtains the surface $X_K = X^{-q}_K$
from the surface $X$ in Example \ref{examp}
as the patchwork of 
the stretched rectangles in Figure \ref{Kl} and Figure \ref{Kkk}, 
respectively.
\end{example}

\subsubsection[$S$ as patchwork of double annuli]{$S$ as patchwork 
of double annuli\index{patchwork!of double annuli}\\[2mm]}
\label{dan}

Recall that, in \ref{patch}, we used $\mu$ to identify the horizontal 
cylinder $Z_i$ on 
$X$  with the euclidean cylinder defined by 
the rectangle $R_i$ in $\CC$; we did so by adding the vertical boundary edges and 
identifying them by a translation.
It turns out to be easier to describe the end point of the 
Strebel ray, if we identify the $Z_i$ with so called double annuli 
$A_i$.

\begin{definition} \label{dani}
A cylinder $Z$ of length $a$ and height $b$ 
defines  a {\em double annulus}\index{double annulus} $A$ as follows:
\begin{itemize}
\item Take two disjoint open annuli $A^1$ and $A^2$ given as
\[A^1 = A^2 = \{z \in \CC |\, r \leq |z| < 1\} 
    \;\;\mbox{ with } r = e^{-\pi\frac{b}{a}}.\]
\item Glue their inner boundary lines $\{|z| = r\}$ 
by the map  $z \mapsto \frac{1}{z}\cdot r^2$.
\item  We call the resulting surface $A$ the {\em double annulus}
of $Z$.
\end{itemize}
\end{definition}
\begin{remark}
$A$ is biholomorphic to $Z$. 
\end{remark}
The identification is given explicitly as follows:
\begin{itemize} 
\item $Z$ is biholomorphic
 to the Euclidean cylinder defined by the rectangle 
\[\{z \in \CC|\, 0 \leq \re(z) \leq a,\;\; 0 < \im(z) < b\}.\]
\item 
Decompose the rectangle into two halves of height $\frac{b}{2}$,\\
a lower half 
$R^1 = \{z \in \CC|\;\; 0 \, \leq \,  \re(z) \, \leq \, a,\;\;
                         0 \, < \, \im(z) \, \leq \, \frac{b}{2}\} $\\
and  an upper  half
$R^2 = \{z \in \CC|\;\; 0 \, \leq \,  \re(z) \, \leq  \, a,\;\; 
                         \frac{b}{2} \, \leq \, \im(z) \, < \,  b\}$.
\item 
The cylinder defined by
$R^1$ is mapped  to $A^1$
by \; $z \; \mapsto \; e^{2\pi i \frac{z}{a}}$.\\[1mm]
The cylinder defined by 
$R^2$ is mapped to $A^2$ 
by $z \; \mapsto \; e^{2\pi i \frac{a+bi-z}{a}}$.\\[1mm]
These maps respect the identifications and define a biholomorphic
map from $Z$ to $A$, as shown in Figure \ref{abab}.
\end{itemize}

\vspace*{-.5mm}

\begin{minipage}{\linewidth}
\begin{center}
\setlength{\unitlength}{1cm}
\begin{picture}(11,3.5)
\put(0,0.5){\framebox(4,2)}
\qbezier[40](0,1.5)(2,1.5)(4,1.5)
\put(1.8,.8){$R^1$}
\put(1.8,1.9){$R^2$}
\put(-.2,.2){$0$}
\put(3.8,.2){$a$}
\put(-.45,2.3){$bi$}
\put(-.45,1.35){$\frac{b}{2}i$}
\put(9,.7){\circle{1.3}}
\put(9,.7){\circle{.6}}
\put(8.9,.1){$A_1$}
\put(9,2.5){\circle{1.3}}
\put(9,2.5){\circle{.6}}
\put(8.82,2.87){$A_2$}
\put(9,1){\vector(0,1){1.23}}
\put(9,1.15){\vector(0,-1){.15}}
\put(9.9,2.1){$z$}
\put(9.8,1){$\begin{array}{l}
                e^{-2\pi \frac{b}{a}}\cdot\frac{1}{z}
               \end{array}$}
\put(9.9,1.89){\line(1,0){.2}}
\put(10,1.9){\vector(0,-1){0.7}}

\qbezier(3.5,2)(5.5,3.5)(8.5,2.7)
\put(8.3,2.75){\vector(4,-1){.3}}
\put(5.5,3.1){$z \; \mapsto \; e^{2\pi i \frac{a+bi-z}{a}}$}
\qbezier(3.5,1)(5.5,.5)(8.5,.7)
\put(8.5,.7){\vector(1,0){.1}}
\put(5,.8){$ z \mapsto
             \hspace*{5mm} e^{2\pi i \frac{z}{a}}$}
\qbezier[15](9,.7)(9.3,.7)(9.6,.7)
\put(9.2,.5){$r$}
\put(9.7,.5){$1$}
\end{picture}
\end{center}
\begin{center}
\refstepcounter{diagramm}{\it Figure \arabic{diagramm}}
\label{abab}
\end{center}
\end{minipage}

Consider the double annuli $A_1,\ldots ,A_p$
defined by the cylinders $Z_1, \ldots, Z_p$. 
The biholomorphic map $Z_i \to A_i$ extends to a continuous 
map from the closure of 
$Z_i$ to the closure $\bar{A}_i$ of $A_i$. The  
zeroes of $q$ on the boundary of $Z_i$ define marked points on the 
boundary of $A_i$ and decompose it into segments. 
The surface $X$ can now be described as a
{\em patchwork of the closed double cylinders} 
$\bar{A}_1$, \ldots, $\bar{A}_p$. 
The  identification maps between the segments on the boundary of the
$A_i$ are essentially the same as in \ref{patch}.

\subsubsection[Contracting the central lines]{Contracting the central lines\\[2mm]}
\label{contract}
Suppose that $X$ is given as a patchwork of double annuli 
$\bar{A}_1$, \ldots, 
$\bar{A}_p$ as in \ref{dan}. We may describe 
the points $(X_K,f_K) = (X_K^{-q},f_K^{-q})$
on the Strebel ray\index{Strebel ray} also as a patchwork of double annuli:\\
Let  $A_i(K) = A_i^1(K) \cup A_i^2(K)$ $(i \in \{1, \ldots, p\})$ be the double annulus from  
Definition~\ref{dani} with  $r = r_i(K) = r_i^K$
and define $X_K = X_K^{-q}$ to be the surface 
obtained by gluing the closures 
$\bar{A}_1(K)$, \ldots,  $\bar{A}_p(K)$ with the same 
maps as $\bar{A}_1$, \ldots,
$\bar{A}_p$.
Furthermore, define the diffeomorphism $f_K = f_K^{-q}$ on $A_i$
by 
\begin{eqnarray*}
&f_K^{-q}:\;  A_i^1 \to A_i^1(K) \; \mbox { and } \; A_i^2 \to A_i^2(K),&\\ 
&\hspace*{8mm} 
   z=r\cdot e^{i\varphi} \; \mapsto \; r^K\cdot e^{i\varphi} \quad  
     \mbox{on both parts.}
\end{eqnarray*}
Then the following diagram is commutative:\\[1mm]
\begin{minipage}{\linewidth}
\setlength{\unitlength}{1cm}
\begin{picture}(15,7.5)
\put(1,1.5){\circle{2}}
\put(.65,1.02){$\sscs A_i^2(K)$}
\put(1,1.5){\circle{.4}}
\qbezier(1.2,1.5)(2.25,1.7)(3.3,1.5)
\put(3.02,1.57){\vector(4,-1){.3}}
\put(1.85,1.75){$\sscs z \mapsto \frac{\scs r_i^{\sscs 2K}}{\scs z}$}
\put(3.5,1.5){\circle{2}}
\put(3.5,1.5){\circle{.4}}
\put(3.15,1.02){$\sscs A_i^1(K)$}
\put(1,4.4){\vector(0,-1){1.8}}
\put(2.1,3.8){$f_K$}
\put(3.5,4.4){\vector(0,-1){1.8}}

\put(9,0){\framebox(2,3)}
\qbezier[40](9,1.5)(10,1.5)(11,1.5)
\put(9.8,.5){$R_i^1(K)$}
\put(9.8,2.1){$R_i^2(K)$}
\put(9.5,4.6){\vector(0,-1){1.2}}
\put(9.6,3.8){$f_K$}
\put(10.8,-.3){$a_i$}
\put(8.5,3.1){$Kb_i$}

\qbezier[100](9.5,2.26)(6,5)(1.3,2)
\put(1.5,2.05){\vector(-4,-1){.3}}
\put(5.2,3.7){$z \mapsto e^{2\pi i\frac{a_i + Kb_ii-z}{a_i}}$}
\qbezier[70](9.3,1)(6,.7)(3.99,1.3)
\put(4.25,1.23){\vector(-4,1){.3}}
\put(5.5,1.13){$z \mapsto e^{2\pi i\frac{z}{a_i}}$}

\put(1,5.5){\circle{2}}
\put(1,5.5){\circle{.7}}
\put(.8,4.92){$\sscs A_i^2$}
\qbezier(1.35,5.5)(2.25,5.7)(3.142,5.5)
\put(2.865,5.57){\vector(4,-1){.3}}
\put(1.85,5.75){$\sscs z \mapsto \frac{\scs r_i^2}{\scs z}$}
\put(3.5,5.5){\circle{2}}
\put(3.5,5.5){\circle{.7}}
\put(3.3,4.92){$\sscs A_i^1$}

\put(9,5){\framebox(2,1)}
\qbezier[40](9,5.5)(10,5.5)(11,5.5)
\put(9.7,5.65){$R_i^2$}
\put(9.7,5.12){$R_i^1$}
\put(10.8,4.7){$a_i$}
\put(8.8,6.15){$b_i$}

\qbezier[100](9.3,5.7)(6,8)(1.3,6)
\put(1.5,6.05){\vector(-4,-1){.3}}
\put(5.8,7){$z \mapsto e^{2\pi i\frac{a_i+b_i i - z}{a_i}}$}
\qbezier[70](9.3,5.3)(6,4.7)(3.99,5.3)
\put(4.25,5.23){\vector(-4,1){.3}}
\put(5.5,5.2){$z \mapsto e^{2\pi i\frac{z}{a_i}}$}

\end{picture}
\begin{center}
\refstepcounter{diagramm}{\it Figure \arabic{diagramm}}
\end{center}
\end{minipage}\\[3mm]
where \; $f_K = (re^{i\varphi} \mapsto r^Ke^{i\varphi})$ \;
on the left side and \; $f_K = (x+iy \mapsto x+Kiy)$ \;
on the right side of the diagram.
Thus, in particular, we have defined here with $(X_K,f_K) = (X_K^{-q},f_K^{-q})$
the same surface (up to isomorphism) and the same diffeomorphism
as in \ref{stretch}.

\subsubsection[The end point of the Strebel ray]{The end point of the Strebel ray\\[2mm]}
\label{endpoint}
We use the description of the Strebel ray in \ref{contract} to 
obtain its end point\index{Strebel ray!end point of} $(X_{\infty}, f_{\infty}) \in \tgquer$. 
Recall from \ref{v-tgnq} that a point in $\tgquer$
consists of a stable Riemann surface $X_{\infty}$
and a deformation $f_{\infty}: X \to X_{\infty}$.\\ 

If $K \to \infty$ 
in \ref{contract}, the
interior radius $r_i(K) = r_i^K$ of the two annuli $A_i^1(K)$ and $A_i^2(K)$
that form the double annulus $A_i(K)$ 
tends to $0$ ($i \in \{1,\ldots, p\}$). $A_i(K)$ tends to a double 
cone $A_i(\infty)$ and the whole surface $X_K$ to a stable
Riemann surface $X_{\infty}$.
More precisely, we define $A_i(\infty)$ and $X_{\infty}$
as complex spaces in the following way.
\begin{definition}
Let $A_i^1(\infty)$ and $A_i^2(\infty)$  both  be
the punctured disk
\[\{z \in \CC|\, 0<|z| < 1\},\]
and let $\pt$ be an arbitrary point.
The disjoint union 
\[A_i(\infty) = A_i^1(\infty) \cup A_i^2(\infty) \cup \{\pt\}\]
becomes a complex cone by the following chart:
\begin{eqnarray*}
&&\varphi:\,\; A_i(\infty) \;\; \to \;\;  \{(z_1,z_2) \in \CC^2|\, z_1\cdot z_2 = 0,
                                        |z_1|, |z_2| < 1\}\\
&&\varphi|_{A_i^1(\infty)}:\, z \mapsto (0,z),\; \quad 
  \varphi|_{A_i^2(\infty)}:\, z \mapsto (z,0), \; \quad
  \varphi(\pt) = (0,0)
\end{eqnarray*}
The closures of the double cones $\bar{A}_1(\infty)$, \ldots, 
$\bar{A}_p(\infty)$  are glued to each other
by the same identification maps as in the 'finite' case in 
\ref{contract}. We call the resulting stable Riemann surface 
$X_{\infty}$. Topologically, $X_{\infty}$ is obtained from the surface $X$ 
by a contraction $f_{\infty}$ of  the middle curves of the cylinders.
\end{definition}

\begin{minipage}{\linewidth}
\begin{center}
\setlength{\unitlength}{.73cm}
\begin{picture}(6.5,6.5)
\put(1,1.5){\circle{2}}
\put(1,1.5){\circle*{.08}}
\put(.65,1.02){$\sscs A_i^2(\infty)$}
\put(.8,1.65){$\pt$}
\qbezier(1.15,1.53)(3.5,2)(5.85,1.53)
\put(5.7,1.57){\vector(4,-1){.3}}
\put(1.3,1.57){\vector(-4,-1){.3}}
\put(6,1.5){\circle{2}}
\put(6,1.5){\circle*{.08}}
\put(5.65,1.02){$\sscs A_i^1(\infty)$}
\put(5.8,1.65){$\pt$}
\put(1,4.4){\vector(0,-1){1.8}}
\put(1.1,3.5){$\begin{array}{l} 
                f_{\infty}:\\ 
                r\cdot e^{i\varphi} \, \mapsto \,
                h_{i,\infty}(r)\cdot e^{i\varphi}
               \end{array}$}
\put(6,4.4){\vector(0,-1){1.8}}

\put(1,5.5){\circle{2}}
\put(1,5.5){\circle{.7}}
\put(.8,4.92){$\sscs A_i^2$}
\qbezier(1.35,5.5)(3.25,6)(5.642,5.5)
\put(1.65,5.57){\vector(-4,-1){.3}}
\put(5.365,5.57){\vector(4,-1){.3}}
\put(2.5,6){$z \mapsto \frac{r_i^2}{z}$}
\put(6,5.5){\circle{2}}
\put(6,5.5){\circle{.7}}
\put(5.8,4.92){$\sscs A_i^1$}
\end{picture}
\end{center}
\begin{center}
\refstepcounter{diagramm}{\it Figure \arabic{diagramm}}
\end{center}
\end{minipage}\\

\noindent
We now define the contraction $f_{\infty}$ 
as the following map:
Let $A_i^1$ and $A_i^2$ be the two annuli in Definition \ref{dani}
that form the double annulus $A_i$ ($i\in\{1,\ldots,p\}$). 
Then $f_{\infty}$ is 
given by
\[
\begin{array}{lllcl}
   f_{\infty}:&  A_i^j                 &\to& 
       A_i^j(\infty) &  
       \mbox{ for } j \in \{1,2\}\\
              &z = r\cdot e^{i\varphi} &\mapsto& 
       f_{\infty}(z) = h_{i,\infty}(r)\cdot e^{i\varphi}& 
\end{array}
\]
with an arbitrary monotonously increasing diffeomorphism $
h_{i,\infty}: [r_i,1) \to [0,1)$. The isotopy class  of $f_{\infty}$ is
independent of the choices of $h_{i,\infty}$.\\

\subsubsection[Convergence]{Convergence\\[2mm]}
\label{conv}
We now show that, in the above notation, the Strebel ray $\gamma_{-q}$ 
converges to the point 
$(X_{\infty},f_{\infty})$
on the boundary of $T_g$. Recall from
(\ref{uke}) in Chapter \ref{volker} that
a base of open neighbourhoods of $(X_{\infty},f_{\infty})$ is given
by the open sets
\[U_{V,\varepsilon}(X_{\infty},f_{\infty}) = \{(X',f')| 
     \begin{array}[t]{l}
        \exists \; \varphi: X' \to X_{\infty}
        \mbox{ s.t. }
         \varphi \mbox{ is deformation},  \\
         \varphi\circ f' \mbox{ is isotopic to } f_{\infty}  \mbox{ and }\\
         \varphi|_{X' \backslash\varphi^{-1}(V)}
              \mbox{ has dilatation } < 1+\varepsilon
         \},
     \end{array}
\]
for all compact neighbourhoods $V$ of the singular points of $X_{\infty}$
and for all $\varepsilon > 0$.
We may restrict to open neighbourhoods $V$ of the form
\[V = V(\kappa) = V_1 \cup \ldots \cup V_p, \quad
\kappa = (\kappa_1, \ldots, \kappa_p),\;\; 0< \kappa_i < 1\]
where $V_i$ is a double cone defined by
\[V_i = V_i^1 \, \cup \, V_i^2 \, \cup \, \{\pt\} \; \mbox{ with } 
  V_i^j = \{0 < |z| \leq \kappa_i \} \, \subseteq \,  A_i^j(\infty) 
    \quad (j \in \{1,2\} ).\]
\begin{lemma}\label{lemconv}
For each such $V = V(\kappa)$ 
and each $\varepsilon > 0$, there is some $K_0 \in \RR_{>0}$ 
such that all points $(X_K,f_K) = (X_K^{-q},f_K^{-q})$ with $K > K_0$
are in $U_{V,\varepsilon}(X_{\infty},f_{\infty})$.
\end{lemma}

\begin{proof}
Choose $K_0$ such that $r_i^{K_0} < \kappa_i$ for all $i \in \{1,\ldots, p\}$ 
and suppose that $K > K_0$. Define the diffeomorphism 
$\varphi: X_K \to X_{\infty}$ on $\bar{A}_i^j(K)$ by
\[\varphi:\; z = r\cdot e^{i\theta} \mapsto 
    \left\{
    \begin{array}{ll}
     z \;\, \in \, A_i^j(\infty), \; &\mbox{ if } 1 >   |z| \geq \kappa_i\\
     h^i_K(r)\cdot e^{i\theta}  \;\, \in \, A_i^j(\infty), \;
                        &\mbox{ if } \kappa_i \geq   |z| > r_i^K\\ 
     \pt \;\, \in \, A_i^j(\infty), \; &\mbox{ if } |z| = r_i^K
    \end{array} \right.
\]
with an arbitrary monotonously increasing diffeomorphism 
$h^i_K: (r_i^K,\kappa_i) \to (0,\kappa_i)$.
Then $\varphi\circ f_K$ is isotopic to
$f_{\infty}$  and 
$\varphi|_{X_K \backslash \varphi^{-1}(V)}$ is holomorphic, hence  its
dilatation is
$1$. Thus $(X_K, f_K)$ is in $U_{V,\varepsilon}(X_{\infty},f_{\infty})$.
\begin{center}
\setlength{\unitlength}{1cm}
\begin{picture}(5,3)

\put(.5,2.1){$A_i^j(K)$}
\put(1,1){\circle{4}}
\put(1,1){\circle{.4}}
\qbezier[10](.55,1)(.55,1.45)(1,1.45)
\qbezier[10](1,1.45)(1.45,1.45)(1.45,1)
\qbezier[10](1.45,1)(1.45,.55)(1,.55)
\qbezier[10](1,.55)(.55,.55)(.55,1)

\qbezier(1,0)(.9,0)(0.35,0)
\put(.82,0.05){\line(0,-1){.15}}
\put(.95,-.5){$r_i^K$}
\put(.55,0.05){\line(0,-1){.15}}
\put(.5,-.5){$\kappa_i$}
\put(.35,0.05){\line(0,-1){.15}}
\put(.2,-.5){$1$}

\put(5.6,2.1){$A_i^j(\infty)$}
\put(6,1){\circle{4}}
\put(6,1){\circle*{.03}}
\qbezier[10](5.55,1)(5.55,1.45)(6,1.45)
\qbezier[10](6,1.45)(6.45,1.45)(6.45,1)
\qbezier[10](6.45,1)(6.45,.55)(6,.55)
\qbezier[10](6,.55)(5.55,.55)(5.55,1)

\qbezier(6,0)(5.9,0)(5.35,0)
\put(5.55,0.05){\line(0,-1){.15}}
\put(5.55,-.5){$\kappa_i$}
\put(5.35,0.05){\line(0,-1){.15}}
\put(5.2,-.5){$1$}

\qbezier(1,1.6)(3.5,2.5)(6,1.6)
\put(5.8,1.7){\vector(2,-1){.3}}
\put(3.5,2.15){id}

\qbezier(1.27,1)(4,1)(5.8,1)
\put(5.8,1){\vector(1,0){.03}}
\put(2.2,1.1){$ re^{i\theta} \mapsto h^i_K(r)e^{\sscs i\theta}$}
\end{picture}
\end{center}
\begin{center}
\refstepcounter{diagramm}{\it Figure \arabic{diagramm}}
\end{center}
\end{proof}
With Lemma \ref{lemconv} we have obtained the desired result 
and completed the proof of  
Proposition~\ref{rayminusq}.

\begin{corollary}
The Strebel ray\index{Strebel ray} defined by $-q$ converges to the point 
$(X_{\infty},f_{\infty})$ on the boundary of $\tg$.
\end{corollary}

\subsection{Boundary points of Teichm\"uller disks}\label{bpofdisks}
\index{Teichm\"uller disk!boundary points of}

In this section we study the boundary points of a Teichm\"uller disk 
$\Delta = \Delta_{\iota}$ in the 
bordification $\Tgq$ of the Teichm\"uller space; in particular, we consider
the case that $\Delta_{\iota}$ projects to an affine curve in the moduli 
space $\mg$. For convenience, we use the upper half plane model
and consider Teichm\"uller embeddings as maps from $\HH$ to $\tg$.
We will obtain Theorem \ref{thm} as our final result. We proceed
in two steps:
\begin{itemize}
\item In \ref{iotaquer}, we show that a Teichm\"uller embedding 
$\iota:\HH \hookrightarrow \tg$ 
has a natural extension 
\[\bar{\iota}: \HH \cup \{\mbox{cusps of $\Gammaquerm_{\iota}$}\} 
     \hookrightarrow 
    \Tgq,\]
\item 
In \ref{igitter}, we show that the image of 
$\bar{\iota}$  is the whole closure of $\Delta_{\iota}$ in $\tgquer$,
if the Teichm\"uller disk $\Delta_{\iota}$ 
projects onto a Teichm\"uller curve in $\mg$.\\  
It will follow from this that one obtains the boundary points of $\Delta_{\iota}$
precisely by contracting the central lines of the cylinders in ``parabolic 
directions''. The parabolic directions correspond to the cusps of the 
projective mirror Veech group $\Gammaquerm_{\iota}$.
\end{itemize}

Throughout this section, we assume that 
$\iota: \HH \into \tg$ is a Teichm\"uller embedding to
a fixed holomorphic quadratic differential $q$ on $X = \Xref$ and that
$\mu$ is the translation structure defined by $q$ as in
Section~\ref{deform}. Recall from Section \ref{tc} that the associated
projective Veech group $\Gammaquer_{\iota} = \Gammaquer(X,\mu)$
and its mirror image $\Gammaquerm_{\iota} = R\Gammaquer_{\iota} R^{-1}$ (with 
$R$ as in Remark \ref{achtionstogether}) are both Fuchsian groups
in $\pslzwei(\RR)$.

\subsubsection[Extending Teichm\"uller embeddings to the cusps of the
mirror Veech group]{Extending Teichm\"uller embeddings to the cusps of 
$\Gammaquerm$\\[2mm]}
\label{iotaquer}

Let $\tilde{s}\in \RR^{\infty} = \RR \cup \{\infty\}$ be a cusp of the 
Fuchsian group $\Gammaquerm_{\iota}$, i.e. $\tilde{s}$ is  a fixed point of 
some parabolic element $\tilde{A}$ of $\Gammaquerm_{\iota}$. We associate
to $\tilde{s}$ a point 
$\iotaquer(\tilde{s}) = (X_{\infty}(\tilde{s}),f_{\infty}(\tilde{s}))$ 
on the boundary of $\tg$ in the following way:
\begin{itemize}
\item In a natural way we associate to $\tilde{s}$ a Strebel ray\index{Strebel ray}
\item We show that this Strebel ray 
is the image in $\tg$   of the hyperbolic ray in $\HH$ from 
$i$ to $\tilde{s}$  under $\iota$.
\item
$\iotaquer(\tilde{s}) = (S_{\infty}(s),f_{\infty}(s))$ is defined to 
be the end point of the 
Strebel ray 
\end{itemize}

\noindent
{\em The Strebel ray associated to $\tilde{s}$:}\\[1mm]
$A = R^{-1}\tilde{A}R$ is a parabolic element in the projective Veech group 
$\Gammaquer_{\iota}$.
Let $v$ be its unit eigenvector. \\
By Proposition 2.4 in \cite{V},
the direction $v$ is fixed by some affine diffeomorphism $h$ of $(X,\mu)$.
The  
derivative of $h$ is $A$ and $v$ is a  Strebel direction. 
More precisely: The
trajectories in the direction of $v$ are preserved by $h$ and each leaf is 
either closed or a saddle connection\index{saddle connection}, i.~e. connects two 
critical points.\\
As in  \ref{patch}, $X$ decomposes into maximal cylinders\index{cylinder!maximal} 
of closed leaves parallel to $v$ and the cylinders are bounded by 
saddle 
connections. The affine diffeomorphism $h$ can be described nicely 
as follows: Passing to a power of $h$ if necessary, one may assume that
$h$ fixes all critical points of $q$. Then $h$ is the composition of 
Dehn twists along the core curves of the cylinders. Each trajectory is
mapped by $h$ to itself and the saddle connections are fixed pointwise.\\
Now, let us take the matrix 
\[U \; = \;  U_{\theta} \; \in \, \sozwei(\RR) \;\; \mbox{ such that } \;
 U\cdot v  =  \vec{e_1}  =   \bpm 1\\0\epm\]
with $U_{\theta}$ defined as in (\ref{decompose}).\\
Consider the affine deformation 
$\id: (X,\mu) \to (X,\mu_U) = (X,\mu)\circ U$ as in
Definition \ref{defdeform}. The
vector $v$ is mapped to $\vec{e_1}$.
Thus the same trajectories are now the horizontal ones.\\ 
Recall from \ref{affdeforms} that the flat structure $(X,\mu_U)$
is defined by the quadratic differential 
$e^{2\theta i}\cdot q$. Thus  $e^{2\theta i}\cdot q$
is Strebel. The ray is by (\ref{qminuszwei}) given as:
\begin{eqnarray*}
   \gamma_{\tilde{s}} = \gamma_{-e^{2\theta i}\cdot q}\; : \;
  [0,\infty) &\to& \tg \\ 
          t  &\mapsto & [\bpm 1&0\\0&K\epm \circ (X,\mu_{U_{\theta}}), \id]
                   \;\; =  \;\;[(X,\mu_{A_K}), \id]  
\end{eqnarray*}
with $K = e^t$ and 
$A_{K} = \bpm 1 & 0\\ 0& K \epm\cdot U_{\theta}$\\

\noindent
{\em The Strebel ray $\gamma_{\tilde{s}}$ is the image of the geodesic ray 
in $\HH$ from $i$ 
to the cusp $\tilde{s}$:}\\[1mm]
From 
Remark~\ref{ioeinszwei} (see also Figure \ref{gb})  one obtains that 
\[ \gamma_{\tilde{s}}(t) =   
          [(X,\mu_{A_K}), \id]
        =  \hat{\iota}(A_K) = \iota(-\overline{A_K^{-1}(i)}).\]
Furthermore, we have
\[-\overline{A_K^{-1}(i)})
        = -\overline{U_{\theta}^{-1}(K\cdot i)} = -U_{\theta}^{-1}(-Ki)
        = RU_{\theta}^{^-1}R^{-1}(Ki).\]
Thus the image of $\gamma_{\tilde{s}}$
is equal to the image of the the 
ray $RU_{\theta}^{-1}R^{-1}(Ki)$ ($K \in [1,\infty)$)
under $\iota$. But the latter one 
is the geodesic ray in $\HH$ 
from $i$ to $RU^{-1}R^{-1}(\infty) = -U^{-1}(\infty)$.\\
Observe finally that $-U^{-1}(\infty) = \tilde{s}$:\;
Since $U\cdot v = \vec{e_1}$ for the eigenvector $v$ of
$A$, one has for the fixed point $s$ of $A$ that $U(s) = \infty$. Hence,
one has for the fixed point $\tilde{s}$ of $\tilde{A} = RAR^{-1}$ that
$\tilde{s} = -s = -U^{-1}(\infty)$. Thus the Strebel ray defined
by $\gamma_{\tilde{s}}$ is the image of the 
geodesic ray from $i$ to $\tilde{s}$ in $\HH$ under $\iota$.\\

Finally we define 
$\iotaquer(\tilde{s}) 
  = (X_{\infty}(\tilde{s}),f_{\infty}(\tilde{s})) \in \tgquer$ 
to be the end point of the Strebel ray $\gamma_{\tilde{s}}$.
We then define the map $\iotaquer$ as follows.

\begin{definition}
$\iotaquer$\index{Teichm\"uller embedding!extension of} is
the extension of $\iota$  
defined by
\begin{eqnarray*}
\iotaquer: \quad \HH \cup \{\mbox{cusps of $\Gammaquerm_{\iota}$}\} 
              &\to& \tgquer,\\ 
            t &\mapsto& \left \{
               \begin{array}{l} 
                 \iota(t) , \mbox{ if } t \in \HH\\
                 \iotaquer(t) =  (X_{\infty}(\tilde{s}),f_{\infty}(\tilde{s}))
                 \mbox{ if } t = \tilde{s} \mbox{ is a cusp of  $\Gammaquerm$}
               \end{array} \right .
\end{eqnarray*}
\end{definition}

We consider $\HH \cup \{\mbox{cusps of $\Gammaquerm_{\iota}$}\}$
as topological space endowed with the horocycle topology
as in Example~\ref{horo}.

\begin{proposition} \label{iquer}
$\iotaquer$ is a continuous embedding. 
\end{proposition}

\begin{proof}

{\em $\iotaquer$ is continuous:}\\[2mm]
Let $s$ be a cusp of $\Gammaquerm_{\iota}$, i.e. $s$ is a fixed point
of some parabolic element $\tilde{A} \in \Gammaquerm_{\iota}$, and
$c:[0,\infty) \to \HH$ an arbitrary  path in $\HH$ 
 converging to $s$ in the horocycle
topology. \\

By Remark \ref{achtionstogether}, the action of $\tilde{A}$
on $\HH$ fits together with the action of $\rho(A) \in \Gamma_g$ on 
$\Tgq$. Both actions may be extended continously to 
$\HH_s = \HH \cup \{s\}$ (endowed with the horocycle topology)
and to $\tgquer$, respectively, and one obtains the following commutative 
diagram:
\[ \xymatrix{
     \HH_{s} = 
      \HH \cup \{s\} 
        \ar[rr]^{\iotaquer} \ar[d]_{p_{\tilde{A}}}  &&
      \tgquer\ar[d]^{p} \\
      \overline{\HH_s/<\tilde{A}>} \ar[rr]^{\iatilde}  && 
         \Mgq  
  }\]
Here the map $ \iatilde:\overline{\HH/<\tilde{A}>} \to \Mgq$ 
is the map induced by
$\iotaquer$ and $ \overline{\HH/<\tilde{A}>}$ \, is a disk with center 
$p_{\tilde{A}}(s)$.\\
Let $W$ be a neighbourhood of 
\[\bar{P}_{\infty} = \iatilde(\patilde(s)) = p(\iotaquer(s)).\]
For $i$ in an index set  $I$,
let 
$P_{\infty}^i$  be the preimages of 
$\bar{P}_{\infty}$ in $\tgquer$ under $p$. One of them is $\iotaquer(s)$, again by the 
commutativity of the
diagram.\\
Since $\{P_{\infty}^i|\; i \in I\}$ is discrete we may choose the neighbourhood $W$ 
in such a manner that 
its preimage under $p$  is of the form:
\[V = p^{-1}(W) = \bigcup_{i\in I} V_i \;\; \subseteq \, \tgquer\]
where the $V_i$ are the connected components of $V$ with 
$P_{\infty}^i \in V_i$ and $V_i$ is invariant under
the stabilizer of $P_{\infty}^i$ in the mapping class group $\Gamma_g$.\\
Furthermore, we may choose $W$ such that the preimage of $W$ 
under $\iatilde$ is a simply connected neighbourhood of $\patilde(s)$.
Then, again, the preimage 
\[U = \patilde^{-1}(\iatilde^{-1}(W))\]
is a neighbourhood of $s$ in the horocycle topology.\\
Thus an end piece of the path $c$ is completely contained in $U$, i.e. 
there is some $l \in \RR_{>0}$ such that $c([l,\infty))$
is contained in $U$.\\
Since the above diagram  is commutative and the $V_i$
are disjoint, the image of $U$ is one of the $V_i$.
This $V_i$ then contains $\iotaquer(c[l,\infty))$.
In addition, $V_i$ has to contain the end piece of the Strebel ray 
that leads to $s$ used to define $\iotaquer(s)$. Hence,
$V_i$ is the component that contains $\iotaquer(s)$.\\
Making $W$ arbitrarily small, the neighbourhood U of $s$ becomes
arbitrarily small. Thus $\iota\circ c$ converges to $\iotaquer(s)$.\\

\noindent{\em $\iotaquer$ is injective:}\\[2mm]
Suppose there are two cusps $s_1$ and $s_2$ with 
$P_{\infty} = \iotaquer(s_1) = \iotaquer(s_2)$. Thus we have two Strebel rays
defined by the negative  of the Strebel differentials $q_1 = e^{i\theta_1}\cdot q$ and 
$q_2 = e^{i\theta_2}\cdot q$
with initial point $P_{0} = \iota(i)$ and the same end point $P_{\infty}$ 
in $\tgquer$. Let $(X_{\infty},f_{\infty})$ and $(Y_{\infty},g_{\infty})$
be the two marked stable Riemann surfaces defined by the two Strebel rays, respectively.
Since they define  the same point in $\tgquer$ the following diagram
is commutative up to homotopy with some biholomorphic $h$:
\[ \xymatrix{
                            & X_{\infty}\\
     X_{\mbox{\fs ref}} \ar[ru]^{f_{\infty}} \ar[rd]^{g_{\infty}}&           \\
                            & Y_{\infty} \ar[uu]_{h} 
  }\]
The core curves of the cylinders relative to the flat structure 
on $X$ defined
by $q_1$ are mapped by $f_{\infty}$ to the singular points of 
$X_{\infty}$. Similarly
the core curves coming from $q_2$ are mapped to the singular points of $Y_{\infty}$.
Since the diagram is commutative up to isotopy, the two  
systems of core curves are homotopic. Thus the two Strebel rays 
are similar by definition,
using the terminology in \cite[Section 5]{M}. From Theorem~2 in 
\cite{M}
it follows that there is some constant $M < \infty$  such that for
two points $Q \neq R$ lying on the two Strebel rays 
which are equidistant from the 
initial point $P_0$, one has
$d(Q,R) \leq M$. But then, since $\iota$ is an isometric embedding,
$M$ would have to be an upper bound for the distance of equidistant points
on two different geodesic rays in $\HH$ starting from $i$. 
This cannot be true. 
\end{proof}

\subsubsection[Boundary of Teichm\"uller disks that
lead to Teichm\"uller curves]{Boundary of Teichm\"uller disks 
that lead to Teichm\"uller curves\\[2mm]}\label{igitter}\index{Teichm\"uller disk!boundary points of}\index{Teichm\"uller curve!boundary points of}

Let now $\iota: \HH \into \tg$  be a Teichm\"uller embedding such that
its image $\Delta_{\iota}$ projects to a Teich\-m\"uller curve 
$C$ in the moduli space $\mg$. 

\begin{proposition}\label{iquerfortc}
In this situation, the extended embedding from \ref{iotaquer}
\[\iotaquer:\; \HH \cup \{\mbox{\em cusps of }\Gammaquerm_{\iota}\} 
   \;\;\hookrightarrow\;\; \Deltaquer_{\iota} \; \subseteq \; \tgquer \]
is surjective onto the closure $\Deltaquer_{\iota}$
of $\Delta_{\iota}$ in $\tgquer$.
\end{proposition}

\begin{proof}
Recall from  Corollary~\ref{latticeproperty} that 
if $\iota$ leads to a Teichm\"uller curve then
the projective Veech group $\Gammaquer = \Gammaquer_{\iota}$ 
is a lattice in $\pslzwei(\RR)$, 
$\HH/\Gammaquerm$ is a complex algebraic curve and  
$\HH/\Gammaquerm \to C \subset\mg$ is
the normalization of $C$. 
Thus it extends to a surjective morphism 
\[\varphi: \overline{\HH/\Gammaquerm} 
\;\; \to \;\;  \Cquer \quad \subseteq \; \Mgq,\]
where $\overline{\HH/\Gammaquerm}$ and $\Cquer$ are  
the projective closure of $\HH/\Gammaquerm$ and the closure of 
$C$ in $\Mgq$, respectively.\\
Furthermore, the map $\HH \to \HH/\Gammaquerm$ extends continuously to 
a surjective map $p_{\Gammaquer}: \HH \cup \{\mbox{cusps of } \Gammaquerm\} 
\to \overline{\HH/\Gammaquerm}$, since $\Gammaquerm$ is a lattice
in $\pslzwei(\RR)$. 
Here we use  the horocycle topology on 
$\HH \cup \{\mbox{cusps of } \Gammaquerm\}$.\\
Thus one has the following commutative diagram of continuous maps:
\[ \xymatrix{
     \HH \cup \{\mbox{cusps of } \Gammaquerm\} 
        \ar[rr]^{\iotaquer} \ar[d]_{p_{\Gammaquer}}  &&
    \;\;\Deltaquer_{\iota} \;\;  \subseteq  \; \Tgq 
          \ar@<-3.5ex>[d]_{p|_{\,\Deltaquer_{\iota}}} 
          \ar@<3.5ex>[d]^{p} \\
     \overline{\HH/\Gammaquerm} \ar[rr]^{\varphi} && 
    \;\;\cquer \;\; \subseteq \; \Mgq  
  }\]

Let now $P_{\infty}$ be a point on the boundary of $\Delta_{\iota}$. Similarly
as in the proof of the continuity of $\iotaquer$ we may choose
a neighbourhood $W$ of $p(P_{\infty})$ in $\cquer$ 
such that all connected components
$V_i$ of the preimage $p^{-1}(W)$ contain only one preimage of 
$p(P_{\infty})$.
One of them, let's say $V_0$, contains  of course $P_{\infty}$ itself.\\
We choose an arbitrary path 
$c_{\iota}:[0,\infty) \to W\backslash\{p(P_{\infty})\} \, 
\subseteq C$ that converges to $p(P_{\infty})$.
Let $\hat{c}_{\iota}:[0,\infty) \,\to\, V_0$  be an arbitrary lift 
of $c_{\iota}$ via $p$ in $V_0$. Since we may choose $W$ arbitrarily small, 
$V_0$
may  become arbitrarily small and
$\hat{c}_{\iota}$ converges to $P_{\infty}$.\\
Now let 
$c:[0,\infty) \,\to\, \HH$ be the preimage of $\hat{c}_{\iota}$ 
under $\iota$, i.~e.
the path such that $\iota\circ c = \hat{c}_{\iota}$.
We project it by $p_{\Gammaquer}$ to $\overline{\HH/\Gammaquerm}$,
i.~e. we take the path  $p_{\Gammaquer}\circ c$.
Its image under $\varphi$ is  
$\varphi \circ p_{\Gammaquer}\circ c = p\circ \hat{c}_{\iota} =
c_{\iota}$ 
and converges to $p(P_{\infty})$ in $\cquer$.
Thus $p_{\Gammaquer}\circ c$ converges
in $\overline{\HH/\Gammaquerm}$, 
since $\varphi$
is an open map.\\ 
Since also  $p_{\Gammaquer}$ is open,
$c$ converges to some 
$t_{\infty} \in \HH \cup \{\mbox{cusps of } \Gammaquerm\}$. By continuity
of $\iotaquer$ one has
$\iotaquer(t_{\infty}) = P_{\infty}$. Thus $\iotaquer$ 
is surjective onto $\Deltaquer_{\iota}$.
\end{proof}

One obtains immediately the following conclusions.

\begin{corollary} \label{finalcor} If $\iota:\HH \into \tg$ 
leads to a Teichm\"uller curve $C$, then
\begin{enumerate}
\item[a)] the boundary points of the Teichm\"uller disk $\Delta_{\iota}$
are precisely the end points of the Strebel rays in $\Delta_{\iota}$
with initial point $\iota(i)$.
\item[b)] These boundary points correspond to the fixed points of parabolic elements 
in the projective Veech group\index{Veech group!parabolic elements of}.
\item[c)] Each boundary point of the Teichm\"uller curve $C$ is obtained
by contracting the core curves of the cylinders\index{cylinder!core curve of} in 
the direction of  $v$, where
$v$ is the eigenvector of a parabolic element in the Veech group. 
\end{enumerate}
\end{corollary}

This finishes the proof of Theorem \ref{thm}.


\section{Schottky spaces}
In this chapter we first recall the construction of Schottky coverings for
smooth and stable Riemann surfaces. We use them to define markings called
Schottky structures. In the smooth case they are classified by the well known Schottky space
$S_g$, a complex manifold of dimension $3g-3$ (if $g\ge 2$). In \cite{GH} it
was shown that also the Schottky structures on stable Riemann surfaces are
parameterized by a complex manifold $\Sgq$. Here we show how to obtain $\Sgq$
from Braungardt's extension $\Tgq$ of the Teichm\"uller space introduced in
Chapter~\ref{volker}. In the last section of this chapter we study the image of a Teichm\"uller disk in the Schottky space.

\subsection{Schottky coverings}
\label{s-struc}
We recall the basic definitions and properties of Schottky uniformization of
Riemann surfaces. We introduce the Schottky space $S_g$ and sketch, following
\cite{GH}, the construction of a universal family over it.
\begin{definition}
\label{s-group}
A group $\Gamma\subset\pslzwei(\CC)$ of M\"obius transformations on
$\PP^1(\CC)$ is called 
a {\it Schottky group}\index{Schottky group} if there are, for some $g\ge1$, disjoint closed simply connected domains
$D_1$, $D_1',\dots,D_g$, $D_g'$ bounded by Jordan curves $C_i = \partial D_i$,
$C_i' = \partial D_i'$, and generators $\gamma_1,\dots,\gamma_g$ of $\Gamma$
such that $\gamma_i(C_i)=C_i'$ and $\gamma_i(D_i) = \PP^1(\CC)-\bar D_i'$ for
$i=1,\dots,g$. The generators $\gamma_1,\dots,\gamma_g$ are called a {\it
  Schottky basis\index{Schottky basis}} of $\Gamma$.
\end{definition}
In Schottky's original paper \cite{S}, the $D_i$ in the definition were disks. With the same notation let
\[F=F(\Gamma)=\PP^1(\CC)-\cup_{i=1}^g(\bar D_i\cup\bar D_i')\ \ \mbox{and}\ \ 
\Omega=\Omega(\Gamma)=\cup_{\gamma\in\Gamma}\gamma(F).\]
It is well known, see e.\,g.\ \cite[X.H.]{Mask} that $\Gamma$ is a Kleinian
group, free of rank $g$ with free generators $\gamma_1,\dots,\gamma_g$, that
$\Omega$ is the region of discontinuity of $\Gamma$, and that
$X=\Omega/\Gamma$ is a closed Riemann surface of genus $g$. The quotient map
$\Omega\to X$ is called a {\it Schottky covering}\index{Schottky covering}.\\

An important fact is the following uniformization\index{Schottky uniformization} theorem:
\begin{proposition}
\label{s-unif}
Every compact Riemann surface $X$ of genus $g\ge1$ admits a Schottky covering
by a Schottky group of rank $g$.
\end{proposition}
\begin{proof}
The proof is based on the following construction that we shall extend to
stable Riemann surfaces in Section \ref{s-ext}: choose disjoint simple loops
$c_1,\dots,c_g$ on $X$ which are independent in homology, i.\,e.\
$F=X-\cup_{i=1}^gc_i$ is connected. Then $F$ is conformally equivalent to a
plane domain that is bounded by $2g$ closed Jordan curves. For $i=1,\dots,g$
denote by $C_i$ and $C_i'$ the two boundary components of $F$ that result from
cutting along $c_i$. Now let $\Phi_g$ be a free group on generators
$\varphi_1,\dots,\varphi_g$, and take a copy $F_w$ of $F$ for every element
$w\in\Phi_g$. The $F_w$ are glued according to the following rule: if $w$ and
$w'$ are reduced words in $\varphi_1,\dots,\varphi_g$ and if $w=w'\varphi_i$
then the boundary component $C_i$ on $F_{w'}$ is glued to $C_i'$ on $F_w$; if
$w$ ends with $\varphi_i^{-1}$ the roles of $C_i$ and $C_i'$ are
interchanged. By this construction we obtain a plane domain $\Omega$
together with a holomorphic action of $\Phi_g$ on it: an element
$\varphi\in\Phi_g$ maps the copy $F_w$ to $F_{w\varphi}$. The crucial step in
the proof now is to show that this action extends to all of $\PP^1(\CC)$,
i.\,e.\ $\Phi_g$ acts by M\"obius transformations. For this we refer to \cite[Ch.~IV, Thm.~19\,F]{AS}.
\end{proof}
\begin{definition}
\label{sg}
Let $\Sgs$ be the set of all $(\gamma_1,\dots,\gamma_g)\in\pslzwei(\CC)^g$
that generate a Schottky group $\Gamma$ and form a Schottky basis for
$\Gamma$. The set $S_g$ of equivalence classes of $g$-tuples
$(\gamma_1,\dots,\gamma_g)\in\Sgs$ under simultaneous conjugation is called
the {\it Schottky space}\index{Schottky space} of genus $g$.
\end{definition}
For a point $s=(\gamma_1,\dots,\gamma_g)\in\Sgs$ let $\Gamma(s)$ be the
Schottky group generated by $\gamma_1,\dots,\gamma_g$, $\Omega(s)$ the region
of discontinuity of $\Gamma(s)$, and $X(s)=\Omega(s)/\Gamma(s)$ the associated
Riemann surface. This leads to an alternative description of the Schottky space:
\begin{remark}
\label{sg-alt}
$S_g$ is the set of equivalence classes of pairs $(X,\sigma)$, where $X$ is a Riemann surface of genus $g$ and $\sigma:\Phi_g\to\empslzwei(\CC)$ is an injective homomorphism such that $\Gamma:=\sigma(\Phi_g)$ is a Schottky group and $\Omega(\Gamma)/\Gamma\cong X$.\\
$(X,\sigma)$ and $(X',\sigma')$ are equivalent if there is some $A\in\empslzwei(\CC)$ such that $\sigma'(\gamma)=A\sigma(\gamma)A^{-1}$ for all $\gamma\in\Phi_g$. Note that then $X'$ is isomorphic to $X$.
\end{remark}
To endow $S_g$ with a complex structure we proceed as follows: Taking the fixed points and the multipliers of the $\gamma_i$ we obtain an embedding of $\Sgs$ as an open subdomain of $\PP^1(\CC)^{3g}$. For $g=1$ each equivalence class contains a unique M\"obius transformation of the form $z\mapsto \lambda z$ for some $\lambda\in\CC$, $0<|\lambda|<1$.
If $g\ge 2$ we find in each equivalence class in $\Sgs$ a unique representative
$(\gamma_1,\dots,\gamma_g)$ such that $\gamma_1$ and $\gamma_2$ have
attracting fixed points 0 and 1, respectively, and $\gamma_1$ has repelling
fixed point $\infty$. This defines a section to the projection $\Sgs\to S_g$ 
and embeds $S_g$ as a closed subspace of $\Sgs$ which, moreover, 
lies in $\{0\}\times\{\infty\}\times\{1\}\times \CC^{3g-3} \subseteq 
\PP^1(\CC)^{3g}$. Thus we have shown.
\begin{proposition}
\label{s-open}
{\bf a)} $S_1$ is a punctured disk.\\[1mm]
{\bf b)} For $g\ge2$, $S_g$ carries a complex structure as an open
subdomain of $\CC^{3g-3}$.
\end{proposition}
Our next goal is to show that this complex structure on $S_g$ is natural. The main step in this direction is
\begin{proposition}
\label{mu}
The forgetful map $\mu:S_g\to M_g$, that sends $s=(X,\sigma)$ to the isomorphism class of $X$, is analytic and surjective.
\end{proposition}
\begin{proof} The surjectivity of $\mu$ follows from Prop.~\ref{s-unif}. 
To show that $\mu$ is analytic we use the fact that $M_g$ is a coarse moduli 
space for Riemann surfaces. Therefore it suffices  to find a holomorphic
family $\pi:\Cg\to S_g$ of Riemann 
surfaces\index{family of Riemann surfaces} over $S_g$ which induces $\mu$ in the sense that for $s\in S_g$,
$\mu(s)$ is the isomorphism class of the fibre $C_s=\pi^{-1}(s)\subset{\cal
  C}_g$. \\[1mm]
The family $\Cg$ is obtained as in Section
\ref{v-structure}: Let 
\[\Omega_g=\{(s,z)\in S_g\times\PP^1(\CC):z\in\Omega(s)\}.\]
$\Omega_g$ is a complex manifold on which the free group $\Phi_g$ acts
holomorphically by $\varphi(s,z) = (s,\sigma(\varphi)(z))$ for $s=(X,\sigma)\in S_g$, $\varphi\in\Phi_g$ and $z\in\Omega(s)$.\\[1mm]
The projection pr$_1:\Omega_g\to S_g$ onto the first component factors through
the orbit space $\Cg = \Omega_g/\Phi_g$, and the induced map $\pi:\Cg\to S_g$
is the family of Riemann surfaces we were looking for.
\end{proof}
The family $\Cg$ is in fact universal for Riemann surfaces with {\it Schottky
  structure}\index{Schottky structure}, a kind of marking that we now recall from \cite[Section~1.3]{GH}:
\begin{definition}
\label{s-structure}
{\bf a)} Let ${\cal U}\to S$ be an analytic map of complex manifolds and $\Gamma\subset\mbox{Aut}({\cal U}/S)$ a properly discontinuous subgroup. Then the analytic quotient
map ${\cal U}\to{\cal U}/\Gamma={\cal C}$ is called a {\it Schottky covering}
if the induced map ${\cal C}\to S$ is a family of Riemann surfaces and
if for every $x\in S$ the restriction $U_x\to C_x$ of the quotient map to the fibres is a Schottky
covering.\\[1mm]
{\bf b)} A {\it Schottky structure} is a Schottky covering \,${\cal U}\to{\cal
  U}/\Gamma={\cal C}$ together with an equivalence class of isomorphisms
$\sigma:\Phi_g\to\Gamma$, where $\sigma$ and $\sigma'$ are considered
equivalent if they differ only by an inner automorphism of $\Phi_g$.
\end{definition}
Note that the construction in the proof of Proposition \ref{mu} endows the
family $\Cg/S_g$ with a Schottky structure. \\

A Schottky structure on a single Riemann surface $X$ is given by a Schott\-ky covering $\Omega\to\Omega/\Gamma=X$ and an isomorphism $\sigma:\Phi_g\to\Gamma$. Comparing the respective equivalence relations we find that the points $(X,\sigma)$ in $S_g$ correspond bijectively to the isomorphism classes of Riemann surfaces with Schottky structure. In fact a much stronger result holds:
\begin{theorem}
\label{s-fine}
$S_g$ is a fine moduli space\index{fine moduli space!for Riemann surfaces with Schottky structure} for Riemann surfaces with Schottky structure.
\end{theorem}
\begin{proof}
Let ${\cal C}/S$ be a family of Riemann surfaces and $({\cal U}\to{\cal
  U}/\Gamma={\cal C},\sigma:\Phi_g\stackrel{\sim}{\longrightarrow}\Gamma)$ a
  Schottky structure on ${\cal C}$. Then we have a map $f:S\to S_g$ which maps
  a point $x$ to the isomorphism class of the Schottky covering $U_x\to C_x$. We
  have to show that $f$ is analytic. Then the other properties of a fine
  moduli space follow easily from the definitions, namely that ${\cal C}$ is
  the fibre product $\Cg\times_{S_g}S$ and that ${\cal U}$ is isomorphic to
  $\Omega_g\times_{\Cg}{\cal C} = \Omega_g\times_{S_g}S$ such that the
  projection ${\cal U}\to\Omega_g$ onto the first factor is equivariant for
  the actions of $\Gamma$ and $\Phi_g$ via the isomorphism $\sigma$.\\
The universal property of $M_g$ as a coarse moduli space gives us, as above for
  $\mu$, that the composition $\mu\circ f$ is analytic. Since $\mu$ has
  discrete fibres, it therefore suffices to show that $f$ is continuous. This
  is quite subtle, see \cite[\S\,3]{GH}.
\end{proof}

\subsection{Relation to Teichm\"uller space}
\label{s-teich}
In this section we explain that Schottky space\index{Schottky space} can 
be obtained as a quotient
space of the Teichm\"uller space which was introduced in
Section \ref{v-structure}. For this purpose we first endow the universal 
family
$\Cgnull$ over the Teichm\"uller space $T_g=T_{g,0}$ with a 
Schottky structure as follows: \\

Let $a_1,b_1,\dots,a_g,b_g$ be a set of standard generators\index{standard generators} of $\pi_g$, the
fundamental group of the reference surface $\Cref$; this means that they
satisfy the relation $\Pi_{i=1}^g a_ib_ia_i^{-1}b_i^{-1}=1$. Then
$b_1,\dots,b_g$ are homologically independent, hence the construction in the proof of Prop.~\ref{s-unif} provides us with a corresponding
Schottky covering $\Oref\to\Cref$. The group Aut$(\Oref/\Cref)$ of deck
transformations is isomorphic to the free group on $b_1,\dots,b_g$. Denoting
$\Uref\to\Cref$ the universal covering, there is a covering
map $\Uref\to\Oref$ over $\Cref$. The group Aut$(\Uref/\Oref)$ is the kernel
$N_\alpha$ of the homomorphism $\alpha:\pi_g\to\Phi_g$ which maps $b_i$ to
$\varphi_i$ and $a_i$ to 1; in other words, $N_\alpha$ is the normal closure in
$\pi_g$ of the subgroup generated by $a_1,\dots,a_g$.\\

In Section \ref{v-structure} we described the family
$\Omega^+_{g,0}\to{\cal C}_{g,0}$ of universal coverings 
of the surfaces in the family ${\cal C}_{g,0}$; the fundamental group $\pi_g$ and hence also $N_\alpha$ acts on the fibres of this covering, and we obtain:
\begin{remark}
\label{ssonteich}
The induced map $\Omega^+_{g,0}/N_\alpha\to{\cal C}_{g,0}$ is a Schottky covering, and the universal Teichm\"uller structure\index{Teichm\"uller structure!universal}
$\tau:\pi_g\stackrel{\sim}{\longrightarrow}\mbox{\em Aut}(\Omega^+_{g,0}/{\cal
  C}_{g,0})$ (cf.~Theorem~\ref{fein}) descends via $\alpha$ to a Schottky structure
$\sigma_\alpha:\Phi_g = \pi_g/N_\alpha
\stackrel{\sim}{\longrightarrow}\mbox{\em Aut}((\Omega^+_{g,0}/N_\alpha)/{\cal C}_{g,0})$
on ${\cal C}_{g,0}$. 
\end{remark}
By Theorem \ref{s-fine} this Schottky structure induces an analytic map
$s_\alpha:T_g\to S_g$. To describe $s_\alpha$ as the quotient map for 
a subgroup of the mapping class group\index{mapping class group} $\Gamma_g$, 
we first identify $\Gamma_g$ with the group
$\mbox{Out}^+(\pi_g)$ of orientation preserving outer automorphisms of
$\pi_g$; then, to a diffeomorphism $f:\Cref\to\Cref$, we associate the induced automorphism $\varphi=f_*:\pi_g\to\pi_g$. It follows from the Dehn-Nielsen theorem that this gives an isomorphism $\Gamma_g\stackrel{\sim}{\longrightarrow}\mbox{Out}^+(\pi_g)$. In this chapter, by $\varphi\in\Gamma_g$ we always mean an element of $\mbox{Out}^+(\pi_g)$.
\begin{proposition}
\label{teichtos}
{\bf a)} $s_\alpha$ is the quotient map for the subgroup
\[\Gamma_g(\alpha)=\{\varphi\in\Gamma_g:\alpha\circ\varphi\equiv\alpha\ \mbox{\em
  mod Inn}(\pi_g)\}\]
of the mapping class group $\Gamma_g$ (where $\mbox{\em Inn}(\pi_g)$ denotes the group of inner automorphisms).\\[1mm] 
{\bf b)} $s_\alpha:T_g\to S_g$ is the universal covering of the Schottky
space.\\[1mm] 
{\bf c)} $s_\alpha$ lifts to maps $\tilde s_\alpha$ and $\omega_\alpha$ that
make the following diagram commutative:
\begin{center}
$\xymatrix@=6ex{\ar[dr]^{/N_\alpha}\ar[dd]_{/\pi_g}\Omega^+_{g,0}\\
&\ar[r]^{\omega_\alpha}\ar[dl]\Omega^+_{g,0}/N_\alpha&\ar[d]^{/\Phi_g}\Omega_g\\
\ar[d]\ar[rr]^{\tilde s_\alpha}{\cal C}_{g,0}&&\ar[d]{\cal C}_g\\
\ar[rr]^{s_\alpha}\ar[dr]_{/\Gamma_g}T_g&&\ar[dl]^\mu S_g\\
&M_g}$
\end{center}
\end{proposition}
\begin{proof}
{\bf a)} Let $x=(X,f)\in T_g$. Recall from Section \ref{v-structure} that the fibre over $x$ in $\Omega^+_{g,0}$ is the component $\Oplusx$ of the region of discontinuity of the quasifuchsian group $G_x$ associated to $x$. The universal Teichm\"uller structure on ${\cal C}_{g,0}$ induces an isomorphism $\tau_x:\pi_g\to G_x=\mbox{Aut}(\Oplusx/X)$. From Remark \ref{ssonteich} we see that the point $s_\alpha(x)=(X,\sigma)\in S_g$ is given by the restriction $\sigma_{\alpha,x}$ of $\sigma_\alpha$ to the fibre over $x$; explicitly,
\[\sigma=\sigma_{\alpha,x}:\Phi_g=\pi_g/N_\alpha\stackrel{\sim}{\longrightarrow}\mbox{Aut}((\Oplusx/\tau_x(N_\alpha))/X)=G_x/\tau_x(N_\alpha).\]
For $\varphi\in\Gamma_g$ we have $s_\alpha(x)=s_\alpha(\varphi(x))$ if and only if $\sigma_{\alpha,x}=\sigma_{\alpha,\varphi(x)}$ up to an inner automorphism.
Since $\tau_{\varphi(x)}=\tau_x\circ\varphi^{-1}$ this happens if and only if $\varphi$ induces an inner automorphism on $\pi_g/N_\alpha$, i.\,e.\ if and only if $\varphi\in\Gamma_g(\alpha)$.\\[1mm]
{\bf b)} This is clear from the fact that $T_g$ is simply connected and
$\Gamma_g(\alpha)$ is torsion free, hence $s_\alpha$ is unramified. Using the
construction in a) one can give a direct proof which in turn provides an
independent proof that $T_g$ is simply connected, see
\cite[Prop.~6]{GH}.\\[1mm] 
{\bf c)} It follows from Remark \ref{ssonteich} that $\Omega^+_{g,0}/N_\alpha
\to {\cal C}_{g,0}$ is a Schottky covering. Therefore, by the universal
property of $S_g$ (Theorem \ref{s-fine}), ${\cal C}_{g,0}$ is the fibre
product $T_g\times_{S_g}{\cal C}_g$, and $\tilde s_\alpha$ is the projection
to ${\cal C}_g$. Moreover the Schottky covering $\Omega^+_{g,0}/N_\alpha
\to {\cal C}_{g,0}$ is a pullback of the universal Schottky covering $\Omega_g
\to{\cal C}_g$, i.\,e.\ $\Omega^+_{g,0}/N_\alpha = {\cal C}_{g,0}\times_{{\cal
    C}_g}\Omega_g$, and again $\omega_\alpha$  is the projection to the second
factor. 
\end{proof}
In fact, the action of $\Gga$ on $T_g$ extends to $\Omega^+_{g,0}$, $\Omega^+_{g,0}/N_\alpha$ and
$\Cgnull$; then $\tilde s_\alpha$ and $\omega_\alpha$ are the quotient maps
for these actions.

\subsection{Schottky coverings of stable Riemann surfaces}
\label{s-stab}
In this and the following section we introduce a partial compactification $\Sgq$ of\index{Schottky space!partial compactification of}
$S_g$ that fits in between $\Tgq$ and $\Mgq$. We have presented two different
ways to define $S_g$, and we shall see that both are suited for extension to
stable Riemann surfaces: The first way is to construct Schottky coverings for
surfaces with nodes, define Schottky structures and find parameters for
them. This approach was pursued in \cite{GH} and will be sketched in this
section. The other possibility is to extend the action of $\Gga$ to (part of)
the boundary of $\Tgq$ and show that the quotient exists and has the desired
properties; this will be done in Section \ref{s-ext}.
\begin{definition}
\label{cut}
Let $X$ be a stable Riemann surface\index{stable Riemann surface} of genus $g$. A {\it cut system}\index{cut system} on $X$ is
a collection of disjoint simple loops $c_1,\dots,c_g$ on $X$, not
passing through any of the nodes, such that $X-\cup_{i=1}^gc_i$ is connected.
\end{definition}
\begin{proposition}
\label{cut-exist}
On any stable Riemann surface there exist cut systems.
\end{proposition}
\begin{proof}
Let $f:\Cref\to X$ be a deformation; we must find
disjoint and homologically independent loops $\tilde c_1,\dots,\tilde c_g$ on
$\Cref$ that are disjoint from the loops $a_1,\dots,a_k$ that are contracted
by $f$. For this we complete $a_1,\dots,a_k$ to a maximal system
$a_1,\dots,a_{3g-3}$ of homotopically independent loops (such a system
decomposes $\Cref$ into pairs of pants). Among the $a_i$ we find
$a_{i_1},\dots,a_{i_g}$ that are homologically independent. If $i_\nu>k$ we
take $\tilde c_\nu = a_{i_\nu}$, and for $i_\nu\le k$ we replace $a_{i_\nu}$
by a loop $\tilde c_\nu$ that is homotopic to $a_{i_\nu}$ and disjoint from
it.
\end{proof}
Once we have found $c_1,\dots,c_g$ as above, we proceed as in the proof of
Proposition \ref{s-unif} to construct a Schottky covering of $X$: Let
$F=X-\cup_{i=1}^gc_i$, take a copy $F_w$ of $F$ for each $w\in\Phi_g$, and
glue these copies exactly as before to obtain a space $\Omega$. Of course,
neither $F$ nor $\Omega$ is planar whenever $X$ has nodes. 
In all cases, the complex structure on $X$ lifts to a structure of
a one-dimensional complex space on $F$. The group
$\Phi_g$ acts on this space by holomorphic automorphisms.
Precisely, there is an isomorphism $\Phi_g\to\Gamma=\mbox{Aut}(\Omega/X)$, and $X$ is isomorphic to $\Omega/\Gamma$ as complex space.
 \begin{definition}
\label{s-covq}
The covering $\Omega\to X$ constructed above for a cut system
$c=(c_1,\dots,c_g)$ on a stable Riemann surface $X$ is called the {\it
Schottky covering}\index{stable Riemann surface!Schottky covering of} of $X$ relative to $c$. 
A covering of $X$ is called a
{\it Schottky covering} if it is the Schottky covering relative to some cut system.
\end{definition}
The next goal is to define a space $\Sgq$ that classifies Schottky coverings in a way
analogous to Definition \ref{sg}. Since the covering space $\Omega$ is in
general not a subspace of $\PP^1(\CC)$ and thus the group of deck
transformations not a subgroup of $\pslzwei(\CC)$, we cannot directly extend
\ref{sg}.\\

A closer look at the construction of a Schottky covering
$\Omega\to\Omega/\Gamma=X$ of a stable Riemann surface $X$ shows the following:\\
Each irreducible component $L$ of $\Omega$ is an open dense subset of a projective
line; more precisely, the stabilizer of $L$ in $\Gamma$ is a Schottky group as in
Definition \ref{s-group}, and $L$ is its region of discontinuity. Moreover the
intersection graph of the irreducible components of $\Omega$ is a
tree (hence $\Omega$ is called a {\it tree of projective lines})\index{tree of projective lines}.\\
Therefore, for each irreducible component $L$, there is a well defined
projection $\pi_L:\Omega\to L$ which is the identity on $L$: For an arbitrary point $x\in\Omega$ there is a unique chain $L_0,L_1,\dots,L_n=L$ of mutually distinct components such that $x\in L_0$ and $L_i$ intersects $L_{i+1}$ for $i=0,\dots,n-1$; then define $\pi_L(x)$ to be the intersection point of $L_{n-1}$ and $L$.\\
An {\it end} of $\Omega$ is an equivalence class of infinite chains $L_0,L_1,L_2,\dots$ of irreducible components as above (i.\,e.\ $L_i\not=L_j$ for $i\not=j$ and $L_i\cap L_{i+1}\not=\emptyset$), where two chains are equivalent if they differ only by finitely many components. Let $\Omega^*=\Omega\cup\{\mbox{ends of}\ \Omega\}$. Clearly the projection $\pi_L$ to a component $L$ can be extended to $\Omega^*$.\\
For any three different points or ends
$y_1, y_2, y_3$ in $\Omega^*$ there is a unique component $L=L(y_1,y_2,y_3)$
(called the {\it median} of the three points) such that the points
$\pi_L(y_1),\pi_L(y_2),\pi_L(y_3)$ are distinct. Now given any four distinct
points or ends $y_1,\dots,y_4$ in $\Omega^*$ we can define a {\it cross ratio}\index{cross ratio}
$\lambda(y_1,\dots,y_4)$ by taking the usual cross ratio of
$\pi_L(y_1),\dots,\pi_L(y_4)$ on the median component $L=L(y_1,y_2,y_3)$ of
the first three of them; note that $\lambda(y_1,\dots,y_4)$ will be 0, 1 or
$\infty$ if $\pi_L(y_4)$ coincides with $\pi_L(y_1)$, $\pi_L(y_2)$ or
$\pi_L(y_3)$. \\

To obtain parameters for the group $\Gamma$ observe that any
$\gamma\in\Gamma$, $\gamma\not=1$, has exactly two fixed points on the
boundary of $\Omega$, where boundary points of $\Omega$ are either points in
the closure of a component, or ends of $\Omega$; one of the fixed points is
attracting, the other repelling. For any four different (primitive) elements
$\gamma_1,\dots,\gamma_4$ in $\Gamma$ we define
$\lambda(\gamma_1,\dots,\gamma_4)$ to be the cross ratio of their attracting
fixed points. It is a remarkable fact that from these cross ratios both the
space $\Omega$ and the group $\Gamma\subset\mbox{Aut}(\Omega)$ can be recovered.
For any parti\-cular Schottky covering finitely
many of them suffice, but for different Schottky coverings we must take, in
general, 
the cross ratios of different elements of $\Phi_g$. To parameterize all
Schottky coverings we therefore have to use infinitely many of these cross
ratios. We consider them as (projective) coordinates on an infinite product
of projective lines $\PP^1(\CC)$. The cross ratios satisfy a lot of
algebraic relations, which define a closed subset $B$ of this huge space.
Every point of $B$ represents a tree of projective lines $\Omega$ as above
together with an action of $\Gamma$ on it. $\Sgq$ is the open subset
of $B$, where this action defines a Schottky covering. For details
and in particular the technical complication caused by the presence 
of infinitely many variables and equations, see 
\cite[\S2]{GH} and \cite{He}. In principle, one can proceed as in
Section~\ref{s-struc} to construct a family of stable Riemann surfaces over
$\Sgq$.\\

Given a family ${\cal C}/S$ of stable Riemann surfaces over a complex manifold
$S$, we can define the notion of a {\it Schottky covering} ${\cal U}/S\to {\cal
C}/S$ and of a {\it Schottky structure}\index{Schottky structure} on ${\cal U}$ exactly as in
Definition \ref{s-structure}, except that now ${\cal U}$ is not assumed to be
a manifold, but only a complex space. It is shown in \cite[\S3]{GH} that
the family over $\Sgq$ carries a universal Schottky structure:
\begin{theorem}
\label{sgq}
$\Sgq$ is a fine moduli space for stable Riemann surfaces with Schottky
structure. 
\end{theorem}

\subsection{$\Sgq$ as quotient of $\Tgq$}
\label{s-ext}
It is not possible to extend the quotient map $s_\alpha:T_g\to S_g$ constructed
in Section~\ref{s-teich} to the whole boundary of $T_g$ in $\Tgq$. Instead we shall, for
each $\alpha$, identify a part $\Tgqa$ of $\Tgq$ to which the action of $\Gga$
and hence the morphism $s_\alpha$ can be extended. It will turn out that the
quotient space is the extended Schottky space\index{Schottky space!extended} $\Sgq$ described in the previous
section.\\

We begin with the definition of the admissible group homomorphisms $\alpha$
and the associated parts $\Tgqa$ of $\Tgq$:
\begin{definition}
\label{sympl}
{\bf a)} A surjective homomorphism $\alpha:\pi_g\to\Phi_g$ is called
{\it symplectic}\index{symplectic homomorphism} if there are standard generators $a_1,b_1,\dots,a_g,b_g$ of
$\pi_g$ (in the sense of Section \ref{s-teich}) such that $\alpha(a_i)=1$ for
$i=1,\dots,g$.\\[1mm]
{\bf b)} Recall from Chapter \ref{volker} that a point in $\Tgq$ can be described as an equivalence class of
pairs $(X,f)$, where $X$ is a stable Riemann surface and $f:\Cref\to X$ is a
deformation (see Corollary \ref{homeo}).\\
For a symplectic homomorphism $\alpha:\pi_g\to\Phi_g$ let
\[\Tgqa =
\{(X,f)\in\Tgq:\mbox{ker}\,(f_*)\subseteq\mbox{ker}\,(\alpha)\}.\] 
\end{definition}
\begin{proposition}
\label{Tgqa}
{\bf a)} $\Tgqa$ is an open subset of $\Tgq$; it contains $T_g$ and is
invariant under the group $\Gga$ introduced in Prop.~\ref{teichtos}.\\[1mm]
{\bf b)} $\Tgq$ is the union of the $\Tgqa$, where $\alpha$ runs through the symplectic homomorphisms.\\[1mm]
{\bf c)} The restriction to $\Tgqa$ of the universal covering $p:\Tgq\to\Mgq$
is surjective for every symplectic $\alpha$.
\end{proposition}
\begin{proof}
{\bf a)} Let $(X,f)$ be a point in $\Tgq$ and $c_1,\dots,c_k$ the loops on
$\Cref$ that are contracted under $f$. Then the kernel of $\pi_1(f):\pi_g\to\pi_1(X)$ is the
normal subgroup generated by $c_1,\dots,c_k$. The local description of $\Tgq$
in Corollary \ref{local} shows that there is a neighbourhood $U$ of $(X,f)$ in
$\Tgq$ such that for every $(X',f')\in U$ the map $f':\Cref\to X'$ contracts a
subset of $\{c_1,\dots,c_k\}$. Hence the kernel of $\pi_1(f')$ is contained in
ker$\,(\pi_1(f))$. Thus if $(X,f)\in\Tgqa$, also $U\subseteq\Tgqa$. The
remaining assertions are clear.\\[1mm]
{\bf b)} Again let $(X,f)$ be a point in $\Tgq$ and $c_1,\dots,c_k$ the
loops on $\Cref$ contracted by $f$. By Proposition~\ref{cut-exist} we can find a cut system
$a_1,\dots,a_g$ on $X$ and a corresponding Schottky covering.
 This covering induces a surjective homomorphism
$\pi_1(X)\to\Phi_g$. Composing this homomorphism with $\pi_1(f)$ yields a
homomorphism $\alpha:\pi_g\to\Phi_g$ which corresponds to a Schottky covering of
$\Cref$ (relative to the cut system $f^{-1}(a_1),\dots,f^{-1}(a_g)$) and hence
is symplectic. By construction, $c_1,\dots,c_k$ are in the kernel of
$\alpha$.\\[1mm] 
{\bf c)} Let $\alpha:\pi_g\to\Phi_g$ be symplectic and $a_1,b_1,\dots,a_g,b_g$
standard generators of $\pi_g$ such that $\alpha(a_i) = 1$ for all $i$. For an
arbitrary stable Riemann surface $X$ choose a deformation $f:\Cref\to X$ and
let $c_1,\dots,c_k$ be the loops that are contracted by $f$. As
in the proof of b) we find standard generators $a'_1,b'_1,\dots,a'_g,b'_g$
such that the $c_j$ are contained in the normal subgroup generated by
the $a'_i$. The map $a_i\mapsto a'_i, b_i\mapsto b'_i$ defines an
automorphism $\varphi$ of $\pi_g$ and thus an element of $\Gamma_g$. Then by construction
$(X,f\circ\varphi)$ lies in $\Tgqa$ and $p(X,f\circ\varphi) = X$.
\end{proof}
As a side remark we note that
$\Tgqa$ is not only invariant under $\Gga$, but also under the larger ``handlebody'' group\index{handlebody group}
\[\Hga = \{\varphi\in\Gamma_g:\varphi(N_\alpha) = N_\alpha\}\]
(where $N_\alpha = \mbox{ker}\,(\alpha)$ as in Section \ref{s-teich}). Note
that $\Hga$ is the normalizer of $\Gga$ in $\Gamma_g$, and that we have an
exact sequence 
\[1\to\Gga\to\Hga\to\mbox{Out}\,(\Phi_g)\to 1.\]
The quotient space $\Sgd=T_g/\Hga=S_g/\mbox{Out}\,(\Phi_g)$ is a parameter
space for Schottky groups of rank $g$ (without any marking).
\begin{proposition}
\label{Sgqa}
For any symplectic homomorphism $\alpha:\pi_g\to\Phi_g$, the quotient space
$\Tgqa/\Gga$ is a complex manifold $\Sgqa$.
\end{proposition}
\begin{proof}
This is a local statement which is clear for points $(X,f)\in T_g$ since
$\Gga$ is torsion free. For an
arbitrary $x=(X,f)\in\Tgq$ we saw in Section \ref{v-tgnq} that the Dehn twists
$\tau_1,\dots,\tau_k$ around the loops $c_1,\dots,c_k$ that are contracted by
$f$ generate a finite index subgroup $\Gamma_x^0$ of the stabilizer $\Gamma_x$
of $x$ in $\Gamma_g$ (the quotient being the finite group Aut$\,(X)$). Let $\alpha$ be a symplectic
homomorphism with respect to standard generators $a_1,b_1,\dots,a_g,b_g$, and
assume $(X,f)\in\Tgqa$. Since the $c_i$ are in the normal subgroup generated
by $a_1,\dots,a_g$, they do not intersect any of the $a_j$ and thus
$\tau_i(a_j) = a_j$ for all $i$ and $j$. This shows $\Gamma_x\subseteq\Gga$.\\
Now choose a neighbourhood $U$ of $x=(X,f)$ in $\Tgqa$ which is precisely
invariant under $\Gamma_x$. Then it follows, from Proposition \ref{tgnq} 
(and
Definition \ref{covmgnq}), that $U/\Gamma_x$ is a complex manifold.
\end{proof}
For any two sets $a_1,b_1,\dots,a_g,b_g$ and $a'_1,b'_1,\dots,a'_g,b'_g$ of
standard generators, $a_i\mapsto a'_i$, $b_i\mapsto b'_i$ defines an
automorphism of $\pi_g$. Therefore for any two symplectic homomorphisms
$\alpha$ and $\alpha'$ there is an automorphism $\psi\in\Gamma_g$ such that
$\alpha = \alpha'\circ\psi$. Then clearly $N_\alpha = \psi(N_{\alpha'})$
and $\Ggas = \psi\Gga\psi^{-1}$. This shows that, as an automorphism of
$\Tgq$, $\psi$ maps $\Tgqa$ to $\Tgqas$ and descends to an isomorphism
$\bar\psi:\Sgqa\to\Sgqas$. We have shown:
\begin{remark}
\label{nureinSgqa}
The complex manifolds $\Sgqa$ are isomorphic for all symplectic homomorphisms
$\alpha$.
\end{remark}
It remains to show that the $\Sgqa$ coincide with the fine moduli space $\Sgq$
of Section~\ref{s-stab}. This is achieved by showing that $\Sgqa$ satisfies
the same universal property as $\Sgq$:
\begin{proposition}
\label{sgqa=sgq}
For any symplectic $\alpha$, $\Sgqa$ is a fine moduli space for stable Riemann
surfaces with Schottky structure and hence isomorphic to $\Sgq$.
\end{proposition}
\begin{proof}
The idea of the proof is to endow the universal family over $\Tgqa$ with
a Schottky structure and to transfer this to a Schottky structure on the image
family over $\Sgqa$.\\[1mm] 
Before explaining this for the whole family we consider
a single stable Riemann surface $X$. Let $d_1,\dots,d_k$ be 
the nodes on $X$, $f:\Cref\to X$ a
deformation and 
$\alpha:\pi_g\to\Phi_g$ a symplectic homomorphism such that
$x=(X,f)\in\Tgqa$. In Section \ref{v-structure} we described the universal
covering $\Oplusxd\to\Oplusxd/G_x=X$ of $X$ with cusps over the nodes. Recall
that  $\Oplusxd$ is the union of the plane region $\Oplusx$ with the common
boundary points
of the doubly cusped regions lying over the nodes $d_i$, and that $G_x$ is
isomorphic to $\pi_g$. 
\begin{remark}
\label{OplusxdmodNa}
Using the above notation, let $\rho:\pi_g\to G_x$ be an isomorphism and
$N_\alpha^{G_x}=\mbox{\em ker}\,(\alpha\circ\rho^{-1})\subseteq G_x$. Then
$\Omega=\Oplusxd/N_\alpha^{G_x}$ is a complex space, $G_x/N_\alpha^{G_x}\cong\Phi_g$ acts
holomorphically on $\Omega$, and $\Omega\to\Omega/\Phi_g=X$ is a Schottky
covering. 
\end{remark}
\begin{proof}
The key observation is that the stabilizer in $G_x$ of a point $\tilde d_i\in\Oplusxd$
lying over $d_i$ is
generated by an element $\gamma_i$ corresponding  under $\rho$ to a conjugate
of the loop
$f^{-1}(d_i)$. Since we assumed $(X,f)\in\Tgqa$, we have $\gamma_i\in
N_\alpha^{G_x}$. This shows that $\Omega$ is a complex space, more precisely: a
Riemann surface with nodes. The other assertions then follow directly from the definitions.
\end{proof}
The above construction can be carried over to families in the following way:
First consider the universal family $\Cgq$ over $\Tgq$ and the universal
Teichm\"uller structure $\Oplusgd\to\Cgq$ on it. Denote by $\Cgqa$ resp.\ 
$\Oplusgda$ the restriction to $\Tgqa$. Then the quotient space
$\Oplusgda/N_\alpha$ is a complex space on which $\Phi_g=\pi_g/N_\alpha$ acts.
The quotient map $\Oplusgda/N_\alpha\to\Cgqa$ is a Schottky covering and the
identification of $\Phi_g$ with the group of deck transformations defines a
Schottky structure. \\[1mm]
The group $\Gga$ acts not only on $\Tgqa$, but also on $\Oplusgda$ 
as follows: for
$\varphi\in\Gga$ and $(x,z)\in \Oplusgda$ with $x\in\Tgqa$ and $z\in\Oplusxd$
we set
\[\varphi(x,z)=(\varphi(x),z).\]
Note that the groups $G_x$ and $G_{\varphi(x)}$ are the same (only the
isomorphism with $\pi_g$ has changed); therefore
$\Oplusxd=\hat\Omega^+(\varphi(x))$. This action, which is trivial on the
fibres,
descends to actions of $\Gga$ on $\Oplusgda/N_\alpha$ and on $\Cgqa$. The
respective orbit spaces give a family $\Cgq = \Cgqa/\Gga$ over $\Sgqa$ and a
Schottky structure on it. Using the universal property of the family over
$\Tgqa$ (see Theorem~\ref{fein}) and the fact that Schottky structures are
locally induced by Teichm\"uller structures, we find that the Schottky
structure on $\Cgq$ is in fact universal.
\end{proof}
The following diagram collects the relations between the spaces introduced and
used in this section. The horizontal maps are open embeddings, the last two
vertical maps are analytic with discrete fibres; all
other maps in the diagram are quotient maps for the groups indicated (to be
precise, the map from $\Tgqa$ to $\Mgq$ is the restriction of the orbit map
for the action of $\Gamma_g$ on $\Tgq$).\\[1mm]
\newpage
\begin{center}
$\xymatrix@=9ex{\ar@{^{(}->}[rr]\ar[drr]^(.55){/\Gga}\ar[ddrr]^(.6){/\Hga}\ar[dddrr]_{/\Gamma_g}T_g\
      &&\ar[drr]^{/\Gga}\ar@{-}'[dr]\ar[dddrr]^{/\Gamma_g}|(.33)\hole|(.666)\hole\Tgqa\\&&\ar@{^{(}->}[rr]\ar[d]^{/\mbox{\scriptsize
      Out}(\Phi_g)}S_g\ &\ar[dr]^(.35){/\Hga}&\ar[d]^{/\mbox{\scriptsize
      Out}(\Phi_g)}\Sgq\\&&\ar@{^{(}->}[rr]\ar[d]\Sgd\ &&\ar[d]\Sgqd\\&&\ar@{^{(}->}[rr]M_g&&\Mgq\ }$\\[5mm]
\end{center}

\subsection{Teichm\"uller disks in Schottky space}
\label{s-geod}
Let $\iota:\HH\to T_g$ be a Teichm\"uller embedding\index{Teichm\"uller embedding} as in Definition \ref{emb} and $\Delta
= \iota(\HH)$ its image in $T_g$. Let $\Stabid$ be the stabilizer of $\Delta$ in $\Gamma_g$. We have seen in Section \ref{lattice} that $\Stabid$ maps surjectively to the projective Veech group $\Grq$ of $\iota$ (see Definition \ref{Veechgroup}); the kernel of this map is the pointwise stabilizer of $\Delta$.\\

In this section we assume that $\Grq$ is a lattice in
PSL$_2(\RR)$, or equivalently that the image $C_\iota$ of $\Delta$ in $M_g$ is a Teichm\"uller curve\index{Teichm\"uller curve} 
(cf.~Corollary \ref{latticeproperty}). As mentioned in the introduction, Veech showed that $C_\iota$ is not a projective curve and thus cannot be closed in $\Mgq$.
\begin{proposition}
\label{durchschnitt}
Let $\iota:\HH\to T_g$ be a Teichm\"uller embedding such that $\Grq$ is a lattice in
$\mbox{\em PSL}_2(\RR)$. Then there exists a symplectic homomorphism
$\alpha:\pi_g\to\Phi_g$ such that
\[{\emStabid}\cap\Gga\not=\{1\}.\]
\end{proposition}
Since $\Gga$ is torsion free, this implies that the intersection is
infinite. As a consequence, the image of the Teichm\"uller disk $\Delta$ in the Schottky space\index{Teichm\"uller disk!image in Schottky space} $S_g$ is the quotient by an infinite group and in particular 
not isomorphic to a disk.
\begin{proof}
Denote by $\overline{\Delta}$ and $\bar
C_\iota$ the closures of $\Delta$ and $C_\iota$ in $\Tgq$ and $\Mgq$,
respectively. Since $C_\iota$ is not closed, we can 
find a point $z\in \bar C_\iota-C_\iota$; let $x\in\overline \Delta$ be a point above $z$. By
Prop.~\ref{Tgqa}\,b) there is a symplectic homomorphism $\alpha$ such that
$x\in\Tgqa$.\\ 
Let $\bar s_\alpha:\Tgqa\to\Sgq$ be the quotient map for $\Gga$ 
(see~Prop.~\ref{Sgqa} and Prop.~\ref{sgqa=sgq}) and let 
$D(\iota)=s_\alpha(\Delta)$ be the image of $\Delta$
in $S_g$. Then the closure  $\bar D(\iota)$ of $D(\iota)$ in $\Sgq$ contains
 $\bar s_\alpha(x)$, and we have $\bar C_\iota=\bar \mu(\bar
D(\iota))$, cf.~the diagram at the end of Section \ref{s-ext}.\\
By our assumption, $\bar C_\iota$ is Zariski closed in $\Mgq$. Therefore
$\bar\mu^{-1}(\bar C_\iota)$ is an analytic subset of $\Sgq$. $\bar
D(\iota)$ is an irreducible component of $\bar\mu^{-1}(\bar C_\iota)$ and
hence also an analytic subset.\\
Recall, from Corollary~\ref{latticeproperty}, that $\Delta/\Stabid$ is the 
normalization of $C_\iota$. Furthermore, by Prop.~\ref{iquerfortc}, $\overline{\Delta}$ is isomorphic to $\HH\cup\{\mbox{cusps of}\ \Gammaquerm_{\iota}\}$. Therefore $\overline{\Delta}/\Stabid$ is the normalization of $\bar C_\iota$. The restriction of the quotient map $\overline{\Delta}\to\overline{\Delta}/\Stabid$ to the intersection
 $\overline{\Delta}_\alpha=\overline
\Delta\cap\Tgqa$ factors through $\bar s_\alpha$. If the intersection 
$\Stabid\cap\Gga$ was trivial,
this restriction  would be an
isomorphism. But then $\overline{\Delta}_\alpha$ would be isomorphic to an
analytic subset of a complex manifold. This is impossible since
$\overline{\Delta}_\alpha$ contains $x\in \Tgq-T_g$ and hence is not a complex space.
\end{proof}


\end{document}